\documentclass[11pt]{article}
\usepackage[pdftex]{graphicx}
\usepackage[small,bf,sf,textfont={small,sf}]{caption}
\usepackage{vmargin}
\usepackage{subfig}
\usepackage{amsmath}
\usepackage{amsmath}
\usepackage{amssymb}
\usepackage{fancyheadings}
\usepackage{color}
\usepackage[colorlinks=true,
            linkcolor=blue,
            urlcolor=blue,
            citecolor=blue]{hyperref}
\usepackage[round]{natbib}
\usepackage{doi}
\usepackage{booktabs}
\usepackage{float}
\usepackage{lineno}
\usepackage{epstopdf}

\floatstyle{ruled}
\newfloat{floatbox}{thp}{lob}
\floatname{floatbox}{Algorithm}

\newcommand{\trp}{^{\scriptsize \textrm{T}}}

\newcommand{\inv}{^{\scriptsize -1}}

\newcommand{\invsqr}{^{\scriptsize -\textrm{1/2}}}

\newcommand{\x}{\mathbf{x}}
\newcommand{\y}{\mathbf{y}}

\newcommand{\m}{\mathbf{m}}

\newcommand{\e}{\mathbf{e}}

\newcommand{\mtrue}{\mathbf{m}_{\textrm{true}}}

\newcommand{\dsim}{\mathbf{d}}
\newcommand{\dobs}{\mathbf{d}_{\textrm{obs}}}

\newcommand{\g}{\mathbf{g}}
\newcommand{\A}{\mathbf{A}}

\newcommand{\C}{\mathbf{C}}

\newcommand{\Ce}{\mathbf{C}_{\mathbf{e}}}

\newcommand{\DM}{\Delta \mathbf{M}}

\newcommand{\DD}{\Delta \mathbf{D}}

\setpapersize{USletter}
\setmarginsrb{1.5in}{1in}{1in}{0.5in}{12pt}{11mm}{12pt}{30pt}

\begin{document}


\title{Analysis of Geometric Selection of the Data-Error Covariance Inflation for {ES-MDA}}

\author{Alexandre A. Emerick}
\maketitle

\begin{center}
\small \noindent Petrobras Research and Development Center -- CENPES \\
\small Av. Hor\'{a}cio de Macedo 950, Cidade Universit\'{a}ria, Rio de Janeiro, RJ 21941-915, Brazil \\
\small \texttt{emerick@petrobras.com.br}
\end{center}

\section*{Abstract}
\label{Sec:Abstract}

The ensemble smoother with multiple data assimilation (ES-MDA) is becoming a popular assisted history matching method. In its standard form, the method requires the specification of the number of iterations in advance. If the selected number of iterations is not enough, the entire data assimilation must be restarted. Moreover, ES-MDA also requires the selection of data-error covariance inflations. The typical choice is to select constant values. However, previous works indicate that starting with large inflation and gradually decreasing during the data assimilation steps may improve the quality of the final models.

This paper presents an analysis of the use of geometrically decreasing sequences of the data-error covariance inflations. In particular, the paper investigates a recently introduced procedure based on the singular values of a sensitivity matrix computed from the prior ensemble. The paper also introduces a novel procedure to select the inflation factors. The performance of the data assimilation schemes is evaluated in three reservoir history-matching problems with increasing level of complexity. The first problem is a small synthetic case which illustrates that the standard ES-MDA scheme with constant inflation may result in overcorrection of the permeability field and that a geometric sequence can alleviate this problem. The second problem is a recently published benchmark and the third one is a field case with real production data. The data assimilation schemes are compared in terms of a data-mismatch and a model-change norm. The first norm evaluates the ability of the models to reproduce the observed data. The second norm evaluates the amount of changes in the prior model. The results indicate that geometric inflations can generate solutions with good balance between the two norms.

\section{Introduction}
\label{Sec:Intro}

Iterative forms of the ensemble smoother (ES) are rapidly becoming popular assisted history matching (AHM) techniques with a growing number of field applications reported in the literature in the last few years. An evidence of the increasing interest of the industry is the recent appearance of commercial AHM softwares based on iterative smoothers. Among these methods, the ES-MDA \citep{emerick:13b} is a popular choice because of the simplicity of its formulation. In a standard implementation of ES-MDA, it is necessary to specify in advance the number of data assimilations, $N_a$, and the set of data-error covariance inflation factors, $\{\alpha_k\}_{k=1}^{N_a}$. In fact, classifying ES-MDA as an iterative ES is debatable because its termination criterion does not rely on a convergence check \citep{mannseth:15a}. A very common choice for ES-MDA is to set $\alpha_k = N_a$ for $k = 1, \ldots, N_a$ \citep{emerick:13b,emerick:13c,emerick:16a}. Unfortunately, if the selected value for $N_a$ is not sufficient to obtain a reasonable data match, it is necessary to restart the entire data assimilation with a larger value. In practice, however, it is desirable to keep $N_a$ as small as possible to reduce computational cost. Moreover, it has been suggested in the literature \citep{emerick:13b,duc:16a,emerick:16a} that starting with large value for $\alpha_k$ and gradually decreasing this value during the data assimilation steps improves the quality of the final models. The rationality for this procedure is based on the conjecture that at the beginning of the data assimilation the predicted data are far from the observations, so it is beneficial to restrict the changes in the models to avoid overcorrection. This conjecture is supported by previous experiences with the application of Gauss-Newton and Levenberg-Marquardt methods for history matching; see, e.g., \citep[Chap.~8]{oliver:08bk}.

The problem of selecting the inflation factors for ES-MDA has been investigated in previous works. \citet{emerick:16a} presented a ES-MDA algorithm where $\alpha_k$ and $N_a$ are selected according to the progress of the average data-mismatch objective function. This algorithm was tested in a field case requiring $N_a = 6$ data assimilations. However, the final results were not significantly superior to the standard method with $\alpha_k = N_a = 4$. Inspired by the work of \citet{iglesias:13b}, \citet{duc:16a} proposed an ES-MDA scheme that adaptively select $N_a$ and $\alpha_k$ using the regularization condition proposed in \citep{hanke:97a,hanke:10a}. They also investigated an alternative procedure based on limiting the maximum change in the model over a data assimilation step. They tested both schemes in a modified version of the PUNQ-S3 case and concluded that the standard ES-MDA with 8 and 16 data assimilations and constant inflation resulted in over and undershooting of the permeability field (very high and low values of permeability). On the other hand, the proposed schemes obtained smoother permeability realizations. Unfortunately, these adaptive schemes tend to be computationally too demanding when applied to large-scale problems because they often require several repetitions of the analysis (inversion) step, which becomes problematic for big models with large number of data points. Moreover, these schemes usually require several data assimilation steps, often more than ten. Recently, \citet{rafiee:17a} proposed to select the inflation factors in a geometrically decreasing order where the first inflation is based on the regularization condition of \citep{hanke:97a}. For a simple test problem, they showed that the geometric scheme resulted in smoother permeability fields compared to constant inflations.

Here, the objective is to investigate the use of geometric inflations. This paper also introduces a novel scheme to select the inflations. The goal is to develop a robust procedure to improve the data assimilation results without compromising the computational efficiency and ease of implementation of the method.

The remaining of the paper is organized as follows: the next two sections reviews basic concepts of regularization for linear and nonlinear inverse problems. The section after that presents the ES-MDA equations emphasizing that the inflation factors plays the same role of the regularization parameter described in the inverse problems literature. This section also reviews the method proposed in \citep{rafiee:17a} and presents the new scheme for selecting a geometric sequence of $\alpha_k$. Both methods are tested in three reservoir history-matching problems. This first problem is similar to the case used in \citep{rafiee:17a} and it was designed to show overcorrection in the permeability field using constant inflation factors. The second one is a more realistic history-matching exercise proposed in \citep{unisim-i-h:13,avansi:15a}. The third problem is the same field case used in \citep{emerick:16a}. The last section of the paper summarizes the main conclusions of the work.

\section{Regularization for Linear Inverse Problems}
\label{Sec:RLS}

Consider the problem of finding $\x \in \mathbb{R}^{N_x}$ such that

\begin{equation}\label{Eq:LP}
  \A \x = \y,
\end{equation}
where $\y \in \mathbb{R}^{N_y}$ is a vector with given observations and $\A$ is a $N_y \times N_x$ matrix. A unique solution exists only if $\y \in \mathfrak{R}(\A)$ and $dim\left(\mathfrak{N}(\A)\right) = 0$, where $\mathfrak{R}(\A)$ and $\mathfrak{N}(\A)$ denote the range and the null space of $\A$, respectively. From a practical perspective, however, it is necessary to find only an approximate solution of (\ref{Eq:LP}) because $\y$ is typically corrupted with noise. An approximate solution of (\ref{Eq:LP}) can be obtained by solving the least-squares problem

\begin{equation}\label{Eq:LSP}
  \x^\star = \min_{\x \in \mathbb{R}^{N_x}} \left\|\A\x -\y \right\|^2.
\end{equation}
All norms in this paper refer to the $L_2$-norm, i.e., $\|\x\| = \sqrt{\x\trp\x}$. The linear least-squares problem (\ref{Eq:LSP}) is ill-posed in the sense that small errors in $\y$ may cause errors of arbitrary size in $\x$. Instead, one can seek for the solution of the regularized least-squares problem

\begin{equation}\label{Eq:RLSP}
  \x_\alpha = \min_{\x \in \mathbb{R}^{N_x}} \left\{ \left\|\A\x -\y \right\|^2 + \alpha\|\x\|^2 \right\},
\end{equation}
where $\alpha > 0$ is the regularization parameter. Eq.~\ref{Eq:RLSP} is known as Tikhonov regularization. The solution of (\ref{Eq:RLSP}) exists; see, e.g., \citep{kaipio:05bk}, and can be written as

\begin{eqnarray}\label{Eq:xalpha}
  \nonumber \x_\alpha & = & \left(\A\trp\A + \alpha \mathbf{I}\right)\inv \A\trp\y \\
  & = & \A\trp\left(\A\A\trp + \alpha \mathbf{I}\right)\inv\y.
\end{eqnarray}

The choice of the regularization parameter controls the quality of the estimates. There are several methods for choosing $\alpha$ presented in the literature; see, e.g., \citep{hamarik:12a} for a detailed discussion. Perhaps the most well-known methods are the ones based on the Morosov's Discrepancy Principle (MDP) \citep{morozov:84a}, which states that it is not reasonable to expect the approximate solution to yield an error smaller than the noise level in the data vector. The MDP can be formulated as selecting $\alpha$ such that

\begin{equation}\label{Eq:MDP}
  \left\|\A\x_\alpha - \y \right\| = \tau \eta.
\end{equation}
In the above equation, $\eta > 0$ is the noise level in $\y$ and $\tau \geq 1$ is an extra parameter to avoid underregularization \citep{kaipio:05bk}. Let $\y$ be corrupted with an unknown additive noise vector $\e$, i.e.,

\begin{equation}\label{Eq:y}
  \y = \y_{\textrm{true}} + \e,
\end{equation}
where $\y_{\textrm{true}}$ is the noiseless observation vector. Assuming that $\e$ is a random vector, then it is convenient to define the square of the noise level as

\begin{equation}\label{Eq:eta2}
  \eta^2 = \textrm{E} \left[ \|\e\|^2 \right].
\end{equation}
Let the noise vector be a sample from a multivariate Gaussian distribution. Without loss of generality, assume that $\e \sim \mathcal{N}\left(\mathbf{0}, \sigma_y^2 \mathbf{I}\right)$ with $\sigma_y$ denoting the standard deviation. Then $\|\e\|^2$ follows as chi-squared distribution with $N_y$ degrees of freedom. Therefore, the noise level becomes

\begin{equation}\label{Eq:eta2}
  \eta = \sqrt{N_y}\sigma_y.
\end{equation}

\section{Regularization for Nonlinear Inverse Problems}
\label{Sec:RNLS}

Consider the problem of computing an approximate solution for the nonlinear least-squares problem

\begin{equation}\label{Eq:NLSP}
 \x^\star =  \min_{\x \in \mathbb{R}^{N_x}} \left\|\mathbf{f}(\x) -\y \right\|^2
\end{equation}
with $\mathbf{f}: \mathbb{R}^{N_x} \rightarrow \mathbb{R}^{N_y}$ being a nonlinear vector function. In order to solve (\ref{Eq:NLSP}), it is necessary to derive an iterative scheme. Consider the first-order Taylor expansion around the current guess $\x^k$

\begin{equation}\label{Eq:Taylor}
  \mathbf{f}(\x) \approx \mathbf{f}(\x^k)  + \A_k \left(\x - \x^k\right)
\end{equation}
where
\begin{equation}\label{Eq:Ak}
  \A_k =\left[
          \begin{array}{ccc}
            \frac{\partial f_1(\x)}{\partial x_1} & \cdots & \frac{\partial f_1(\x)}{\partial x_{N_x}} \\
            \vdots & \ddots & \vdots \\
            \frac{\partial f_{N_y}(\x)}{\partial x_1} & \cdots & \frac{\partial f_{N_y}(\x)}{\partial x_{N_x}}
          \end{array}
        \right]
\end{equation}
is the sensitivity matrix. The linearized version of (\ref{Eq:NLSP}) becomes

\begin{equation}\label{Eq:NLSP2}
  \x^{k+1} = \min_{\x \in \mathbb{R}^{N_x}} \left\| \mathbf{f}(\x^k)  + \A_k \left(\x - \x^k\right) - \y \right\|^2.
\end{equation}
Let $\widehat{\x}^k =  \x - \x^k$ and $\widehat{\y}^k = \y - \mathbf{f}(\x^k)$, then (\ref{Eq:NLSP2}) becomes

\begin{equation}\label{Eq:NLSP3}
  \widehat{\x}^{k+1} = \min_{\widehat{\x} \in \mathbb{R}^{N_x}} \left\| \A_k\widehat{\x}^k -\widehat{\y}^k \right\|^2,
\end{equation}
which has the same form of the linear least-squares problem (\ref{Eq:LSP}). This problem is ill-posed and instead of (\ref{Eq:NLSP3}), consider the following regularized problem

\begin{equation}\label{Eq:RNLSP}
  \widehat{\x}^{k+1}_\alpha = \min_{\widehat{\x} \in \mathbb{R}^{N_x}} \left\{ \left\|\A_k\widehat{\x}^k -\widehat{\y}^k \right\|^2 + \alpha_k\left \|\widehat{\x}^k \right\|^2 \right\},
\end{equation}
which has the same solution presented in Eq.~\ref{Eq:xalpha}. The resulting iterative procedure corresponds to the Levenberg-Marquardt (LM) algorithm for nonlinear least squares; see, e.g. \citep{nocedal:06}.

\citet{hanke:97a} used the MDP to propose a regularized LM scheme (see also \citep{hanke:10a}). In his scheme, the regularization parameter $\alpha_k$ should be selected such that

\begin{equation}\label{Eq:Hanke}
  \rho^2 \left\| \widehat{\y}^k \right\|^2 \leq \alpha_k^2 \left\| \left(\A_k\A_k\trp + \alpha_k \mathbf{I} \right)\inv \widehat{\y}^k \right\|^2
\end{equation}
for some $\rho \in (0,1)$. This condition is the base for the works presented by \citep{iglesias:13b,luo:15b,duc:16a} and it is the starting point for the procedure proposed in \citep{rafiee:17a}, which is discussed in the next sections.

\section{ES-MDA}
\label{Sec:ES-MDA}

The ES-MDA analysis equation to update a vector of model parameter, $\m \in \mathbb{R}^{N_m}$, is

\begin{equation}\label{Eq:ESMDA}
  \m^{k+1}_j = \m^{k}_j + \C_{\m\dsim}^k \left(\C_{\dsim\dsim}^k + \alpha_{k+1} \Ce \right)\inv \left(\dobs + \e^k_j - \g(\m^k_j)\right),~~\text{for}~j=1,\ldots, N_e,
\end{equation}
where $N_e$ is the size of the ensemble. $\dobs \in \mathbb{R}^{N_d}$ is the vector containing the observation and $\g(\m^k_j) \in \mathbb{R}^{N_d}$ is the vector of data predicted by the $j$th model realization at the $k$th data assimilation step of ES-MDA. $\e^k_j \in \mathbb{R}^{N_d}$ is a random vector drawn from $\mathcal{N}(\mathbf{0}, \alpha_{k+1}\Ce)$, where $\Ce$ is the data-error covariance matrix. Here, $N_d$ is used to denote the number of data points. The covariance matrices $\C_{\m\dsim}^k$ and $\C_{\dsim\dsim}^k$ are computed as

\begin{equation}\label{Eq:Cmd}
   \C_{\m\dsim}^k = \DM^k \left(\DD^k\right)\trp
\end{equation}
and
\begin{equation}\label{Eq:Cdd}
   \C_{\dsim\dsim}^k =  \DD^k \left(\DD^k\right)\trp
\end{equation}
where
\begin{equation}\label{Eq:DM}
   \DM^k = \frac{1}{\sqrt{N_e - 1}} \left[ \begin{array}{ccc} \m_1^k - \overline{\m^k} & \cdots & \m_{N_e}^k - \overline{\m^k}  \\ \end{array}  \right]
\end{equation}
and
\begin{equation}\label{Eq:DM}
   \DD^k = \frac{1}{\sqrt{N_e - 1}} \left[ \begin{array}{ccc} \g(\m_1^k) - \overline{\g(\m^k)} & \cdots & \g(\m_{N_e}^k) - \overline{\g(\m^k)} \\ \end{array}  \right].
\end{equation}

In ES-MDA, $\alpha_k > 1$ is referred to as data-error covariance inflation factor. In the standard implementation of this method, Eq.~\ref{Eq:ESMDA} is applied a pre-defined number of times, $N_a$, and the set $\{\alpha_k \}_{k=1}^{N_a}$ must satisfy the condition

\begin{equation}\label{Eq:ESMDACond}
  \sum_{k=1}^{N_a} \frac{1}{\alpha_k} = 1.
\end{equation}
This condition was originally derived by \citet{emerick:13b} using linear algebra to ensure that ES-MDA is consistent with the standard ES for the linear-Gaussian case. However, the same condition can be derived directly from Bayes' rule as discussed in \citep{stordal:15a,emerick:16a}.

\subsection{Geometric Selection of ES-MDA Inflation}
\label{Sec:Geo-ES-MDA}

Here, the objective is to derive a scheme to select $\alpha_k$ for $k = 1, \ldots, N_a$ to be used in the ES-MDA method. Specifically, the goal is to select a geometrically decreasing sequence $\left\{\alpha_k\right\}_{k=1}^{N_a}$ such that

\begin{equation}\label{Eq:alphak}
  \alpha_{k+1} = \gamma \alpha_k = \gamma^{k} \alpha_1
\end{equation}
for $\gamma \in (0, 1]$. For a fixed $\alpha_1$, which can be selected with one of the procedures discussed in the next sections, the computation of $\gamma$ is straightforward. Recall that

\begin{equation}\label{Eq:alpha1}
  \sum_{k=1}^{N_a} \frac{1}{\alpha_k} = \sum_{k=1}^{N_a} \frac{1}{\gamma^{k-1} \alpha_1} = 1.
\end{equation}
Hence, defining the function

\begin{equation}\label{Eq:fgamma}
  f_1(\gamma) = \sum_{k=1}^{N_a} \frac{1}{\gamma^{k-1} \alpha_1} - 1,
\end{equation}
it is easy to find $\gamma \in (0, 1]$ such that $f_1(\gamma) = 0$ using, for example, the bisection method \citep{press:07}.

Before introducing the procedures to compute $\alpha_1$, it is convenient first to consider the ES-MDA equation to update the ensemble mean, $\overline{\m}$, i.e.,

\begin{equation}\label{Eq:MDAmean}
  \overline{\m^{k+1}} = \overline{\m^{k}} + \C_{\m\dsim}^k \left(\C_{\dsim\dsim}^k + \alpha_{k+1} \Ce \right)\inv \left(\dobs - \overline{\g(\m^k)}\right).
\end{equation}
Eq.~\ref{Eq:MDAmean} uses the fact that $\overline{\e^k}=\mathbf{0}$. Using (\ref{Eq:Cmd}) and (\ref{Eq:Cdd}) in (\ref{Eq:MDAmean}) and setting $k = 0$

\begin{eqnarray}\label{Eq:MDAmean1}
  \nonumber \overline{\m^{1}} & = & \overline{\m^{0}} + \DM^0 \left(\DD^0\right)\trp \left(\DD^0 \left(\DD^0\right)\trp + \alpha_{1} \Ce \right)\inv \left(\dobs - \overline{\g(\m^k)}\right)\\
  \nonumber & = & \overline{\m^{0}} + \DM^0 \left(\Ce\invsqr \DD^0\right)\trp \left( \left(\Ce\invsqr \DD^0\right) \left(\Ce\invsqr \DD^0\right)\trp + \alpha_{1} \mathbf{I} \right)\inv \\
    & ~ & \times \Ce\invsqr \left(\dobs - \overline{\g(\m^0)}\right)
\end{eqnarray}
or
\begin{eqnarray}\label{Eq:MDAmean2}
 \nonumber \left(\DM^0\right)^\dag  \left(\overline{\m^{1}} - \overline{\m^{0}} \right) & = &  \left(\Ce\invsqr \DD^0\right)\trp \left( \left(\Ce\invsqr \DD^0\right) \left(\Ce\invsqr \DD^0\right)\trp + \alpha_{1} \mathbf{I} \right)\inv \\
 & ~ & \times \Ce\invsqr \left(\dobs - \overline{\g(\m^0)}\right)
\end{eqnarray}
where $\m^0$ denotes a prior realization and $\left(\DM^0\right)^\dag$ is the pseudo-inverse of $\DM^0$. Note that there will be no need to compute $\left(\DM^0\right)^\dag$ as the objective is only to derive a scheme to compute $\alpha_1$. Now call

\begin{equation}\label{Eq:MDAx}
  \x = \left(\DM^0\right)^\dag  \left(\overline{\m^{1}} - \overline{\m^{0}} \right),
\end{equation}

\begin{equation}\label{Eq:MDAy}
  \y = \Ce\invsqr \left(\dobs - \overline{\g(\m^0)}\right)
\end{equation}
and
\begin{equation}\label{Eq:MDAA}
  \A =  \Ce\invsqr \DD^0,
\end{equation}
so that Eq.~\ref{Eq:MDAmean2} can be written as

\begin{equation}\label{Eq:MDAx1}
  \x = \A\trp\left(\A\A\trp + \alpha_1 \mathbf{I}\right)\inv\y.
\end{equation}
Here, the dimension of the vectors $\x$ and $\y$ are $N_x = N_e$ and $N_y = N_d$, respectively. Eq.~\ref{Eq:MDAx1} assumes the same form of (\ref{Eq:xalpha}), which is the solution of the regularized linear least-squares problem (\ref{Eq:RLSP}). Clearly the inflation factor plays the role of the Tikhonov regularization parameter. It is worth mentioning that the matrix $\A$ defined in Eq.~\ref{Eq:MDAA} is scaled by the square-root of the data-error covariance matrix. This has numerical benefits because the matrix $\DD^0$ may be constructed with data of different orders of magnitude. Hence, scaling $\A$ avoids losing relevant information from data when truncating small singular values. Moreover, $\A$ corresponds to the ensemble equivalent of the dimensionless sensitivity matrix \citep{zhang:02b} as discussed in the Appendix~A of \citep{emerick:12a}. \citet{tavakoli:10a} showed that for a linear-Gaussian problem the singular values of the dimensionless sensitivity govern the reduction in uncertainty and that the smallest singular values have negligible influence on the reduction of uncertainty. In this sense, defining $\A$ as in Eq.~\ref{Eq:MDAA} is optimal.

It is important to note that the data-error covariance matrix, $\Ce$, plays a crucial role in the data assimilation process. The construction of this matrix, however, is not very well discussed in the petroleum literature. Here, it is assumed that the measurement errors in data are independent, such that $\Ce$ is a diagonal matrix. This is a very common assumption and it is present in the large majority of the history-matching publications. Perhaps the main exceptions are the papers that discuss assimilation of seismic data; see, e.g., \citep{emerick:16a,luo:17b}, and the work from \citet{evensen:18b} which discusses the use of rate an cumulative production data. However, even the selection of the variance (diagonal elements of $\Ce$) is not very well documented. A common procedure for production data is to estimate the data variance using a smoothed version of the time series \citep{zhao:06,emerick:11c}. However, the noise in the measured data are not the only source of errors. There are also the ``model errors,'' which includes errors because of coarse discretizations, simplified physics and insufficient parameterizations, just to mention a few. Model errors are difficult to be correctly accounted for and traditionally they are neglected in the majority of the history-matching applications. However, if we make the crude assumption that model errors are additive and follow the Gaussian distribution $\mathbf{q} \sim \mathcal{N}(\mathbf{0},\C_\mathbf{q})$, then it is straightforward to include their effect in the data assimilation by simply replacing $\Ce$ by $\Ce + \C_\mathbf{q}$. In fact, some recent works \citep{sun:17a,oliver:17a} show that selecting a larger variance for the data noise tends to partially compensate for model errors and improve the data assimilation results, especially in terms of uncertainty quantification. For this reason and based on previous experiences with field history-matching problems, in the last two test problems presented in this paper the selected data-error variances are overestimated compared to the variance that would be estimated based solely on the production time series.

The next sections present two procedures to compute $\alpha_1$. The first procedure was originally presented by \citet{rafiee:17a} and the second is proposed in this paper.

\subsubsection{Procedure Proposed by Rafiee and Reynolds}

\citet{rafiee:17a} used the Hanke's regularization condition (\ref{Eq:Hanke}) as a starting point to define a criterion to select $\alpha_1$. First, they wrote Eq.~\ref{Eq:Hanke} as

\begin{equation}\label{Eq:Reynolds1}
  \rho^2 \leq \alpha_1^2 \frac{\left\| \C\inv\y \right\|^2}{\left\| \y \right\|^2}
\end{equation}
where $\C = \left( \A\A\trp + \alpha_1 \mathbf{I}\right)$. Then they noted that

\begin{equation}\label{Eq:Reynolds2}
   \rho^2 \leq \alpha_1^2 \max_i \lambda_i^2 = \alpha_1^2 \max_i \frac{1}{\left(\sigma_i^2 + \alpha_1 \right)^2}
\end{equation}
with $\lambda_i$ denoting the $i$th eigenvalue of $\C$ and $\sigma_i$ the $i$th singular value of $\A$. Instead of selecting $\alpha_1$ that satisfies the condition (\ref{Eq:Reynolds2}), they use

\begin{equation}\label{Eq:Reynolds3}
  \rho^2 \leq \alpha_1^2 \frac{1}{\left(\overline{\sigma}^2 + \alpha_1 \right)^2}.
\end{equation}
The equality occurs for

\begin{equation}\label{Eq:Reynolds4}
  \alpha_1 = \frac{\rho}{1 - \rho} \overline{\sigma}^2,
\end{equation}
where
\begin{equation}\label{Eq:Reynolds5}
  \overline{\sigma} = \frac{1}{N_r}\sum_{i=1}^{N_r} \sigma_i,
\end{equation}
with $N_r = \min\left\{N_y, N_x\right\}$ denoting the number of non-zero singular values of $\A$.

\citet{rafiee:17a} proposed to use $\rho = 0.5$, which leads to $\alpha_1 = \overline{\sigma}^2$. The resulting procedure is summarized in the Algorithm~\ref{Algo:GEO1}. Here, this method is referred to as GEO1.

\begin{floatbox}{}{}
\begin{enumerate}
  \item Specify $N_a$.
  \item Compute
      \begin{equation}
         \nonumber \alpha_1 = \max \left\{ \overline{\sigma}^2, N_a \right\},
      \end{equation}
  where
  \begin{equation}
    \nonumber \overline{\sigma} = \frac{1}{N_r}\sum_{i=1}^{N_r} \sigma_i.
  \end{equation}
  \item Compute $\gamma \in (0, 1]$ by solving
      \begin{equation}\label{Eq:fgamma}
        \nonumber \sum_{k=1}^{N_a} \frac{1}{\gamma^{k-1} \alpha_1} - 1 = 0
      \end{equation}
  using the bisection method.
  \item Apply ES-MDA with $\alpha_{k+1} = \gamma^k \alpha_1$.
\end{enumerate}
\caption{ES-MDA-GEO1}\label{Algo:GEO1}
\end{floatbox}

\subsubsection{Procedure Proposed in this Work}

The procedure to define the inflation factors proposed here is inspired by the ideas in \citep{rafiee:17a} in the sense that it is also based on a geometrically decreasing sequence. However, there are two main differences between the two procedures.

The first difference is that instead of computing $\alpha_1$ independently from the value of $N_a$, the new proposal computes $\alpha_1$ based on a specified inflation factor for the last step of the data assimilation, $\alpha_{N_a}$. Then the value of $N_a$ is checked against the MDP. The idea is to ensure that the geometric scheme $\alpha_{k+1} = \gamma^k\alpha_1$ will not result in a sequence that decays too fast. In fact, during the computational experiments that conducted to develop this procedure, it has been noticed that the procedure proposed in \citep{rafiee:17a} tends to generate good results as long as an appropriate value for $N_a$ is selected in advance. This is particularly evident for a small number of data assimilations, say $N_a = 4$, and a large initial inflation, say $\alpha_1 \geq 10^4$. This situation generates a small value for $\gamma$ and consequently the last inflation is often very close to one. This fact is illustrated in the test cases presented in the next sections. The idea to circumvent this problem is to specify the value of the inflation at the last data assimilation and compute $\alpha_1$ using

\begin{equation}\label{Eq:Alpha1New}
 \alpha_1 = \gamma^{1 - N_a}\alpha_{N_a}.
\end{equation}
Now, define

\begin{equation}
  f_2(\gamma) = \sum_{k=1}^{N_a} \frac{1}{\gamma^{k-N_a} \alpha_{N_a}} - 1
\end{equation}
and compute $\gamma \in (0, 1]$ by solving $f_2(\gamma) = 0$ using the bisection method.

The second difference between \citet{rafiee:17a} and the new proposal is that \citet{rafiee:17a} used the Hanke's regularization condition (\ref{Eq:Hanke}) to derive a procedure to compute $\alpha_1$. In the new proposal the value of $\alpha_1$ computed using Eq.~\ref{Eq:Alpha1New} is tested against the MDP. If $\alpha_1$ does not satisfy the MDP, then the number of data assimilation steps is increased until a suitable value of $\alpha_1$ is obtained.

In order to check the MDP, it suffices to define a discrepancy function as

\begin{equation}\label{Eq:f1}
  h(\alpha) = \left\|\A\x_\alpha - \y \right\|^2 - (\tau\eta)^2,
\end{equation}
and compute $\alpha^\star > 0$ such that $h(\alpha^\star) = 0$ by solving a root-finding problem.

Consider the singular value decomposition of $\A$

\begin{equation}\label{Eq:SVDA}
  \A = \mathbf{U} \mathbf{\Sigma} \mathbf{V}\trp,
\end{equation}
where $\mathbf{U} \in \mathbb{R}^{N_y \times N_y}$ and $\mathbf{V} \in \mathbb{R}^{N_x \times N_x}$ are orthogonal matrices and $\mathbf{\Sigma} \in \mathbb{R}^{N_y \times N_x}$ is a matrix with diagonal element given by

\begin{equation}\label{Eq:singular}
  \sigma_1 \geq \sigma_2 \geq \ldots \geq \sigma_{N_r} > \sigma_{N_r + 1} = \ldots = \sigma_{\min\{N_y, N_x\}} = 0.
\end{equation}
Using (\ref{Eq:SVDA}) and the properties $\mathbf{U}\mathbf{U}\trp = \mathbf{I}$ and $\mathbf{V}\mathbf{V}\trp = \mathbf{I}$  in (\ref{Eq:f1}) results in

\begin{eqnarray}\label{Eq:f2}
  \nonumber h(\alpha) & = & \left\|\A\x_\alpha - \y \right\|^2 - (\tau\eta)^2 \\
  \nonumber & = & \left\| \mathbf{U}\trp \left(\A \mathbf{V}\mathbf{V}\trp \x_\alpha - \y\right) \right\|^2 - (\tau\eta)^2 \\
  \nonumber & = & \left\| \mathbf{\Sigma}\mathbf{V}\trp \x_\alpha - \mathbf{U}\trp\y \right\|^2 - (\tau\eta)^2 \\
  & = & \sum_{i=1}^{N_r} \left( \sigma_i \mathbf{v}_i\trp \x_\alpha - \mathbf{u}_i\trp\y  \right)^2 + \sum_{i=N_r + 1}^{N_y} \left( \mathbf{u}_i\trp\y  \right)^2 - (\tau\eta)^2.
\end{eqnarray}

Recall that $\x_\alpha$ is given by

\begin{eqnarray}\label{Eq:x1}
  \nonumber \x_\alpha & = & \left(\A\trp\A + \alpha \mathbf{I}\right)\inv \A\trp\y \\
  \nonumber & = & \left( \mathbf{V}\mathbf{\Sigma}\trp \mathbf{\Sigma}\mathbf{V}\trp + \alpha \mathbf{I} \right)\inv \mathbf{V}\mathbf{\Sigma}\trp \mathbf{U}\trp\y \\
  \nonumber & = & \mathbf{V}\left( \mathbf{\Sigma}\trp \mathbf{\Sigma} + \alpha \mathbf{I} \right)\mathbf{\Sigma}\trp \mathbf{U}\trp\y \\
  & = & \sum_{j=1}^{N_r} \left( \frac{\sigma_j}{\sigma_j^2 + \alpha} \mathbf{u}_j\trp\y  \right)\mathbf{v}_j + \x_0,
\end{eqnarray}
where $\x_0 \in \mathfrak{N}(\A)$ can be chosen arbitrarily as long as $\A\x_0 = \mathbf{0}$. In particular, let $\x_0 = \mathbf{0}$ in the following. Using (\ref{Eq:x1}) in (\ref{Eq:f2}) results in

\begin{eqnarray}\label{Eq:f3}
  \nonumber h(\alpha) & = & \sum_{i=1}^{N_r} \left[ \sigma_i \mathbf{v}_i\trp \left(\sum_{j=1}^{N_r} \left( \frac{\sigma_j}{\sigma_j^2 + \alpha} \mathbf{u}_j\trp\y  \right)\mathbf{v}_j \right) - \mathbf{u}_i\trp\y  \right]^2 + \sum_{i=N_r + 1}^{N_y} \left( \mathbf{u}_i\trp\y  \right)^2 - (\tau\eta)^2 \\
  \nonumber & = & \sum_{i=1}^{N_r} \left[ \frac{\sigma_i^2}{\sigma_i^2 + \alpha} \mathbf{u}_i\trp\y - \mathbf{u}_i\trp\y  \right]^2 + \sum_{i=N_r + 1}^{N_y} \left( \mathbf{u}_i\trp\y  \right)^2 - (\tau\eta)^2 \\
  & = & \sum_{i=1}^{N_r} \left( \frac{\alpha}{\sigma_i^2 + \alpha} \mathbf{u}_i\trp\y \right)^2 + \sum_{i=N_r + 1}^{N_y} \left( \mathbf{u}_i\trp\y  \right)^2 - (\tau\eta)^2.
\end{eqnarray}

The value of $\alpha^\star$ such that $h(\alpha^\star) = 0$ can be computed using the Newton-Raphson method \citep{press:07}, in which case it is necessary to compute the derivative
\begin{equation}\label{Eq:df}
  \frac{\textrm{d}h(\alpha)}{\textrm{d}\alpha} = \sum_{i=1}^{N_r} \frac{2\alpha\sigma_i^2}{\left(\sigma_i^2 + \alpha\right)^3}\left(\mathbf{u}_i\trp\y\right)^2.
\end{equation}

Note that $\frac{\textrm{d}h(\alpha)}{\textrm{d}\alpha} > 0$ for any $\alpha > 0$. Therefore, $h(\alpha)$ is a strictly increasing function of $\alpha$. Moreover, note that

\begin{equation}\label{Eq:lim0}
  \lim_{\alpha \rightarrow 0^+} h(\alpha) = \sum_{i=N_r + 1}^{N_y} \left( \mathbf{u}_i\trp\y  \right)^2 - (\tau\eta)^2 = \left\| \mathbf{U}_0\trp \y \right\|^2 - (\tau\eta)^2
\end{equation}
and
\begin{equation}\label{Eq:limInf}
  \lim_{\alpha \rightarrow \infty} h(\alpha) = \sum_{i=1}^{N_y} \left( \mathbf{u}_i\trp\y  \right)^2 - (\tau\eta)^2 = \left\| \y \right\|^2 - (\tau\eta)^2.
\end{equation}
Therefore, $h(\alpha) = 0$ has a unique solution as long as the noise level satisfies the following condition

\begin{equation}\label{Eq:Cond1}
  \left\| \mathbf{U}_0\trp \y \right\|^2 \leq (\tau\eta)^2 \leq \left\| \y \right\|^2.
\end{equation}

The matrix $\mathbf{U}_0$ forms a basis for $\mathfrak{N}\left(\A\trp\right)$. The left hand side of the inequality (\ref{Eq:Cond1}) can be interpreted as that one would expect that the components of $\y$ orthogonal to $\mathfrak{R}(\A\trp)$ are due to noise \citep{kaipio:05bk}. However, in the specific case of ensemble data assimilation, $\A$ is given by Eq.~\ref{Eq:MDAA}; i.e., the components of $\A$ are computed from an ensemble of size $N_e$. Therefore, $rank(\A) \leq \left\{N_d, N_e - 1\right\}$. If $N_y = N_d \gg N_e$, which is often the case for history matching with smoothers, then $dim\left(\mathfrak{N}\left(\A\trp\right)\right) \gg N_e$. Hence, it is possible to have the cases where $\left\| \mathbf{U}_0\trp \y \right\|^2 > (\tau\eta)^2$. This effectively means that the MDP would not be violated for any $\alpha > 0$. To avoid this case, a safe simplification was added to the proposed method by simply dropping the term $\sum_{i=N_r + 1}^{N_y} \left( \mathbf{u}_i\trp\y  \right)^2$ in Eq.~\ref{Eq:f3}. This simplification is based on the fact that in practice we always use localization which expands the space for reconstruction of $\x_\alpha$ (the localized version of $\A\A\trp$ is full rank). The right hand side of the inequality (\ref{Eq:Cond1}) can be interpreted as that one would expect the noise level not to exceed the signal \citep{kaipio:05bk}.

Note that the MDP was proposed as a procedure to regularize linear inverse problems. The condition proposed by \citet{hanke:97a}, on the other hand, was developed to generate a sequence of $\alpha_k$ that stabilizes the LM algorithm in nonlinear least-squares problems. However, as discussed in the introduction section of this paper, the use of the sequence of $\alpha_k$ from Eq.~\ref{Eq:Hanke} with ES-MDA generates a method that takes too many iterations for practical applications. This fact motivated \citet{rafiee:17a} to use Eq.~\ref{Eq:Hanke} to define only $\alpha_1$. However, one might argue that there is no reason to enforce Eq.~\ref{Eq:Hanke} to compute $\alpha_1$ if a different sequence $\{\alpha_k\}_{k=1}^{N_a}$ is to be used. In the new proposal, on the other hand, the idea is to use the MDP to check if the selected value of $\alpha_1$ is enough to ensure a stable linear inversion in the first data assimilation step. But, similarly to \citet{rafiee:17a}, a geometric decreasing sequence is imposed to keep the computational cost affordable for large-scale problems.

Algorithm~\ref{Algo:GEO2} summarizes the proposed method, which is referred to as GEO2. Besides $\alpha_{N_a}$ the algorithm requires the specification of a maximum inflation, $\alpha_{\max}$, which is used as a safeguard to ensure that the method will not take too many data assimilations. The fourth step of Algorithm~\ref{Algo:GEO2} requires to solve the discrepancy function which is done using the Algorithm~\ref{Algo:NR}. The cases presented in this paper use the following values for the parameters of Algorithms~\ref{Algo:GEO2} and \ref{Algo:NR}: $\alpha_{N_a} = 1.5$, $\alpha_{\max} = 10^5$, $N_a = 4$, $\tau = 1$ and $L_{\max} = 100$.

\begin{floatbox}{}{}
\begin{enumerate}
  \item Specify $\alpha_{N_a}$, $\alpha_{\max}$ and $N_a$.
  \item Compute $\gamma \in (0, 1]$ by solving
      \begin{equation}
        \nonumber \sum_{k=1}^{N_a} \frac{1}{\gamma^{k-N_a} \alpha_{N_a}} - 1 = 0
      \end{equation}
  using the bisection method.
  \item Compute
      \begin{equation}
        \nonumber \alpha_1 = \gamma^{1-N_a} \alpha_{N_a}.
      \end{equation}
  \item Compute $\alpha^\star \in \left[N_a, \alpha_{\max}\right]$ by solving $h(\alpha) = 0$ using the Newton-Raphson method (Algorithm~\ref{Algo:NR}).
  \item If $\alpha_1 < \alpha^\star$ set $N_a = N_a + 1$ and return to step 2.\\
  Else apply ES-MDA with $\alpha_{k+1} = \gamma^k \alpha_1$.
\end{enumerate}
\caption{ES-MDA-GEO2}\label{Algo:GEO2}
\end{floatbox}

\begin{floatbox}{}{}
\begin{enumerate}
  \item Specify $\alpha_{\min} = N_a$, $\alpha_{\max}$, $\tau$ and $L_{\max}$.
  \item Compute
      \begin{equation}
        \nonumber h\left(\alpha_{\min}\right) = \sum_{i=1}^{N_r} \left( \frac{\alpha_{\min}}{\sigma_i^2 + \alpha_{\min}} \mathbf{u}_i\trp\y \right)^2 - \tau^2N_d
      \end{equation}
  \item If $h\left(\alpha_{\min}\right) \geq 0$ set $\alpha^\star = \alpha_{\min}$ and exit.
  \item Compute
  \begin{equation}
    \nonumber \alpha_{(1)} = \frac{\alpha_{\min} + \alpha_{\max}}{2}.
  \end{equation}
  \item For $\ell = 1$ to $L_{\max}$ do
  \begin{enumerate}
    \item Compute
      \begin{equation}
        \nonumber h\left(\alpha_{(\ell)}\right) = \sum_{i=1}^{N_r} \left( \frac{\alpha_{(\ell)}}{\sigma_i^2 + \alpha_{(\ell)}} \mathbf{u}_i\trp\y \right)^2 - \tau^2N_d
      \end{equation}
    and
      \begin{equation}
        \nonumber h^\prime\left(\alpha_{(\ell)}\right) = \sum_{i=1}^{N_r} \frac{2\alpha_{(\ell)}\sigma_i^2}{\left(\sigma_i^2 + \alpha_{(\ell)} \right)^3}\left(\mathbf{u}_i\trp\y\right)^2.
      \end{equation}
    \item Compute
    \begin{equation}
      \nonumber \alpha_{(\ell+1)} = \alpha_{(\ell)} - \frac{h\left(\alpha_{(\ell)}\right)}{h^\prime\left(\alpha_{(\ell)}\right)}.
    \end{equation}
    \item If $\alpha_{(\ell+1)} \geq \alpha_{\max}$ set $\alpha^\star = \alpha_{\max}$ and exit. \\
    Else if $\left|\alpha_{(\ell+1)} - \alpha_{(\ell)} \right| < 10^{-3}$ set $\alpha^\star = \alpha_{(\ell+1)}$ and exit. \\
  \end{enumerate}
  End for.
  \item Set $\alpha^\star = \alpha_{(\ell)}$.
\end{enumerate}
\caption{Newton-Rapshon to compute $\alpha^\star$}\label{Algo:NR}
\end{floatbox}

\section{Test Cases}
\label{Sec:TestCases}

\subsection{2D Model}
\label{Sec:2DModel}

\begin{figure}
	\centering
	\includegraphics[width=0.4\linewidth]{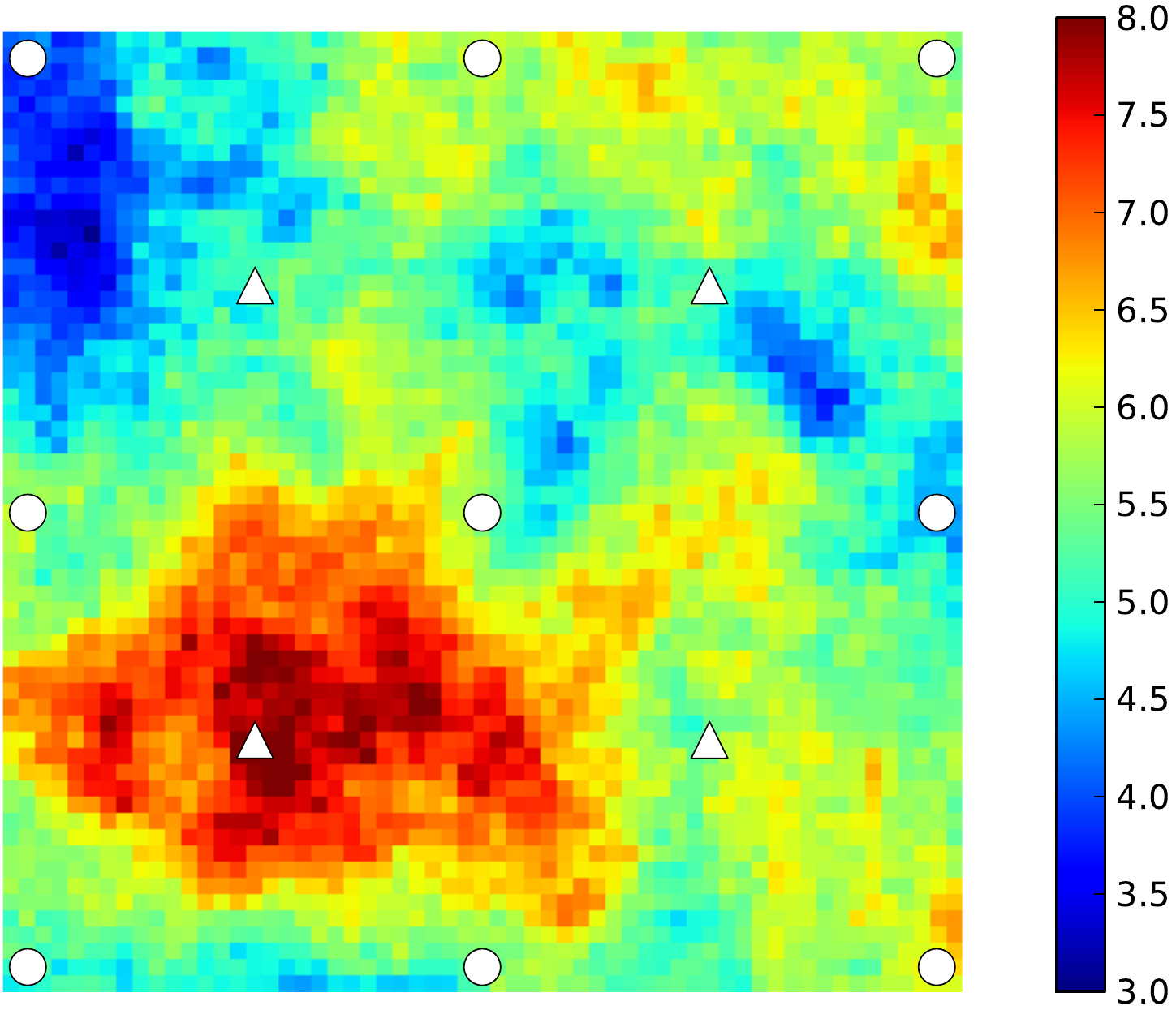}
	\caption{True log-permeability (2D model). Circles correspond to oil producing and triangles to water injection wells.}
	\label{Fig:2DModelTrue}
\end{figure}

The first test case is a small synthetic history-matching problem. The reservoir corresponds to a 2D model discretized into a $59 \times 59$ uniform grid. The only uncertain property is the permeability field. The reference (true) permeability field (Fig.~\ref{Fig:2DModelTrue}) was generated using a spherical covariance function with isotropic correlation length of 40 gridblocks. The value of the prior mean of the natural logarithm of permeability (log-permeability) is constant and equal to 5.5 and the prior standard deviation is equal to one. No hard data were used at well locations to make the problem more challenging. The same setup was used to generate 200 prior realizations for history matching. The reservoir produces with four five-spots (Fig.~\ref{Fig:2DModelTrue}). All wells operate at constant bottom-hole pressure of 20000~kPa for producers and 30000~kPa for water injectors. The observed data correspond to values of oil and water rates every 90~days for a period of five years. Random noise with standard deviation corresponding to 3\% of the actual data were added to form the synthetic measurements. A small value for the measurement errors was intentionally selected because it tends to make the problem more prone to overcorrection. Also note that in this simple problem there are no model errors, so a diagonal matrix $\Ce$ is used with variance equal to the square of the standard deviation of the random noise used to generated the synthetic measurements.

\begin{figure}
	\centering
	\includegraphics[width=0.5\linewidth]{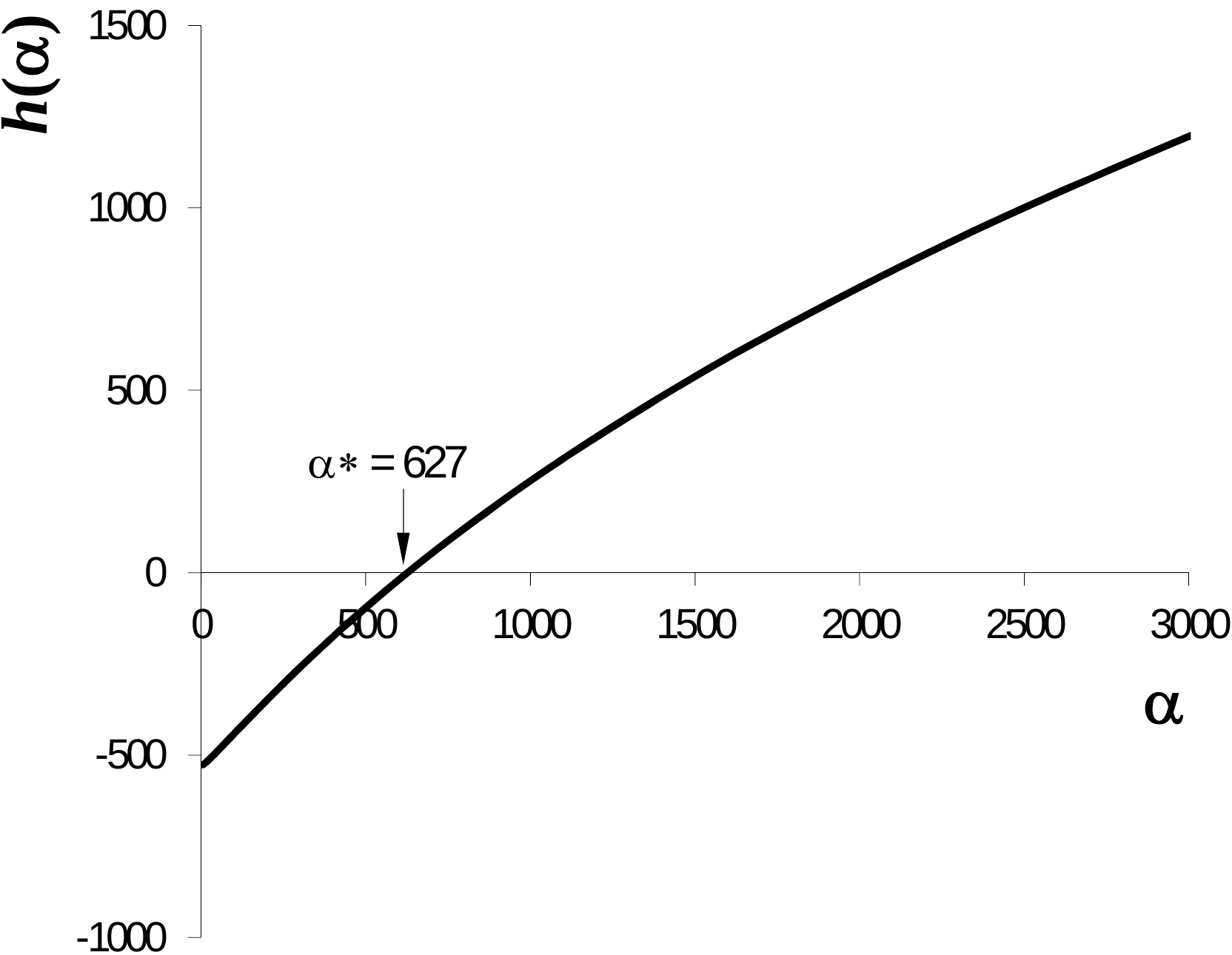}
	\caption{Discrepancy function (2D model).}
	\label{Fig:2DModel-DF}
\end{figure}

The Algorithm~\ref{Algo:NR} was applied to this test problem resulting in $\alpha^\star = 627$, which is illustrated in the plot of the discrepancy function presented in Fig.~\ref{Fig:2DModel-DF}. Using $\alpha^\star = 627$ in the Algorithm~\ref{Algo:GEO2} resulted in $N_a = 7$ and $\alpha_1 = 1087.48$. For comparisons, $N_a = 7$ was also used in the cases with constant inflation and the method proposed in \citep{rafiee:17a}. During all data assimilations a Schur product-based localization \citep{houtekamer:01} was applied to the Kalman gain using a constant correlation length corresponding to the size of 50 gridblocks in the Gaspari-Cohn correlation function. The matrix inversion in Eq.~\ref{Eq:ESMDA} was done using the subspace inversion procedure \citep{evensen:04} keeping 99\% of the sum of the singular values.

\begin{figure}
\centering
    \captionsetup{justification=centering}
    \subfloat[$7\times$-CONST]{
      \includegraphics[width=0.33\textwidth]{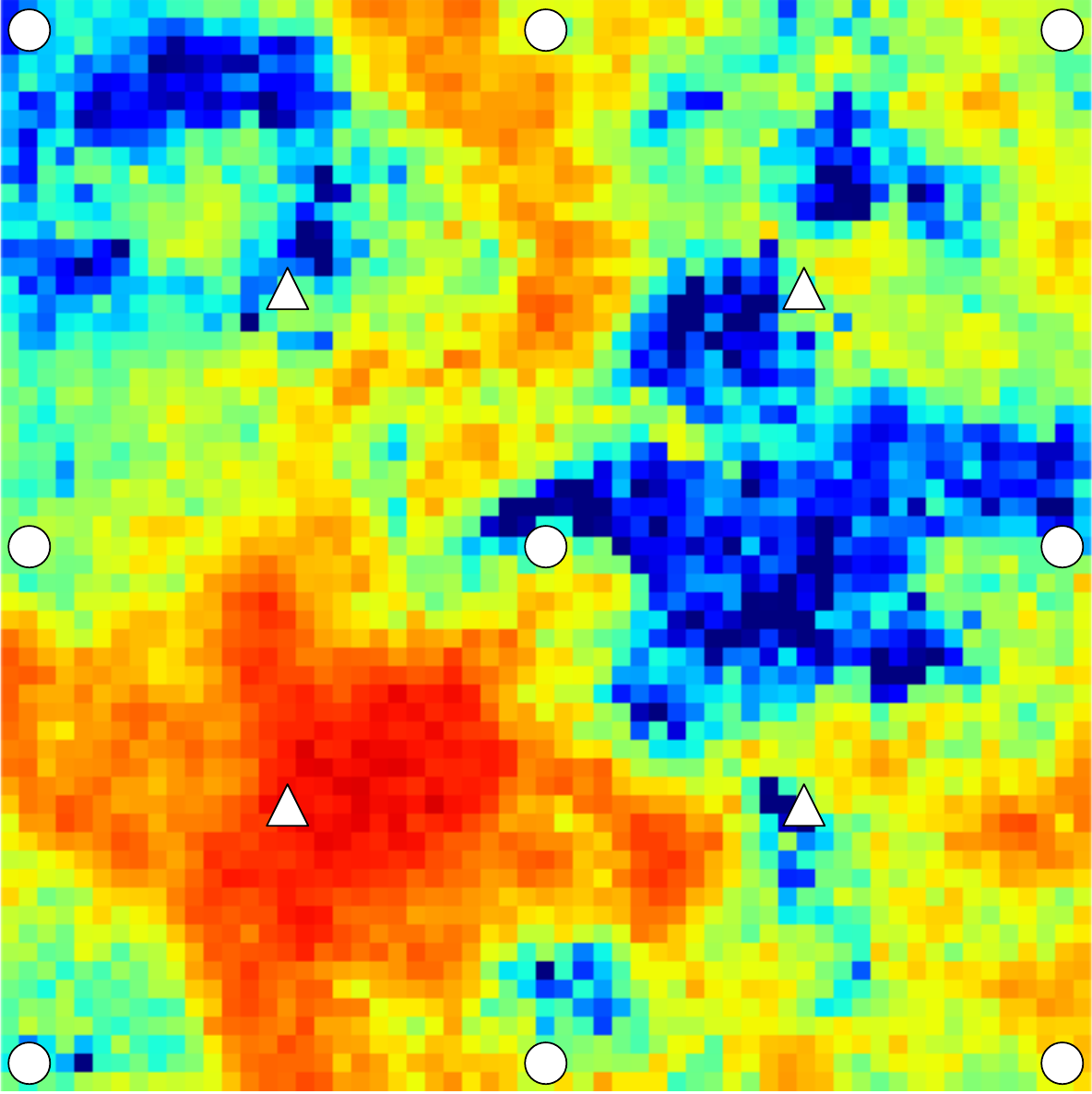}
    }
    \subfloat[$7\times$-GEO1]{
      \includegraphics[width=0.33\textwidth]{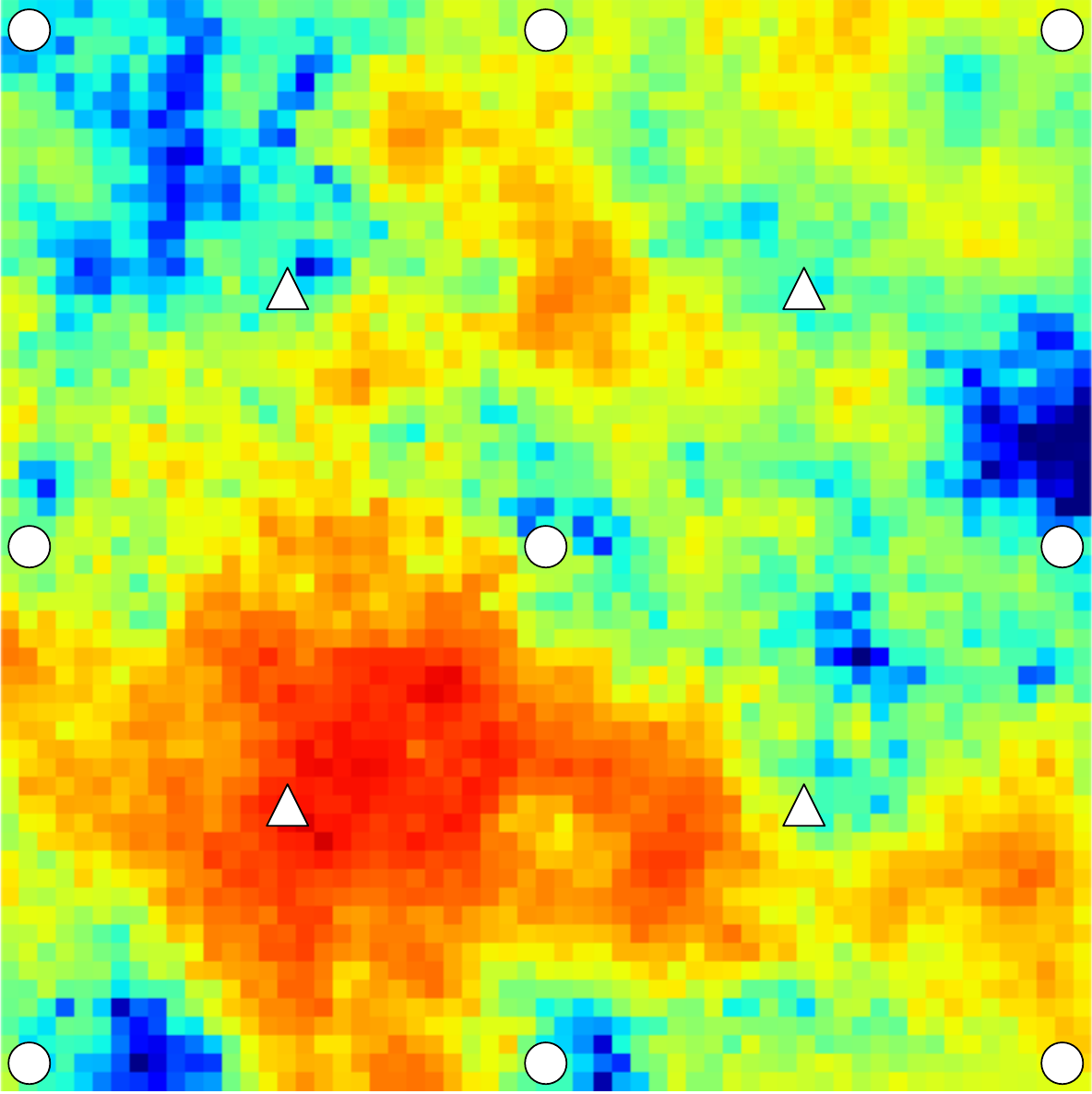}
    }
    \subfloat[$7\times$-GEO2]{
      \includegraphics[width=0.33\textwidth]{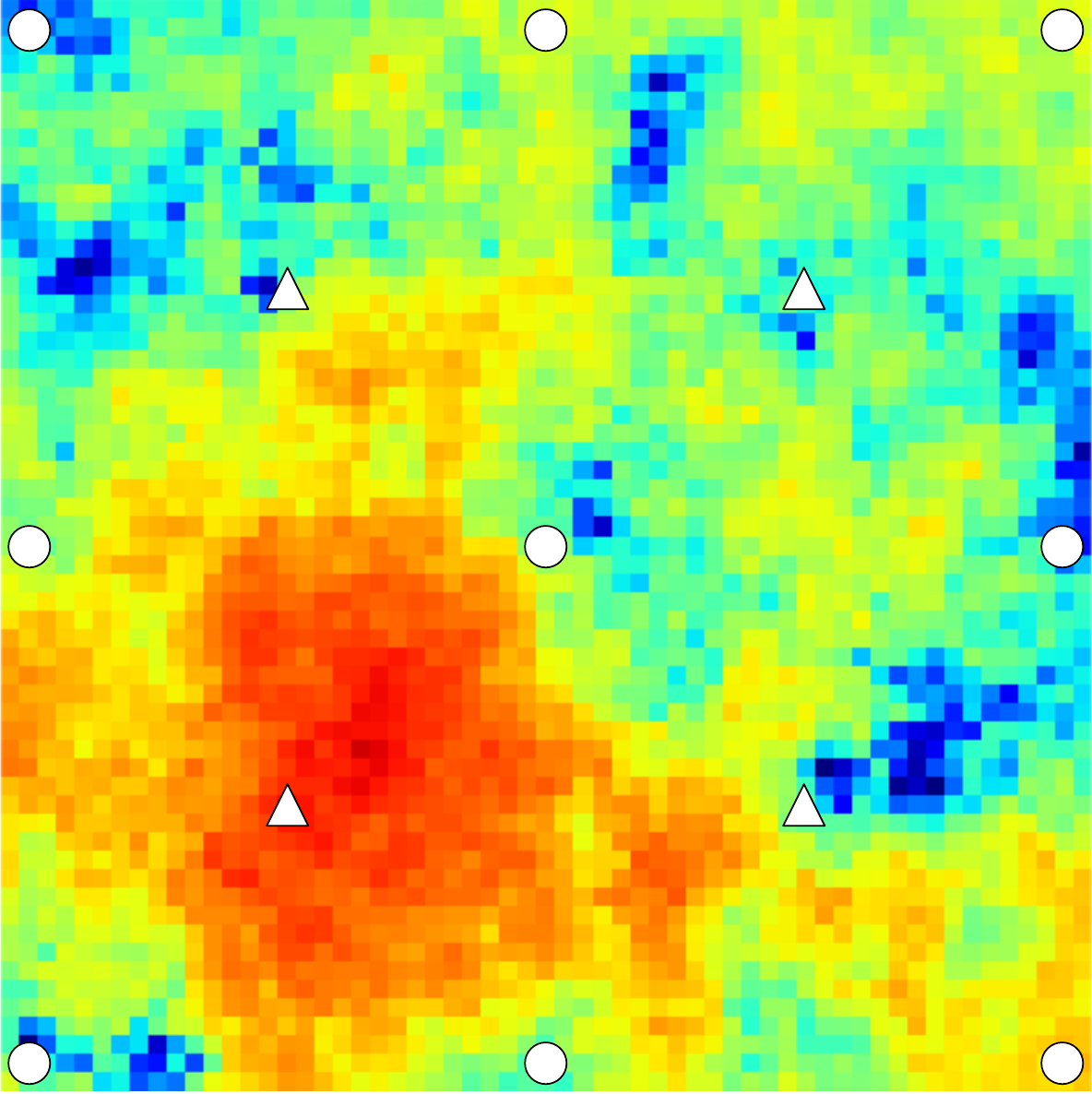}
    }
    \linebreak
      \includegraphics[width=0.5\textwidth]{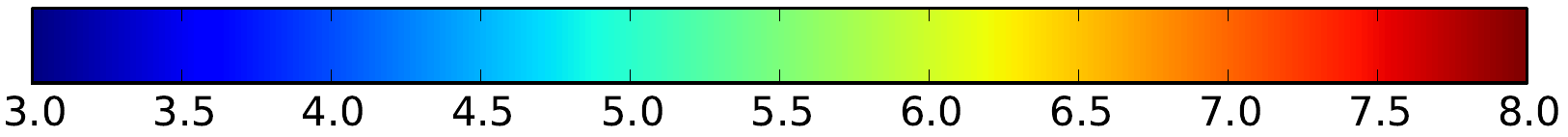}
    \captionsetup{justification=justified}
\caption{Mean log-permeability (2D model).}
\label{Fig:2DModelPearmMean}
\end{figure}

\begin{figure}
\centering
    \captionsetup{justification=centering}
    \subfloat[$7\times$-CONST]{
      \includegraphics[width=0.34\textwidth]{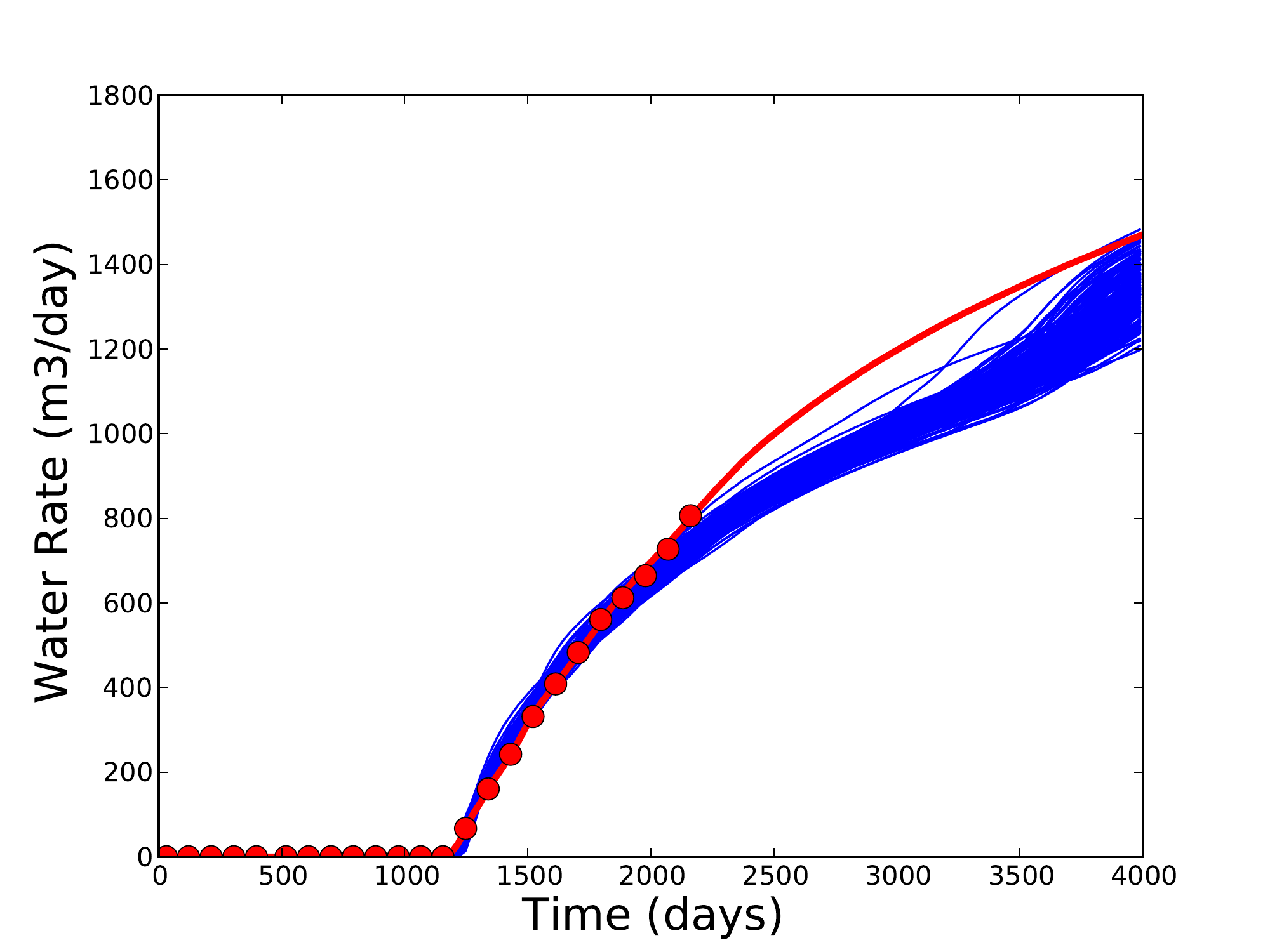}
    }
    \subfloat[$7\times$-GEO1]{
      \includegraphics[width=0.34\textwidth]{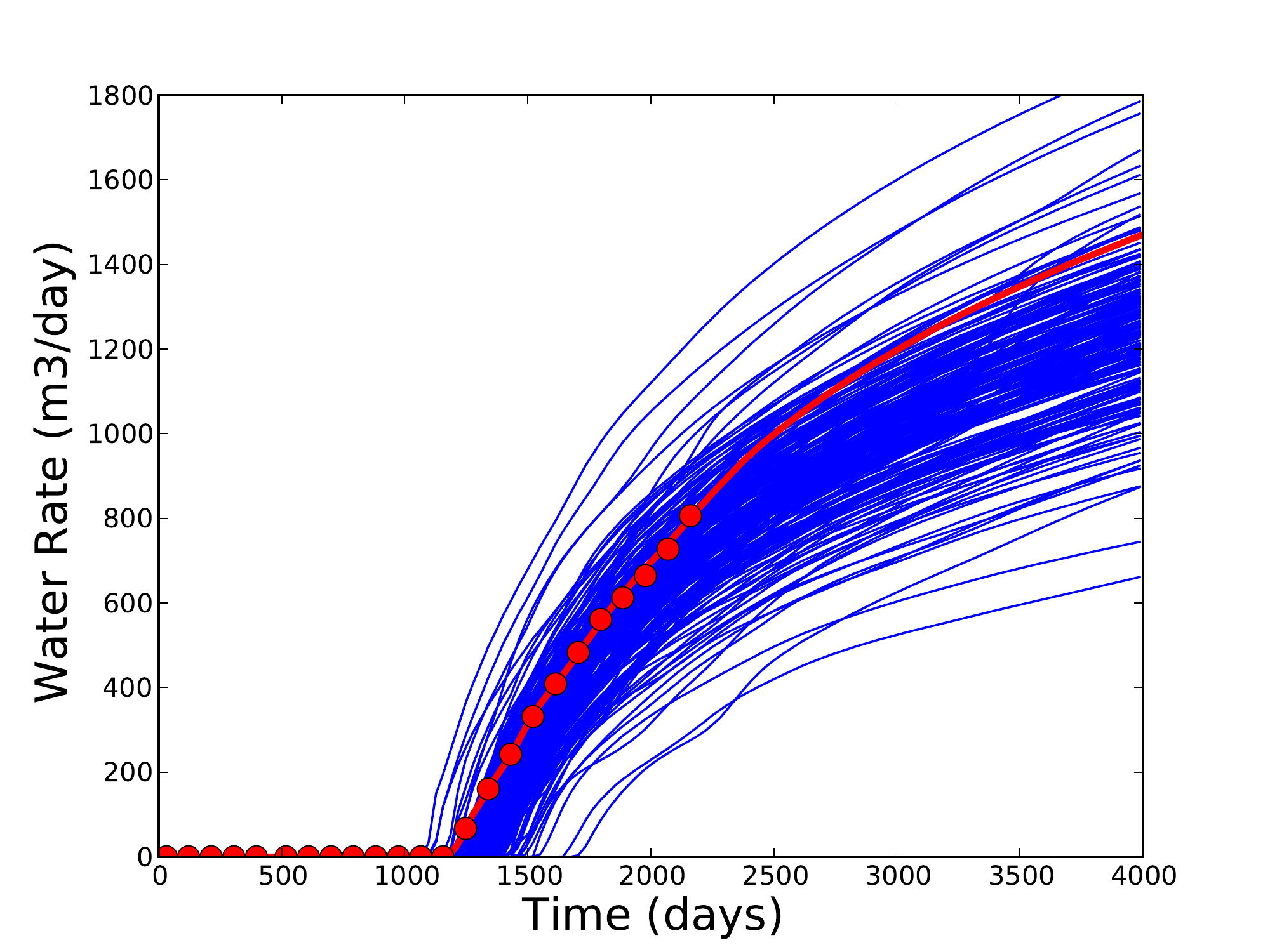}
    }
    \subfloat[$7\times$-GEO2]{
      \includegraphics[width=0.34\textwidth]{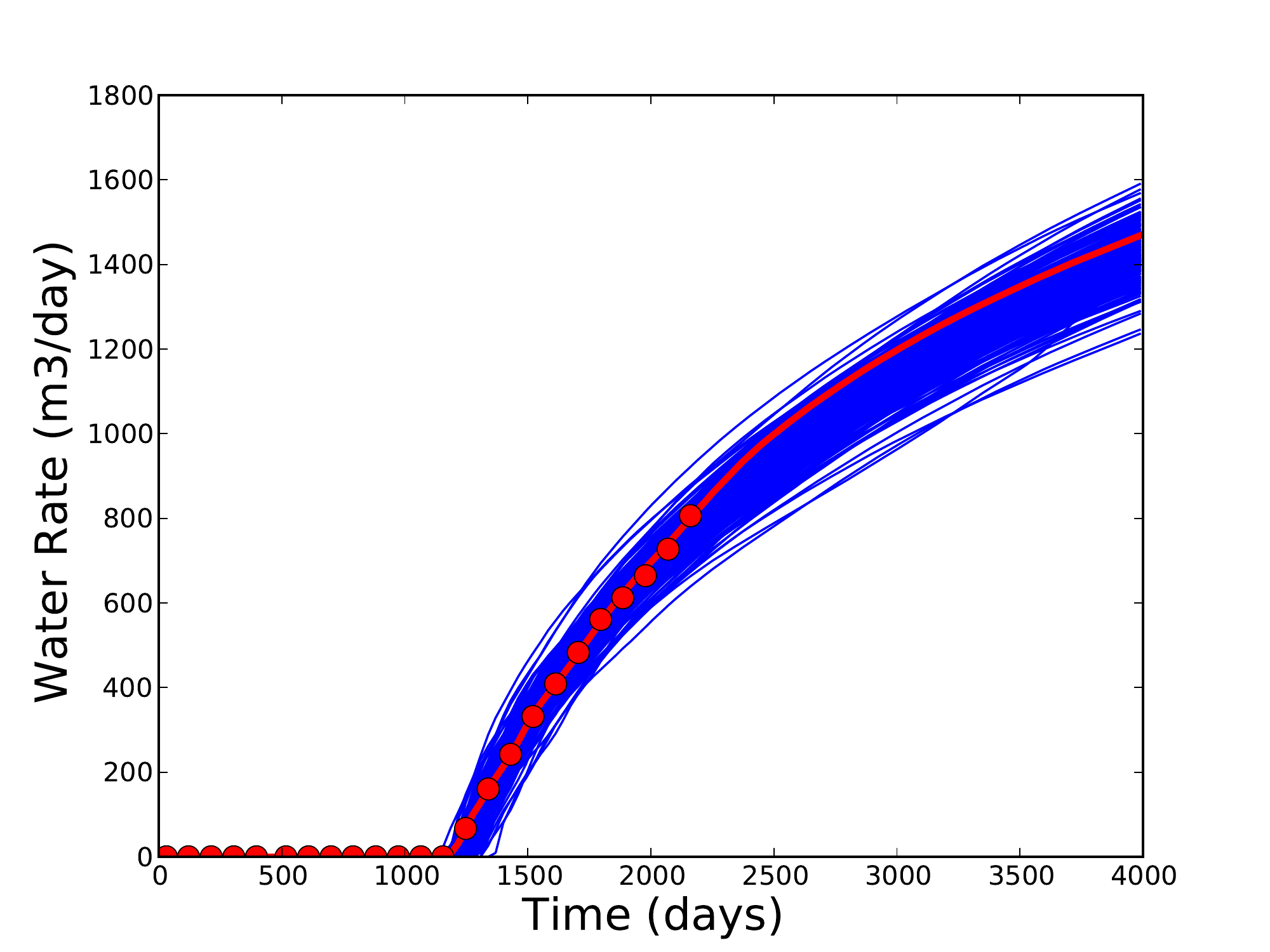}
    }
    \captionsetup{justification=justified}
\caption{Water rate at well P2. Red dots are the observations, red curve is the prediction from the true model and blue curves are the predictions from the posterior ensemble (2D model). }
\label{Fig:2DProdData}
\end{figure}

Table~\ref{Tab:2DModel-Alphas} summarizes the inflation factors for the three cases. The case with constant inflation is labeled as $7\times$-CONST, the method from \citep{rafiee:17a} is labeled as $7\times$-GEO1 while the method proposed in this paper is labeled as $7\times$-GEO2. The method $7\times$-GEO1 resulted in a large initial inflation ($\alpha_1 = 1442941$) and because the number of data assimilations is fixed, the procedure resulted in an inflation of only 1.11 in the last data assimilation step. Fig.~\ref{Fig:2DModelPearmMean} shows the mean log-permeability fields obtained by each method. In all cases, the data assimilation recovered the main permeability features of the true model. However, the permeability resulting from the case $7\times$-CONST shows some regions with overcorrection (very low permeability values). Table~\ref{Tab:2DModel-Results} presents the average values of the root mean squared error (RMSE) and the squared data-mismatch norm divided by the number of data points obtained by each method. The RMSE of the $j$th posterior realization, $\m^{N_a}_j$, is computed as

\begin{equation}\label{Eq:RMSE}
  \text{RMSE}_j = \sqrt{\frac{1}{N_m}} \left\| \m^{N_a}_j - \mtrue \right\|,
\end{equation}
where $\mtrue$ is the true model. The lowest RMSE was obtained by $7\times$-GEO2 indicating a better recovering the true values of log-permeability. The data-mismatch norm is computed as

\begin{equation}\label{Eq:ynorm}
  \left \| \y_j \right\|^2 = \left\| \Ce\invsqr \left(\dobs - \g(\m^{N_a}_j)\right) \right\|^2.
\end{equation}
The case $7\times$-CONST resulted in the lowest average data-mismatch norm indicating that this case obtained the best agreement with the observations. This fact is also illustrated in Fig.~\ref{Fig:2DProdData} which shows the predicted water rate for one well of the field. In addition to the historical period, the plots in Fig.~\ref{Fig:2DProdData} include a five-years forecast period. The results in this figure show that although the case $7\times$-CONST obtained the best data match, it resulted in a narrow span of forecasted water rate, which do not cover the forecast from the true model. The case $7\times$-GEO1 clearly do not match the observations sufficiently well and the case $7\times$-GEO2 resulted in a reasonable compromise between data-match and forecast. Ideally, a successful history-matching method should be able to provide a good sampling of the posterior probability density function of the model parameters. However, it is difficult to compare methods in terms of their sampling performance because generating a reference sampling is very computationally demanding. Even though it is not correct to claim that $7\times$-GEO2 obtained a better sampling than $7\times$-CONST and $7\times$-GEO1, the results show evidences that it resulted in better model estimates and better uncertainty assessment of the production forecast.

\begin{table}
\caption{Inflation factors (2D model)}
\label{Tab:2DModel-Alphas}
\begin{center}
\begin{tabular}{lccc}
\toprule
$k$ & $7\times$-CONST & $7\times$-GEO1 & $7\times$-GEO2  \\
\midrule
1 & 7 & 1442941.18 & 1087.48 \\
2 & 7 & 138031.75 & 362.83 \\
3 & 7 & 13204.12 & 121.05 \\
4 & 7 & 1263.11 & 40.39 \\
5 & 7 & 120.83 & 13.48 \\
6 & 7 & 11.56 & 4.50 \\
7 & 7 & 1.11 & 1.50 \\
\midrule
$\gamma$ & 1 & 0.0957& 0.3336 \\
\bottomrule
\end{tabular}
\end{center}
\end{table}

\begin{table}
\caption{Average metrics (2D model)}
\label{Tab:2DModel-Results}
\begin{center}
\begin{tabular}{lccc}
\toprule
$~$ & $7\times$-CONST & $7\times$-GEO1 & $7\times$-GEO2  \\
\midrule
RMSE & 1.714 & 2.227 & 1.680 \\
Squared data-mismatch norm & 3.390 & 63.359 & 14.828 \\
\bottomrule
\end{tabular}
\end{center}
\end{table}

\subsection{UNISIM-I-H}
\label{Sec:UNISIM}

\begin{figure}
	\centering
	\includegraphics[width=0.7\linewidth]{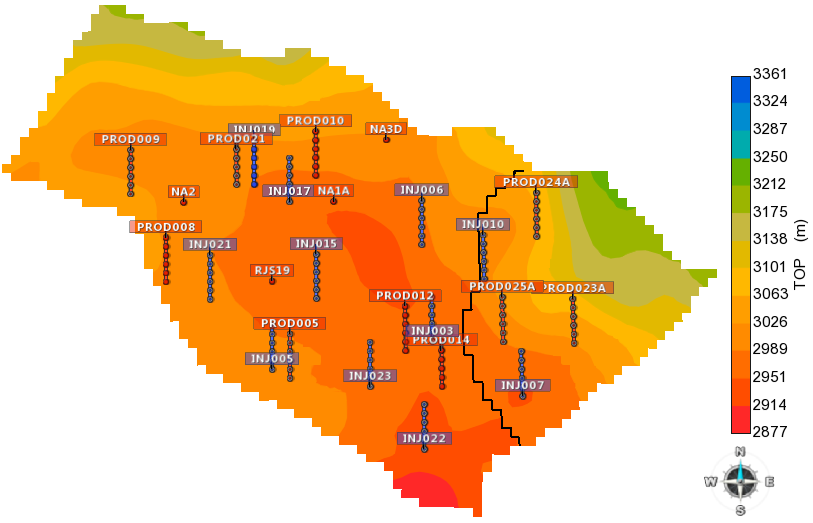}
	\caption{Position of the wells (UNISIM-I-H)}
	\label{Fig:UNISIM-Wells}
\end{figure}

The second test case is the benchmark problem UNISIM-I-H \citep{avansi:15a}, which is available in \citep{unisim-i-h:13}. The benchmark was constructed with actual data from Namorado Field (Campos Basis, Brazil). The UNISIM-I-H model has $81 \times 58 \times 20$ gridblocks, with dimensions $100 \times 100 \times 8$~meters and 37,000 active gridblocks. A total of 500 realizations of petrophysical properties (porosity, permeability in the three orthogonal directions and net-to-gross ratio) is provided with the dataset. There are 14 oil producing and 11 water injection wells. Fig.~\ref{Fig:UNISIM-Wells} shows the position of the wells projected in the first layer of the model. Note that the actual wells are perforated in different layers of the model, typically water injectors are perforated in bottom while oil producing wells are perforated in the top of the reservoir. The synthetic observed data provided with the dataset correspond to monthly ``measurements'' of water cut and bottom-hole pressure (BHP) for the oil producing wells and BHP for the water injection wells for a period of 4018 days. These observations were generated based on a fine-scale model with 3.5 million active gridblocks and are corrupted with noise. Here, the noise in each individual datum was assumed to be an independent sample from a Gaussian distribution with zero mean and standard deviation corresponding to 10\% of the water cut data and 10~$\text{kgf}/\text{cm}^2$ for pressure measurements. Besides the petrophysical properties, the UNISIM-I-H case also includes six global parameters, whose prior uncertainties were modeled as independent triangle distributions with values presented in Table~\ref{Tab:UNISIM-ScalarPar}. A detailed description of the UNISIM-I-H case can be found in \citep{avansi:15a,unisim-i-h:13}.

\begin{table}
\caption{Prior distribution of global parameters (UNISIM-I-H)}
\label{Tab:UNISIM-ScalarPar}
\begin{scriptsize}
\begin{center}
\begin{tabular}{lccc}
\toprule
Parameter & Mode & Min. & Max. \\
\midrule
Rock compressibility in $(\text{cm}^2/\text{kgf})$ & $5.3\times 10^{-5}$ & $1\times 10^{-5}$ & $9.6\times 10^{-5}$ \\
Depth of the oil-water contact at the East block (Fig.~\ref{Fig:UNISIM-Wells}) & 3174 & 3169 & 3179 \\
Maximum water relative permeability & 0.35 & 0.15 & 0.55 \\
Corey exponent of water relative permeability & 2.4 & 2.0 & 3.0 \\
Multiplier for vertical permeability & 1.5 & 0.0 & 3.0 \\
\bottomrule
\end{tabular}
\end{center}
\end{scriptsize}
\end{table}

\begin{figure}
	\centering
	\includegraphics[width=0.5\linewidth]{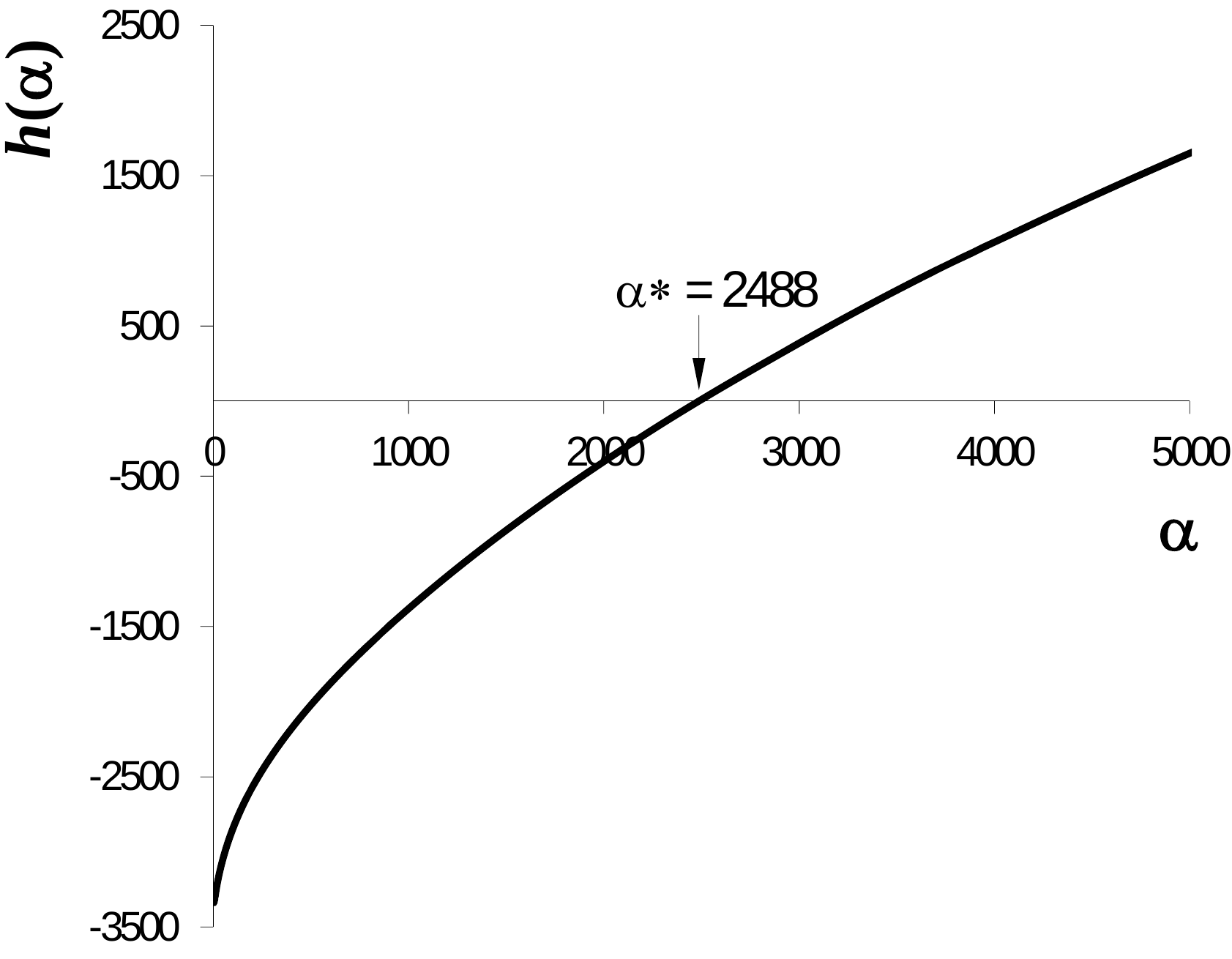}
	\caption{Discrepancy function (UNISIM-I-H).}
	\label{Fig:UNISIM-DF}
\end{figure}

ES-MDA was applied to update the prior ensemble considering 13 different configurations of inflation factors. Table~\ref{Tab:UNSIM-Alphas} summarizes the cases considered. In this table, $4\times$-CONST and $8\times$-CONST stand for ES-MDA with constant inflation factors with four and eight data assimilations, respectively. $4\times$-GEO(1E2), $4\times$-GEO(1E3), $4\times$-GEO(1E4) and $4\times$-GEO(1E5) stand for ES-MDA with four data assimilations and geometric selection of the inflations, where the values of the first inflation are given between the braces, i.e., $\alpha_1 = 10^2$, $\alpha_1 = 10^3$, $\alpha_1 = 10^4$ and $\alpha_1 = 10^5$, respectively. The same initial values for $\alpha_1$ are also used for the cases with eight data assimilations. The method proposed in \citep{rafiee:17a} was tested with four and eight data assimilations ($4\times$-GEO1 and $8\times$-GEO1). The method proposed in this paper required eight data assimilations to satisfy the MDP and it is labeled as $8\times$-GEO2. Fig.~\ref{Fig:UNISIM-DF} shows a plot of the discrepancy function indicating that the inflation value such that $h(\alpha) = 0$ is 1172. Using $\alpha^\star = 1172$ in the Algorithm~\ref{Algo:GEO2} resulted in $N_a = 8$, $\alpha_1 = 3273.79$ and $\gamma = 0.3334$. In all cases, Kalman gain localization was applied using the Gaspari-Cohn correlation function with correlation length of 2000 meters. Subspace inversion keeping 99\% of the sum of the singular values was applied.

\begin{table}
\caption{Inflation factors (UNISIM-I-H)}
\label{Tab:UNSIM-Alphas}
\begin{scriptsize}
\begin{center}
\begin{tabular}{lccccccc}
\toprule
$k$ & $4\times$-CONST & $4\times$-GEO(1E2) & $4\times$-GEO(1E3) & $4\times$-GEO(1E4) & $4\times$-GEO(1E5) & $4\times$-GEO1 & ~ \\
\midrule
1 & 4 & 100 & 1000 & 10000 & 100000 & 4010.30 & ~ \\
2 & 4 & 23.54 & 103.71 & 471.69 & 2172.79 & 258.07 & ~ \\
3 & 4 & 5.54 & 10.76 & 22.25 & 47.21 & 16.61 & ~ \\
4 & 4 & 1.30 & 1.12 & 1.05 & 1.03 & 1.07 & ~ \\
\midrule
$\gamma$ & 1 & 0.2354 & 0.1037 & 0.0472 & 0.0217 & 0.0644 & ~ \\
\bottomrule
\toprule
$k$ & $8\times$-CONST & $8\times$-GEO(1E2) & $8\times$-GEO(1E3) & $8\times$-GEO(1E4) & $8\times$-GEO(1E5) & $8\times$-GEO1 & $8\times$-GEO2  \\
\midrule
1 & 8 & 100 & 1000 & 10000 & 100000 & 4010.30 & 3273.79\\
2 & 8 & 58.64 & 401.08 & 2812.60 & 19929.85 & 1296.25 & 1091.58\\
3 & 8 & 34.39 & 160.87 & 791.07 & 3971.99 & 418.99 & 363.96\\
4 & 8 & 20.17 & 64.52 & 222.50 & 791.61 & 135.43 & 121.36\\
5 & 8 & 11.83 & 25.88 & 62.58 & 157.77 & 43.77 & 40.46\\
6 & 8 & 6.94 & 10.38 & 17.60 & 31.44 & 14.15 & 13.49\\
7 & 8 & 4.07 & 4.16 & 4.95 & 6.27 & 4.57 & 4.50\\
8 & 8 & 2.39 & 1.67 & 1.39 & 1.25 & 1.48 & 1.50 \\
\midrule
$\gamma$ & 1 & 0.5864 & 0.4011 & 0.2813 & 0.1993 & 0.3232 & 0.3334 \\
\bottomrule
\end{tabular}
\end{center}
\end{scriptsize}
\end{table}

Even though a rigorous comparison in terms of the sampling performance is not possible because the extreme computational requirements to generate a reference distribution, it is possible to compare history-matching methods in terms of some desirable properties. One important quality for a good history-matching method is the ability to achieve a good data match at a reasonable computational cost. Here, the data match quality is evaluated in terms of the data-mismatch norm $\left \| \y_j \right\|$. Another important quality of a good history-matching method is the ability to preserve the geological realism. In this sense, it is desirable that the method makes the smallest changes possible in the prior realizations to match data if these realizations were generated from a rigourous geomodeling process. The model changes can be evaluated with the norm $\left\| \delta \m_j  \right\|$, where $\delta\m_j = \mathbf{S}_\m \left( \m^{N_a}_j - \m^0_j \right)$ with $\m^0_j$ denoting the $j$th prior realization and $\mathbf{S}_\m$ denoting a diagonal matrix containing the inverse of the prior standard deviation of $\m$. History matching practitioners usually give more emphasis on the data match quality. However, minimizing the amount of changes in the geological model is also important. Unfortunately, these objectives are often conflicting. In this sense, it is desirable to obtain models with a good balance between the two norms.

\begin{figure}
	\centering
	\includegraphics[width=0.65\linewidth]{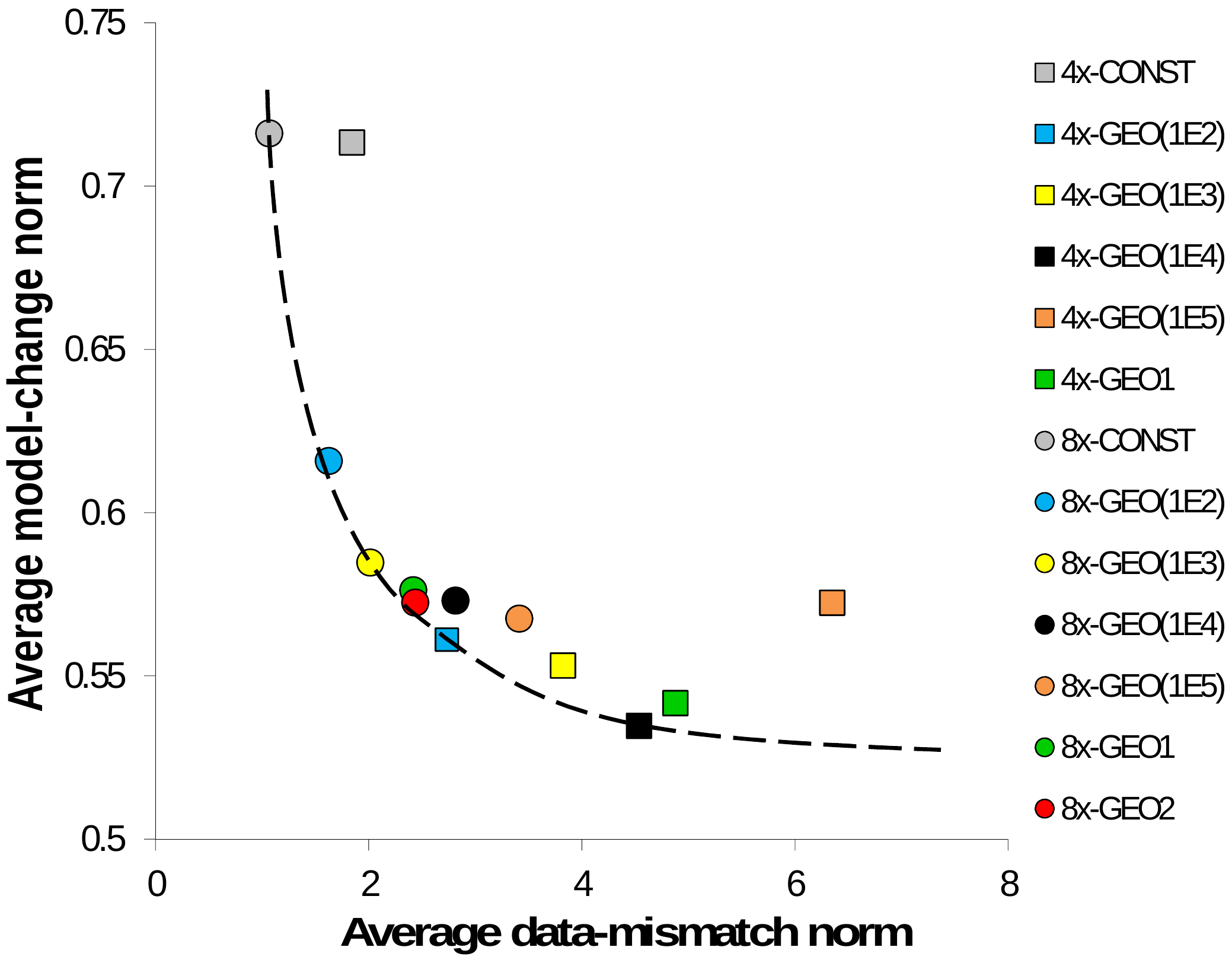}
	\caption{Average data-mismatch and model-change norms (UNISIM-I-H).}
	\label{Fig:UNISIM-Norms}
\end{figure}

Fig.~\ref{Fig:UNISIM-Norms} plots the values of the average of the square of the data-mismatch norm, $\overline{\left \| \y_j \right\|^2}$ divided by the number of data points, against the average squared-norm of the model change, $\overline{\left \| \delta\m_j \right\|^2}$ divided by the number of model parameters. The model-change norm presented in Fig.~\ref{Fig:UNISIM-Norms} was computed only for the log-permeability in the $I$-direction of the grid. Similar plots were obtained for the other petrophysical properties, but they are omitted here. The dashed curve in Fig.~\ref{Fig:UNISIM-Norms} was drawn connecting the ``best'' points considering the two norms to generate the analogous of a Pareto frontier. This procedure is closely related to the L-curve method; see, e.g. \citep{hansen:93a}, which is another method described in the literature to select regularization parameters for linear inverse problems. For example, note that the case $8\times$-CONST resulted in the smallest data-mismatch norm, but at a cost of a large change in the model. The cases closer to the region of maximum curvature correspond to the ones with the best balance between data match and model change. The results in Fig.~\ref{Fig:UNISIM-Norms} indicate that $\alpha_1$ and $N_a$ must be selected in conjunction. For example, for $N_a = 4$, it seems that the best choice of $\alpha_1$ is on the order of $10^2$. However, for $N_a = 8$, the best choice of $\alpha_1$ is on the order to $10^3$. In fact, the proposal from \citet{rafiee:17a} obtained a particularly well-balanced result for the case with eight data assimilations ($8\times$-GEO1). However, the case $4\times$-GEO1 resulted in a larger data-mismatch norm, making this case distant from the maximum-curvature region. The procedure proposed in this paper ($8\times$-GEO2) obtained very similar results to the case $8\times$-GEO1 showing a good compromise between data match and model change.

\begin{figure}
\centering
    \captionsetup{justification=centering}
    \subfloat[Prior]{
      \includegraphics[width=0.5\textwidth]{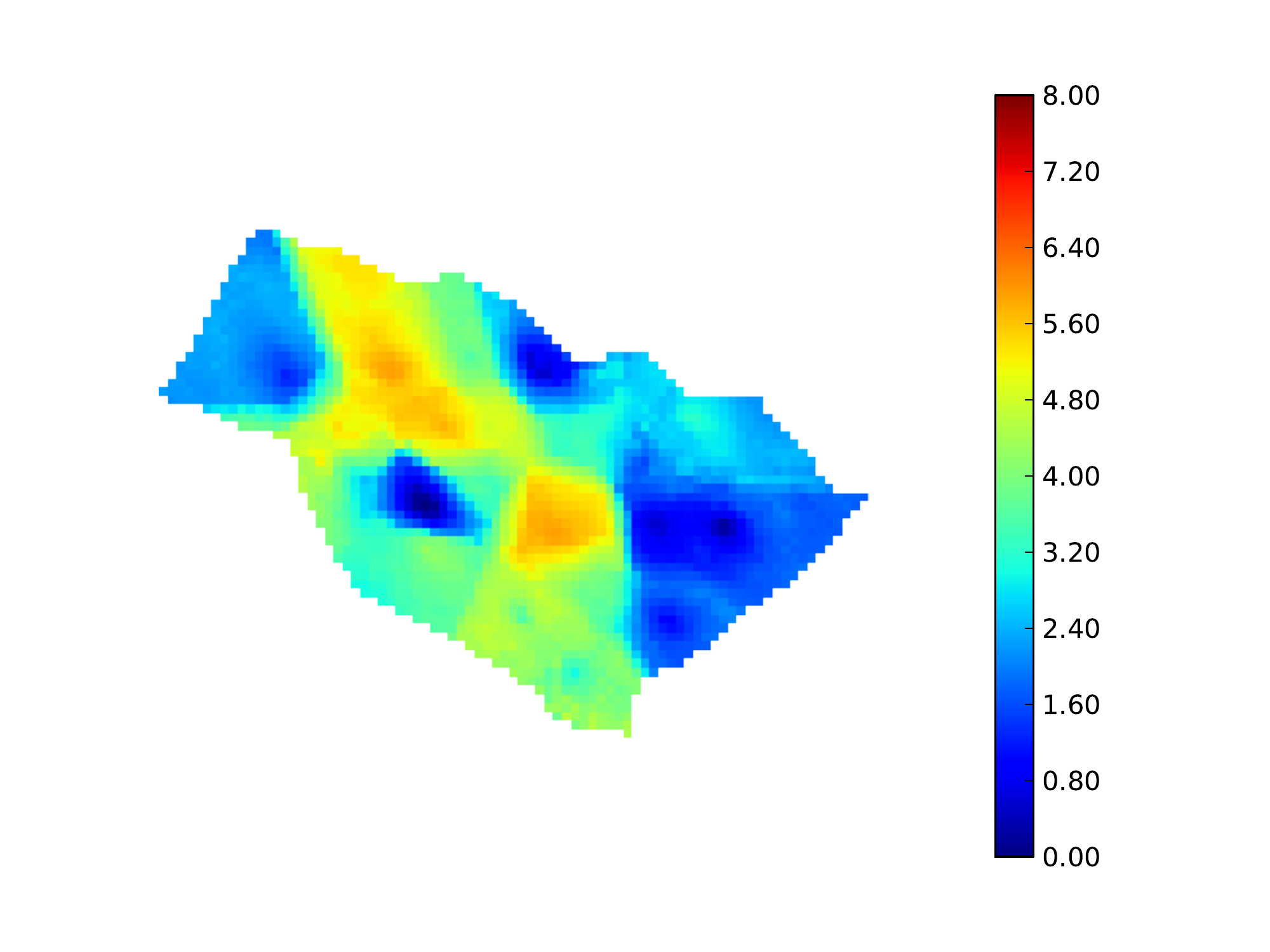}
    }
    \subfloat[$8\times$-CONST]{
      \includegraphics[width=0.5\textwidth]{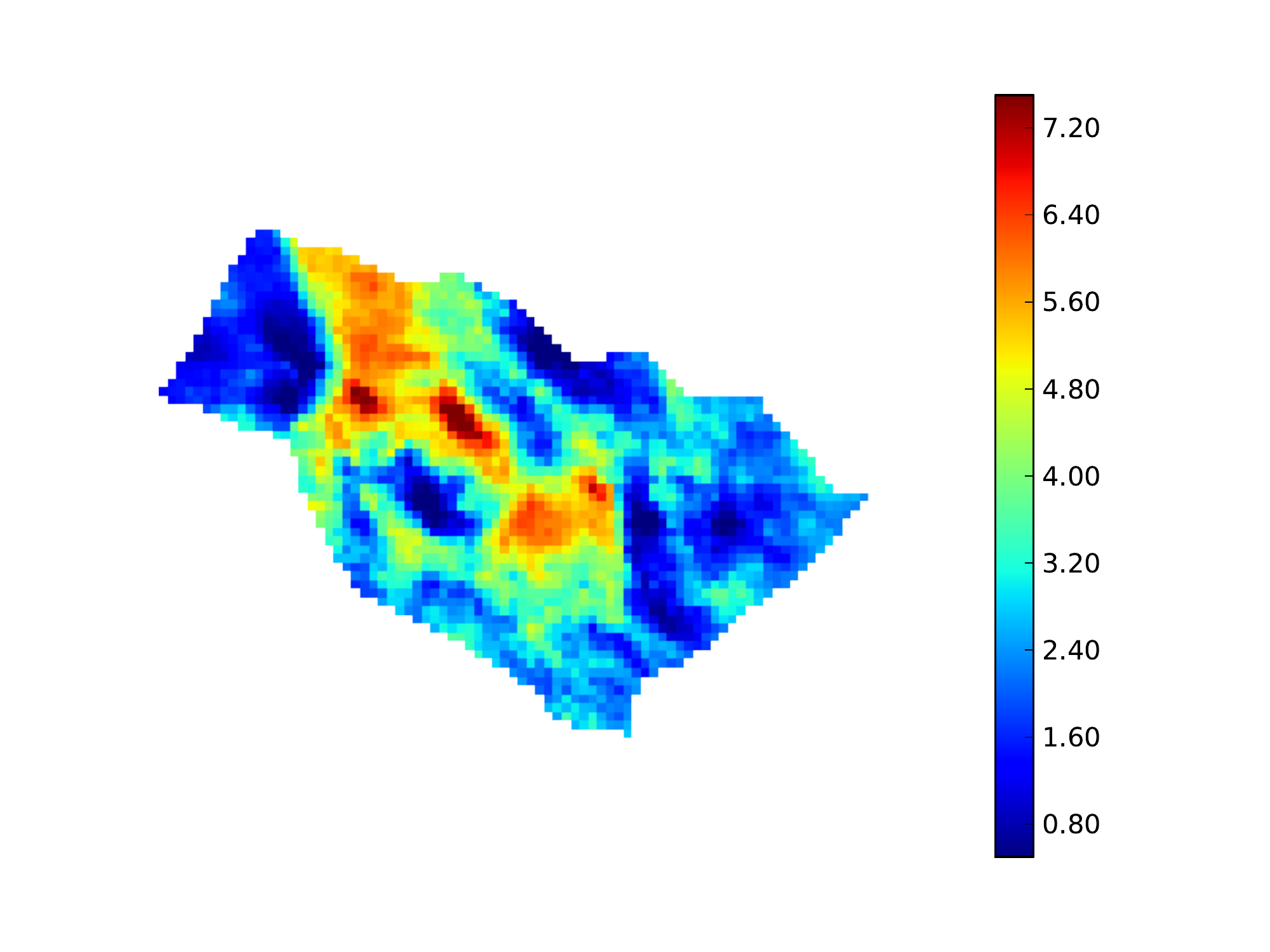}
    }
    \linebreak
    \subfloat[$8\times$-GEO1]{
      \includegraphics[width=0.5\textwidth]{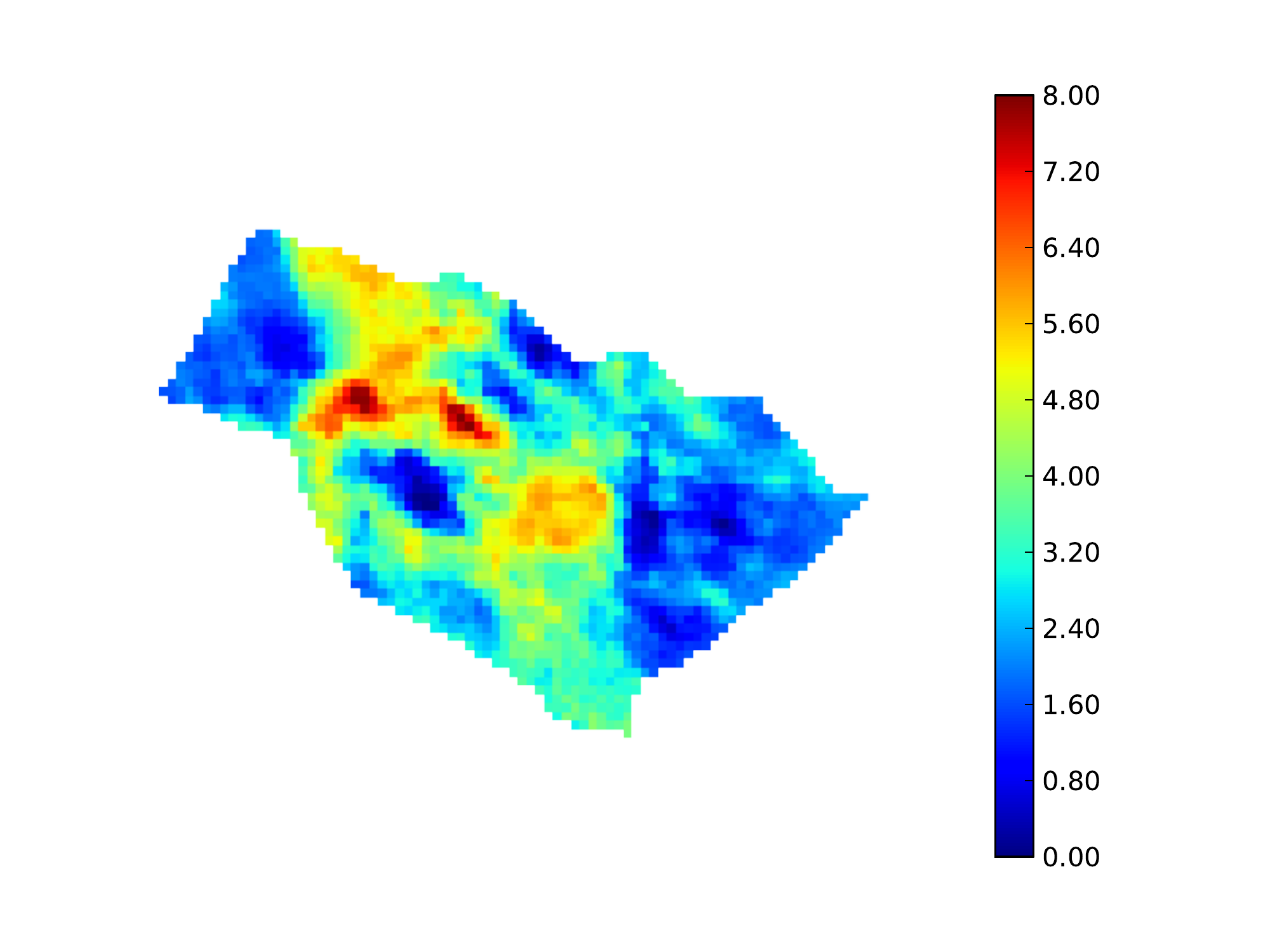}
    }
    \subfloat[$8\times$-GEO2]{
      \includegraphics[width=0.5\textwidth]{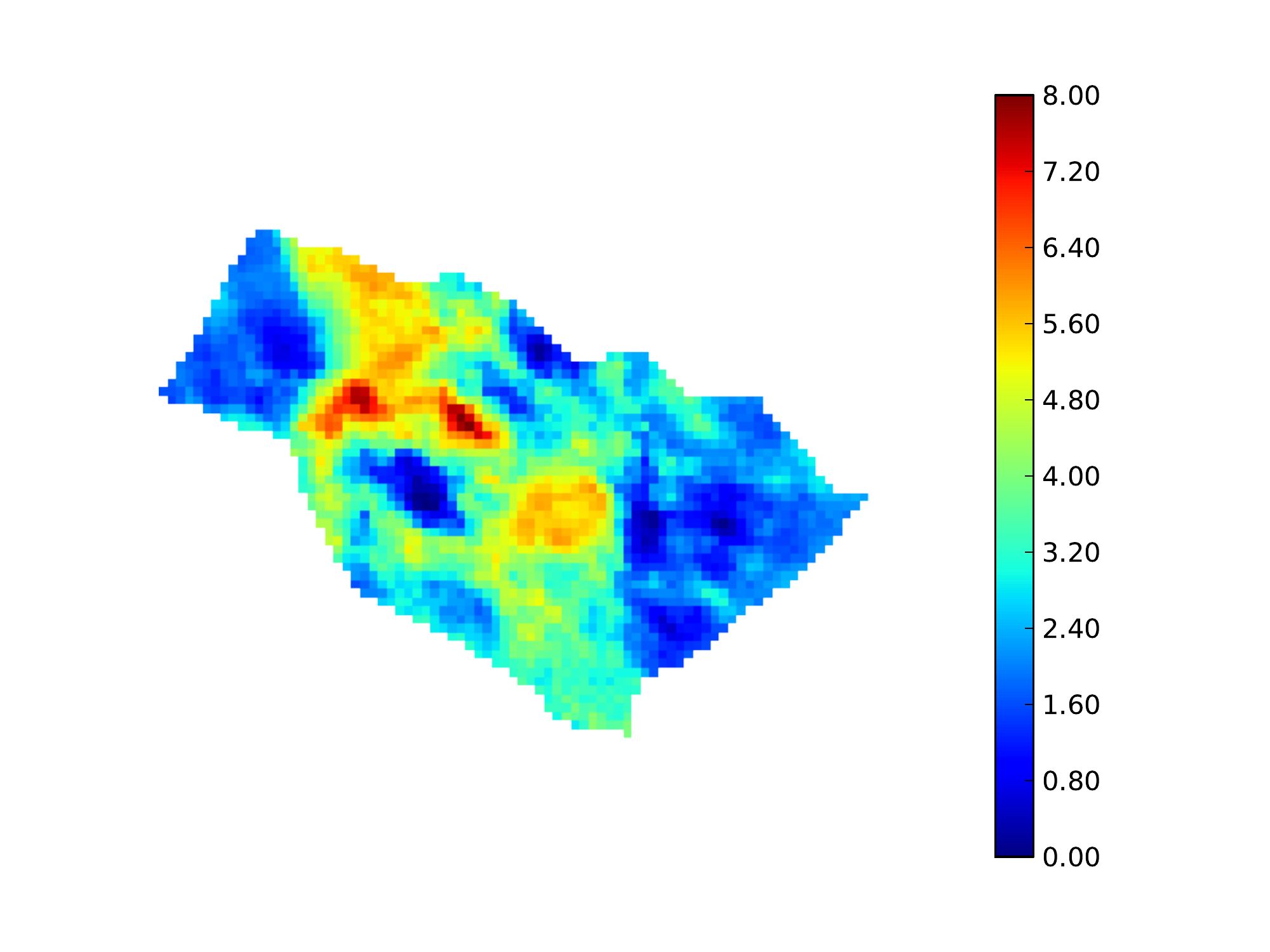}
    }
    \linebreak
      \includegraphics[width=0.5\textwidth]{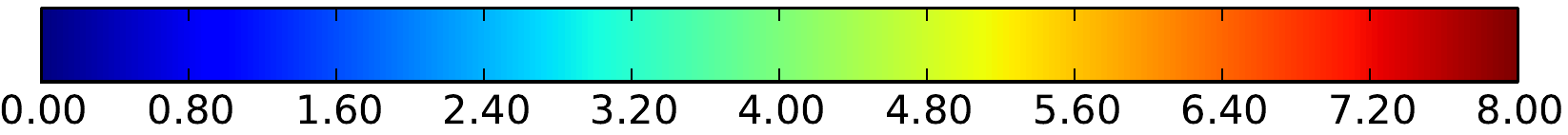}
    \captionsetup{justification=justified}
\caption{Mean log-permeability (UNISIM-I-H, layer 12).}
\label{Fig:UNISIM-PERMI1}
\end{figure}

\begin{figure}
\centering
    \captionsetup{justification=centering}
    \subfloat[$8\times$-CONST]{
      \includegraphics[width=0.5\textwidth]{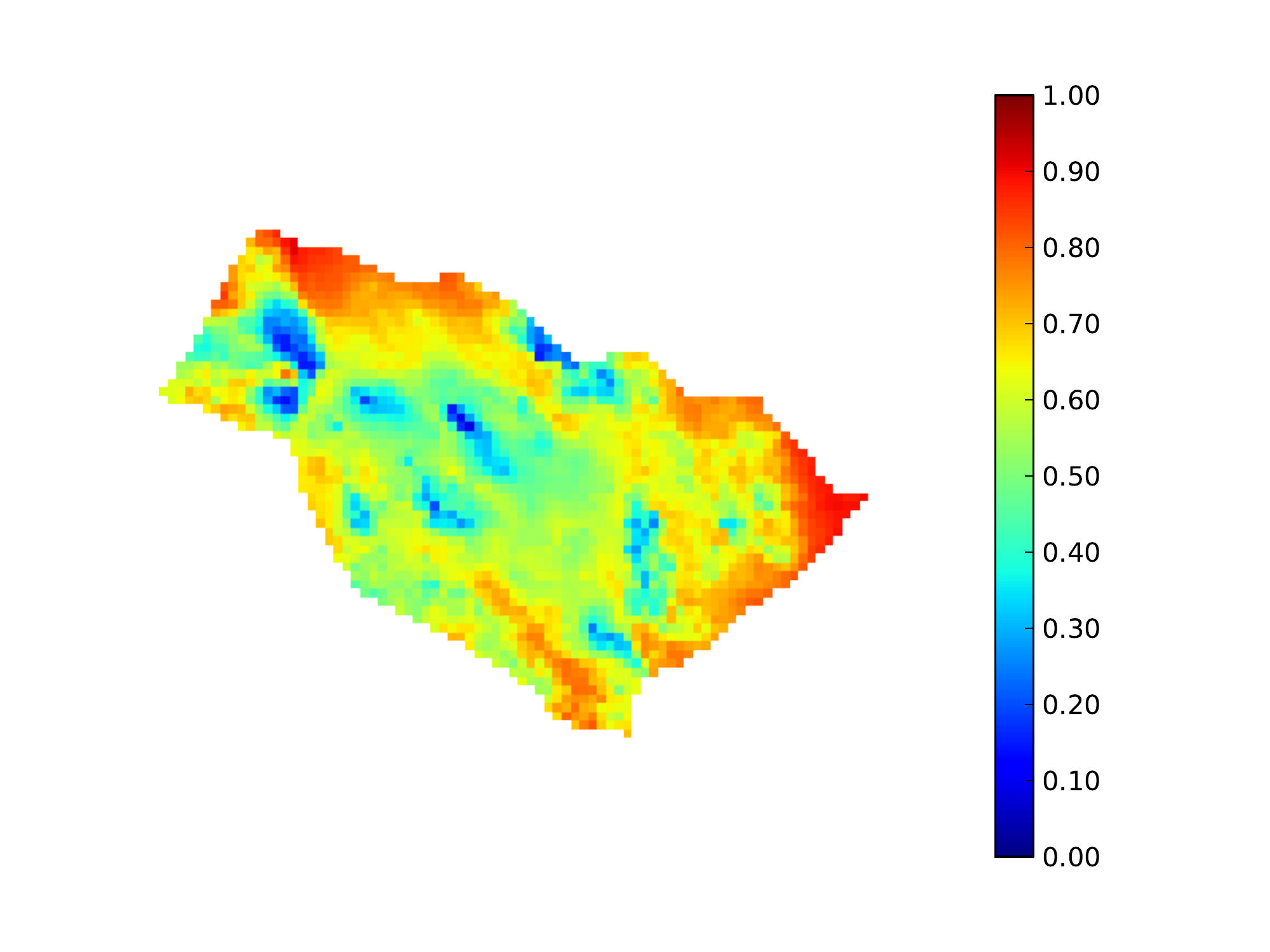}
    }
    \subfloat[$8\times$-GEO1]{
      \includegraphics[width=0.5\textwidth]{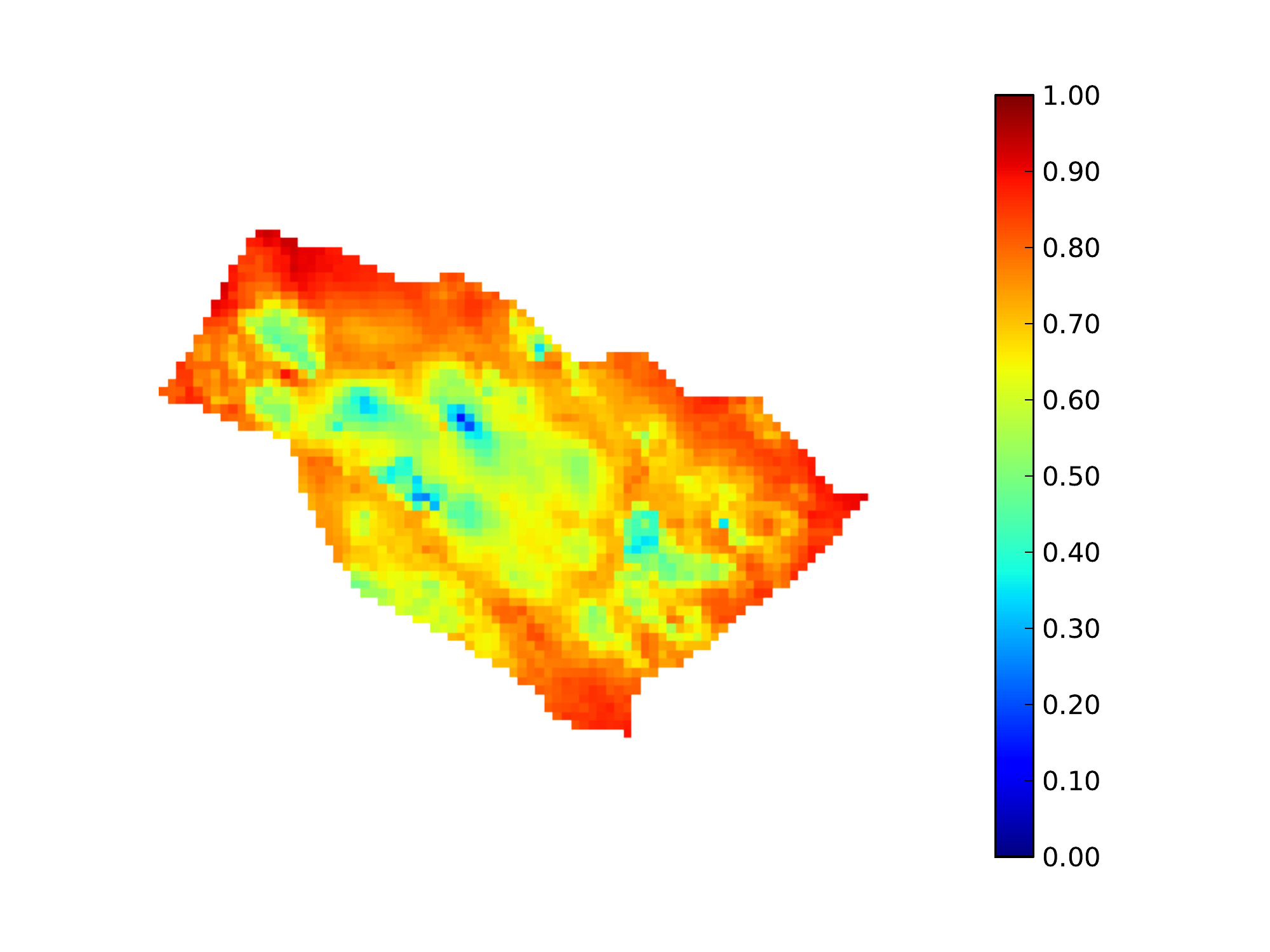}
    }
    \linebreak
    \subfloat[$8\times$-GEO2]{
      \includegraphics[width=0.5\textwidth]{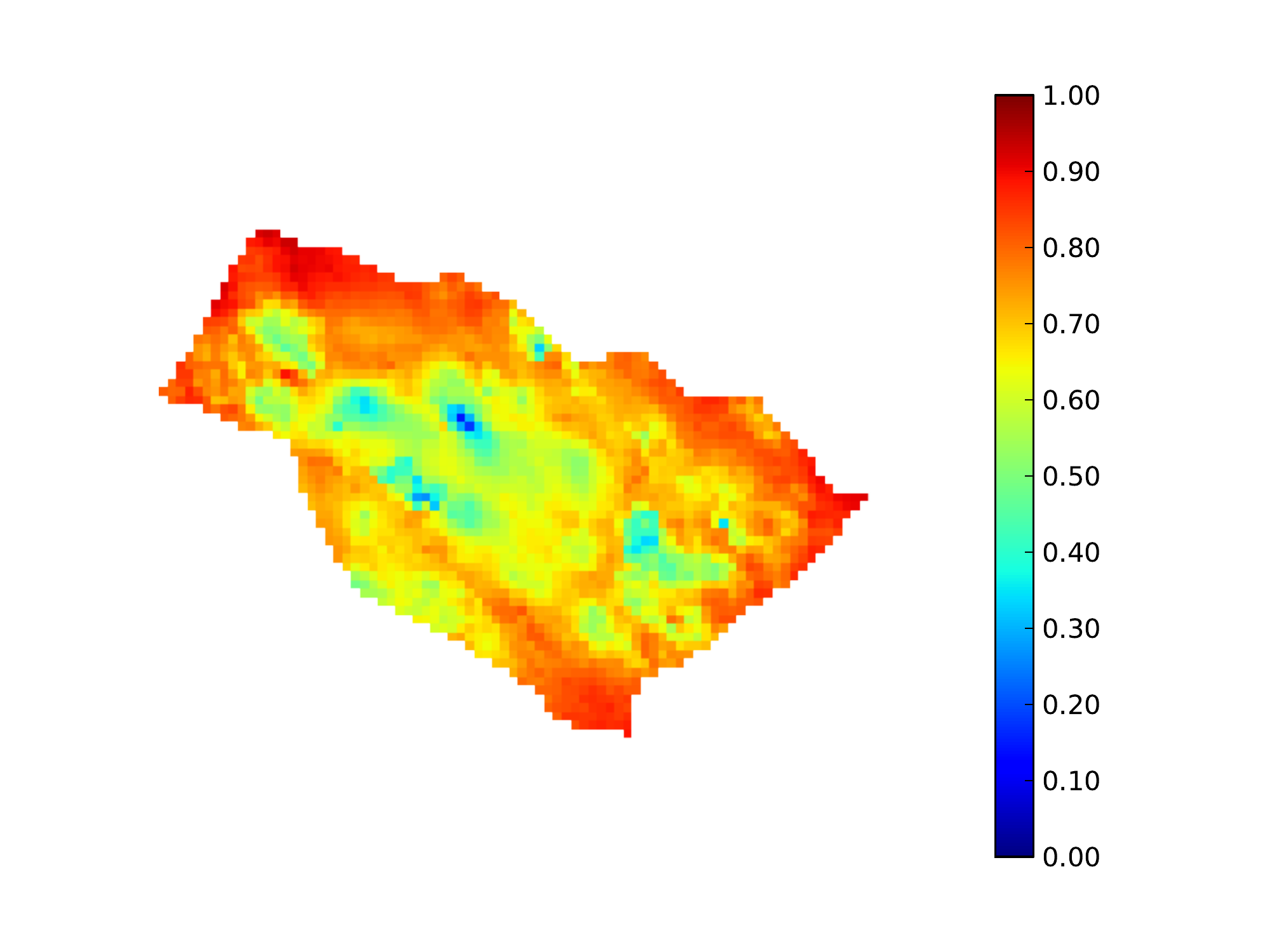}
    }
    \linebreak
      \includegraphics[width=0.5\textwidth]{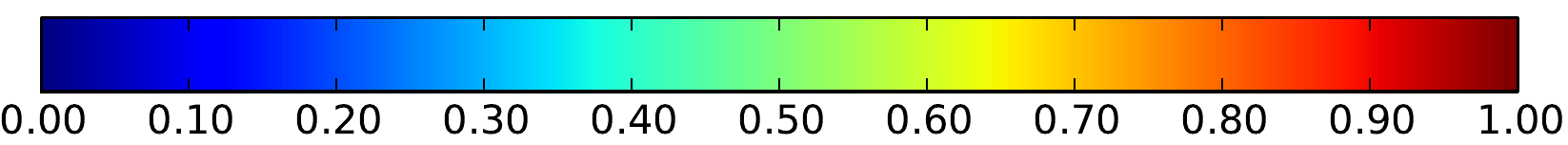}
    \captionsetup{justification=justified}
\caption{Normalized variance of log-permeability (UNISIM-I-H, layer 12).}
\label{Fig:UNISIM-NV}
\end{figure}

\begin{figure}
\centering
    \captionsetup{justification=centering}
    \subfloat[$8\times$-CONST, Water cut]{
      \includegraphics[width=0.4\textwidth]{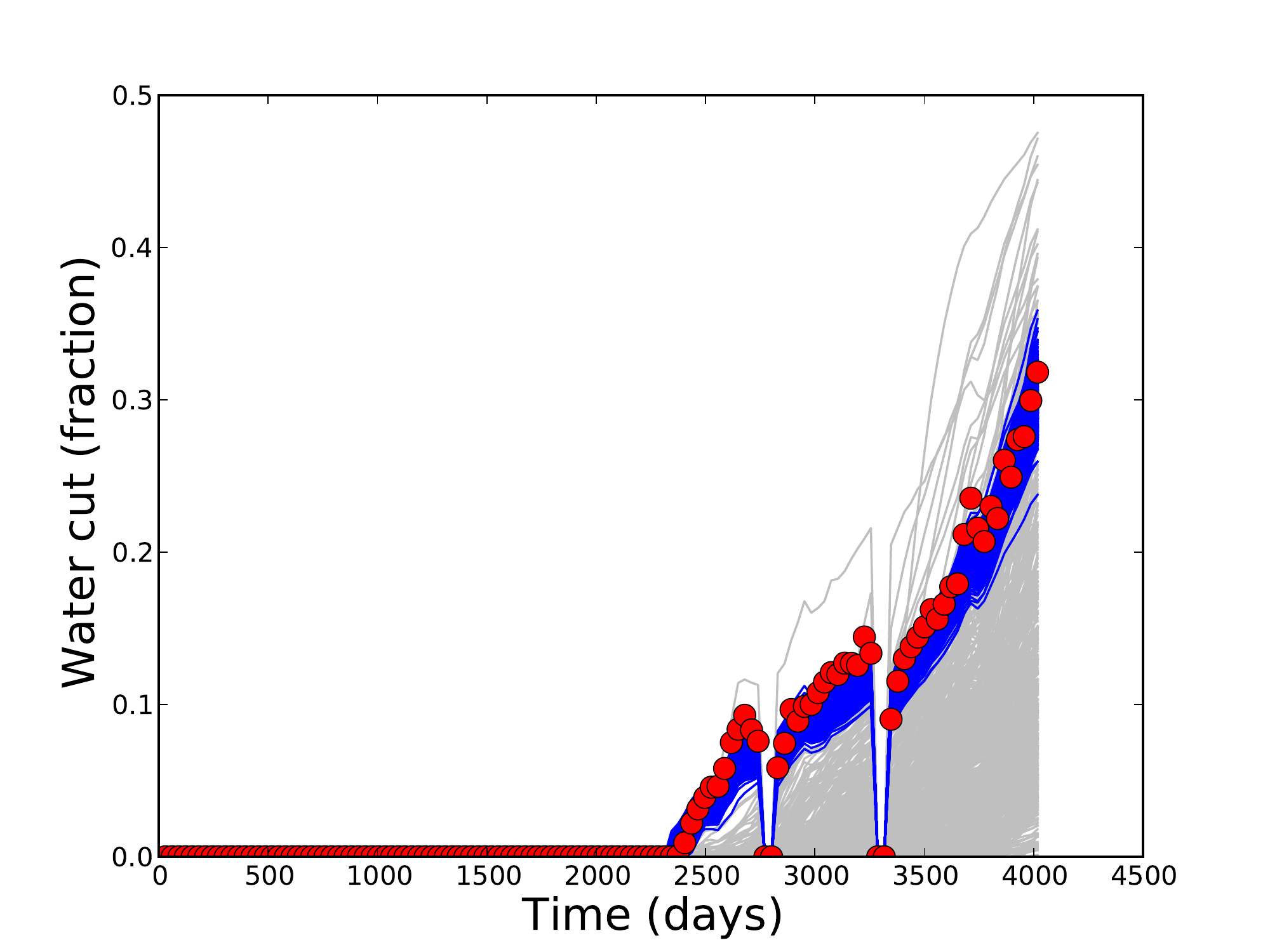}
    }
    \subfloat[$8\times$-CONST, BHP]{
      \includegraphics[width=0.4\textwidth]{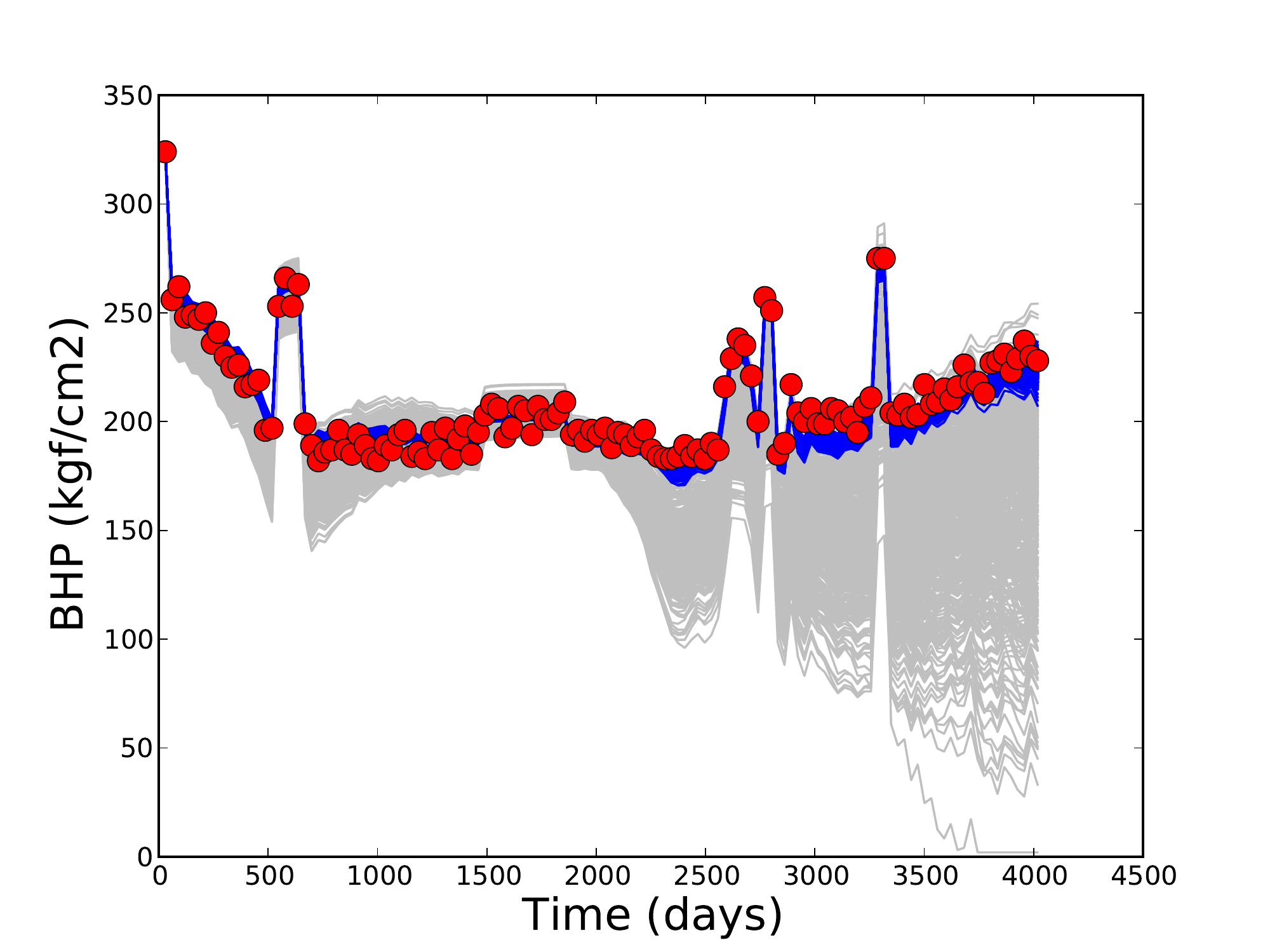}
    }
    \linebreak
    \subfloat[$8\times$-GEO1, Water cut]{
      \includegraphics[width=0.4\textwidth]{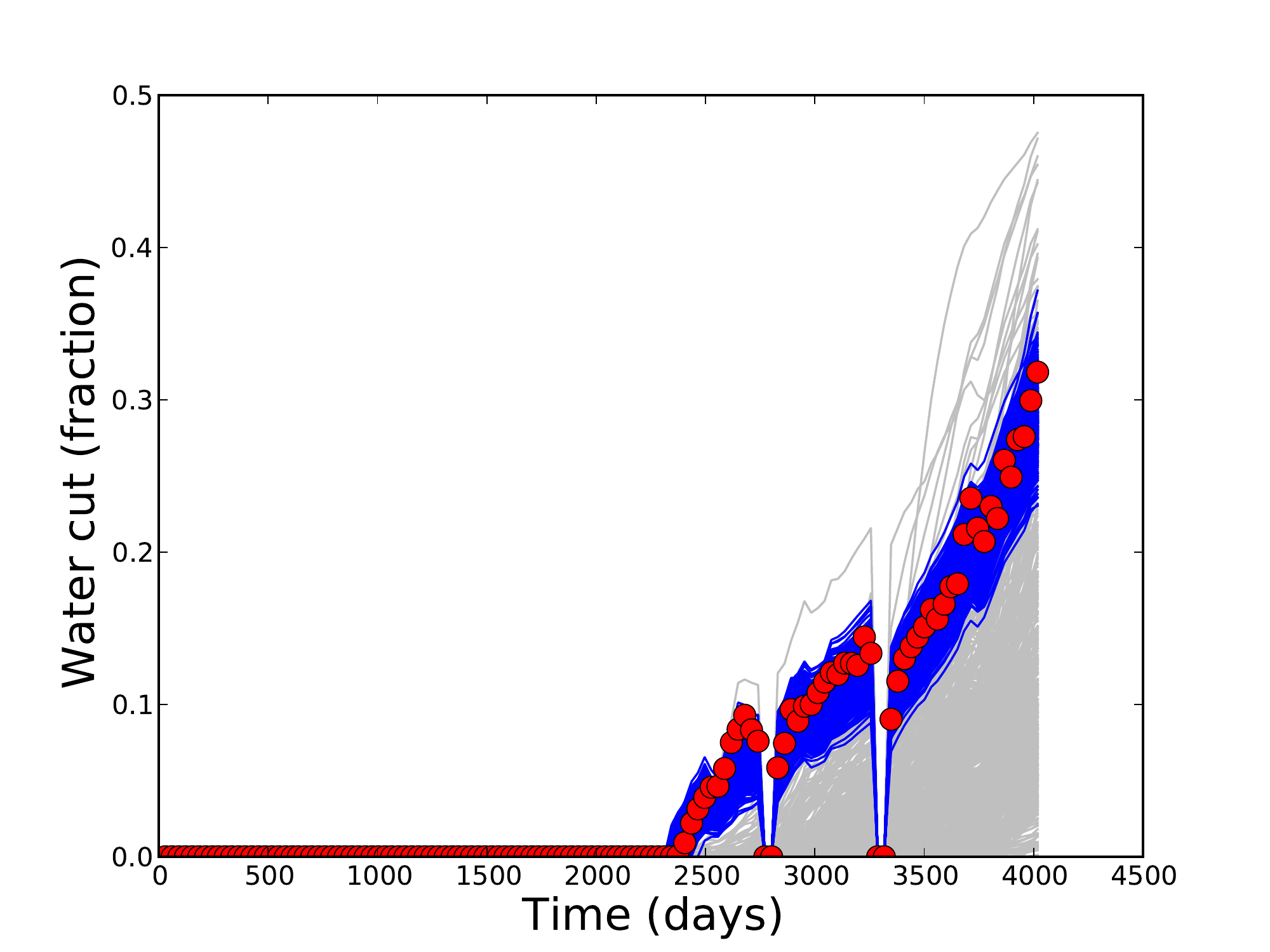}
    }
    \subfloat[$8\times$-GEO1, BHP]{
      \includegraphics[width=0.4\textwidth]{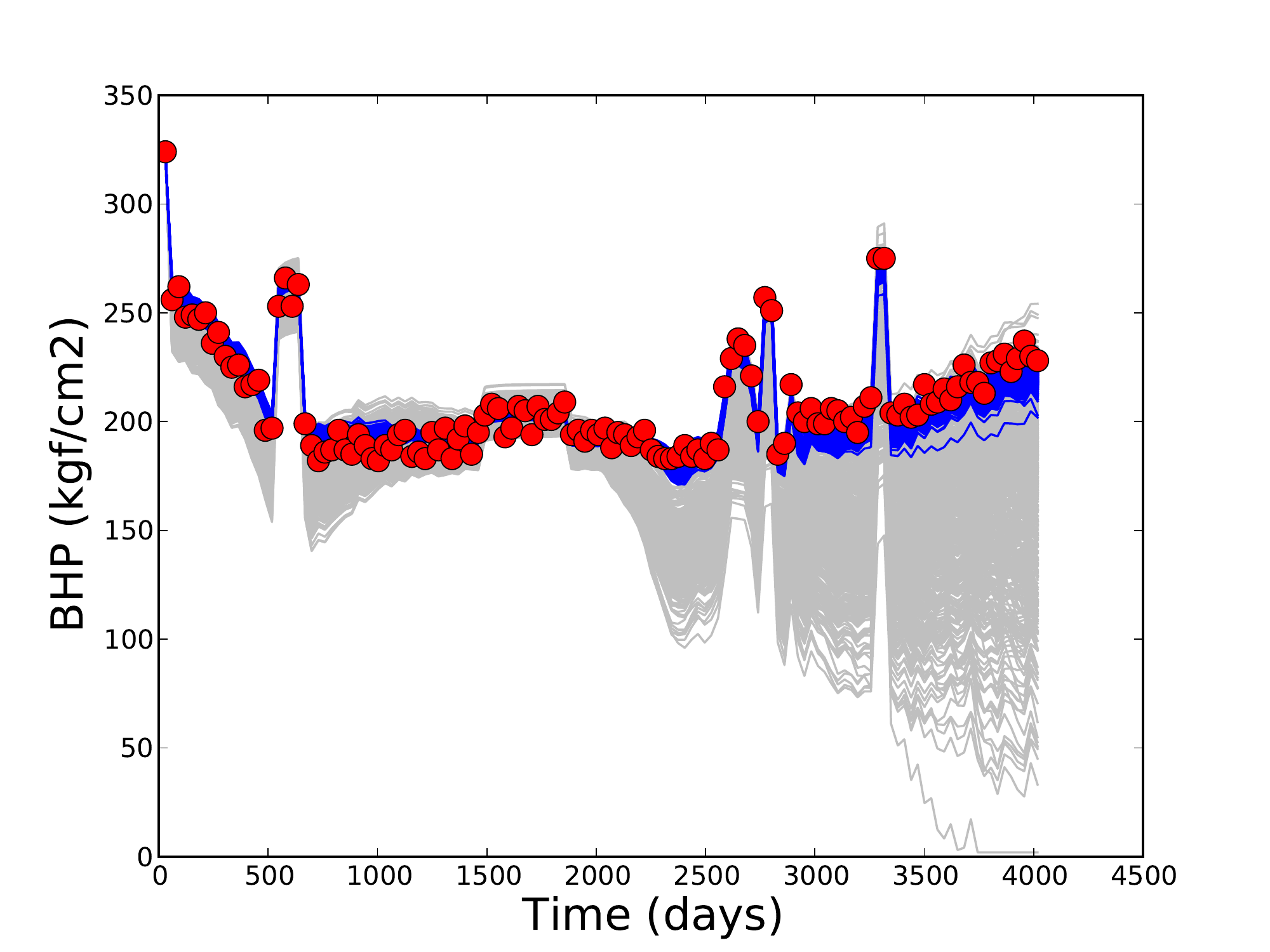}
    }
    \linebreak
    \subfloat[$8\times$-GEO2, Water cut]{
      \includegraphics[width=0.4\textwidth]{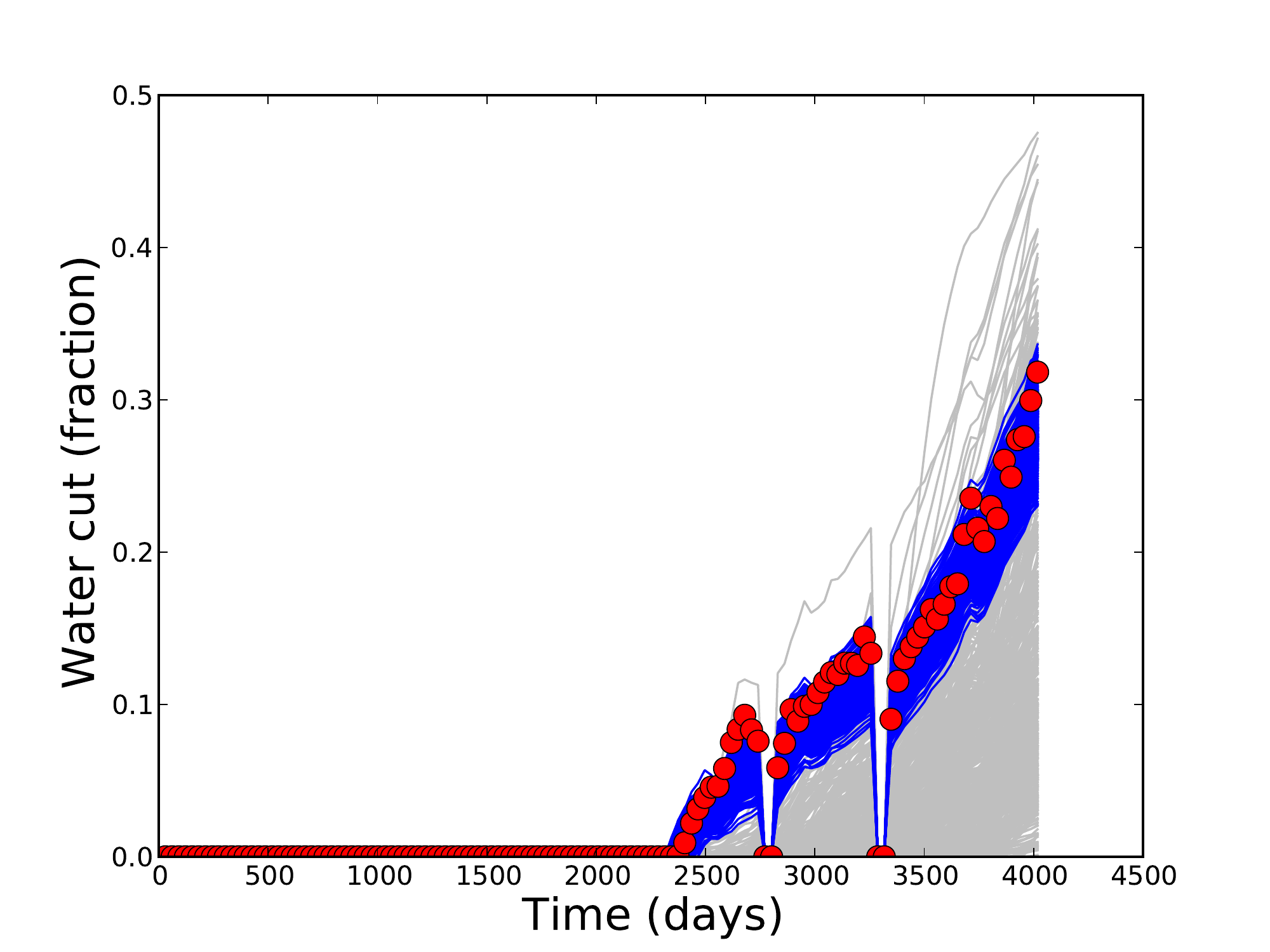}
    }
    \subfloat[$8\times$-GEO2, BHP]{
      \includegraphics[width=0.4\textwidth]{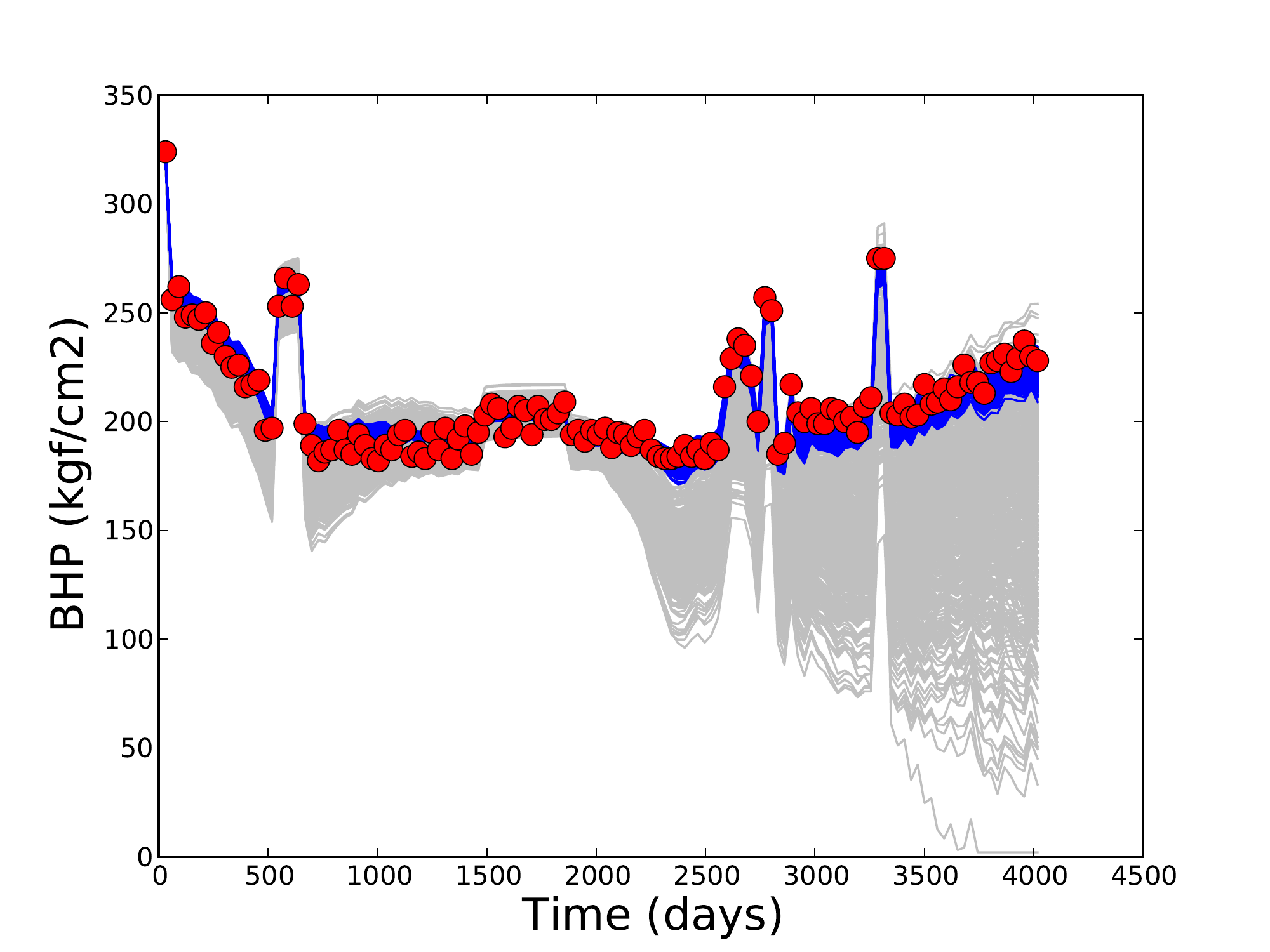}
    }
    \captionsetup{justification=justified}
\caption{Predicted production data for the well NA1A. Red dots are the observations, gray curves are the predictions from the prior and blue curves are the predictions from the posterior ensemble (UNISIM-I-H).}
\label{Fig:UNISIM-NA1A}
\end{figure}

Fig.~\ref{Fig:UNISIM-PERMI1} shows the ensemble mean of log-permeability for a intermediate layer for the cases $8\times$-CONST, $8\times$-GEO1 and $8\times$-GEO2. This figure shows that $8\times$-GEO1 and $8\times$-GEO2 resulted smoother fields with smaller changes in the prior realization. Fig.~\ref{Fig:UNISIM-NV} shows the normalized variance of log-permeability for the same three cases. The normalized variance is computed dividing the variance of the posterior ensemble by the variance of the prior ensemble. Clearly the methods $8\times$-GEO1 and $8\times$-GEO2 preserved more the variability in the posterior ensemble. Even though the correct level of variance is unknown, previous experience with ensemble data assimilation methods show a tendency of variance underestimation even when localization is used; see, e.g., \citep{emerick:16b}. Therefore, the larger normalized variances of the cases $8\times$-GEO1 and $8\times$-GEO2 are a very positive result, especially if the posterior ensemble is to be used for uncertainty quantification of production forecast and to assimilated future data. Fig.~\ref{Fig:UNISIM-NA1A} shows the predicted data for a selected well of the field. This figure illustrates the fact that even though $8\times$-GEO1 and $8\times$-GEO2 resulted in larger data-mismatch norm than $8\times$-CONST, the three procedures results in predicted data with a reasonable agreement with the observations.

\subsection{Field Case}
\label{Sec:FieldCase}

The last test case corresponds to the same field problem presented in \citep{emerick:16a}. This is heavy-oil reservoir in Campos Basis operating with 36 wells. The main recovery mechanism is waterflooding. The model contains 130,000 active gridblocks and the prior ensemble includes 200 realizations of porosity and horizontal permeability. These realizations are the same used in \citep{emerick:16a}. Besides porosity and permeability, the history matching parameters include the transmissibility across several faults, the maximum water relative permeability curve and the Corey exponent of the water-oil relative permeability. The production history includes monthly measurements of water cut, gas-oil ratio (GOR) and BHP. Compared to the original study presented in \citep{emerick:16a}, the production history has been updated with five additional years. The noise in the data was assumed Gaussian with zero mean and standard deviation of 10\% of the water cut, 20\% of the GOR and a constant value of 5~kgf/cm$^2$ for pressure measurements. During data assimilations, Schur-product localization was applied to the Kalman gain using the Gaspari-Cohn correlation function with a constant correlation length of 2000 meters. Similarly to the other cases, subspace inversion was applied retaining 99\% of the singular values.

Unlike the UNISIM-I-H case, for this field problem, the discrepancy function (\ref{Eq:f2}) was greater than zero any $\alpha_1 \geq N_a = 4$. For comparisons, the cases with constant inflations and the method proposed in \citep{rafiee:17a} were also performed using four data assimilations. Table~\ref{Tab:Field-Alphas} presents the inflation factors for each case. Note that $4\times$-GEO1 started with relatively large inflation ($\alpha_1 = 16987$), as a result the inflation at the last data assimilation is only 1.04. For $4\times$-GEO2, on the other hand, because the geometric sequence is selected such that at the last data assimilation $\alpha_{N_a} = 1.5$, the initial inflation was only $\alpha_1 =  37.33$.

\begin{table}
\caption{Inflation factors (Field Case)}
\label{Tab:Field-Alphas}
\begin{center}
\begin{tabular}{lccc}
\toprule
$k$ & $4\times$-CONST & $4\times$-GEO1 & $4\times$-GEO2  \\
\midrule
1 & 4 & 16986.84 & 37.33 \\
2 & 4 & 670.47 & 12.79 \\
3 & 4 & 26.46 & 4.38 \\
4 & 4 & 1.04 & 1.50 \\
\midrule
$\gamma$ & 1 & 0.0395& 0.3425 \\
\bottomrule
\end{tabular}
\end{center}
\end{table}

\begin{figure}
	\centering
	\includegraphics[width=0.55\linewidth]{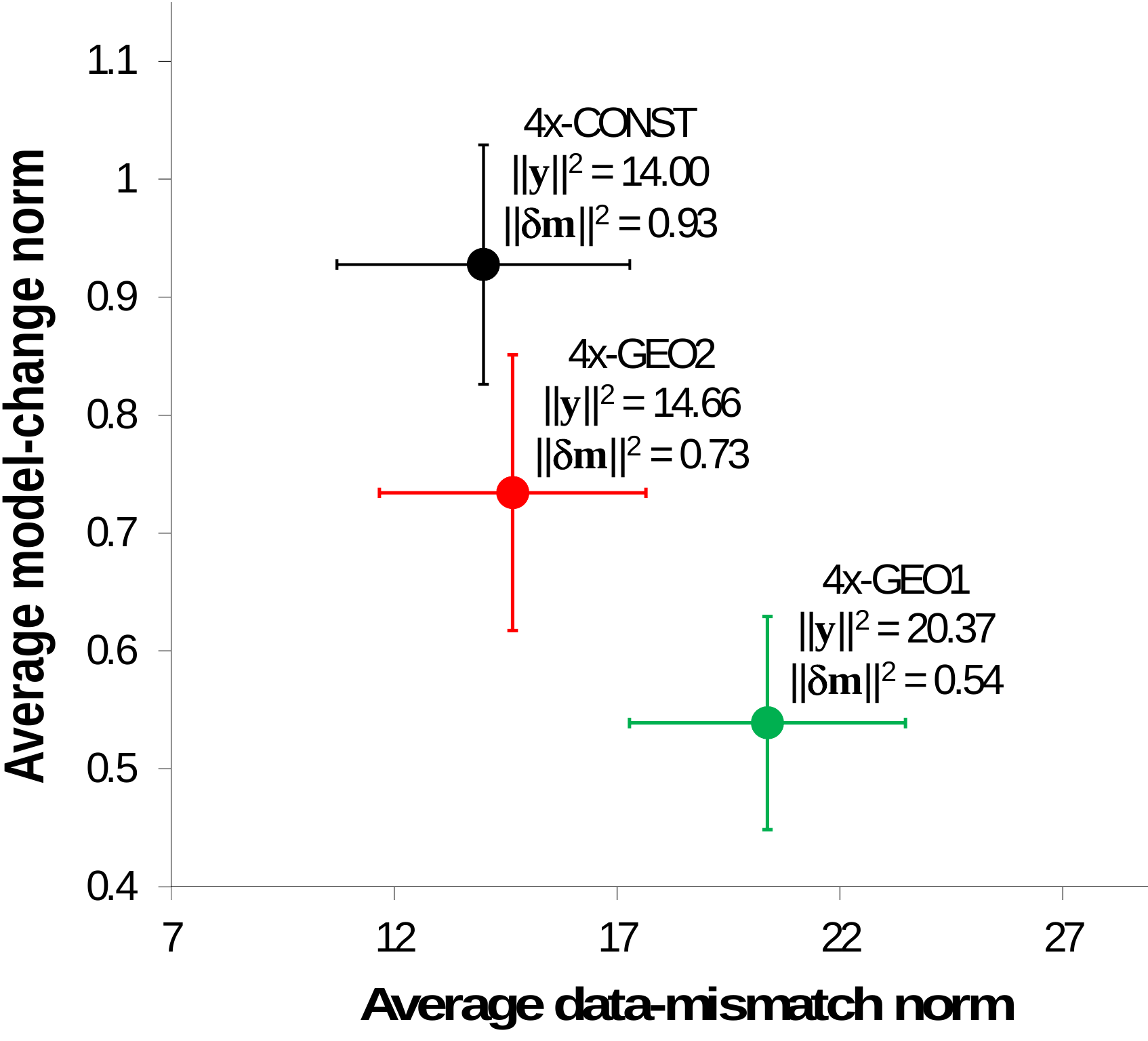}
	\caption{Average data-mismatch and model-change norms (Field Case). The error bars correspond to one standard deviation.}
	\label{Fig:Field-Norms}
\end{figure}

\begin{figure}
\centering
    \captionsetup{justification=centering}
    \subfloat[Prior]{
      \includegraphics[width=0.5\textwidth]{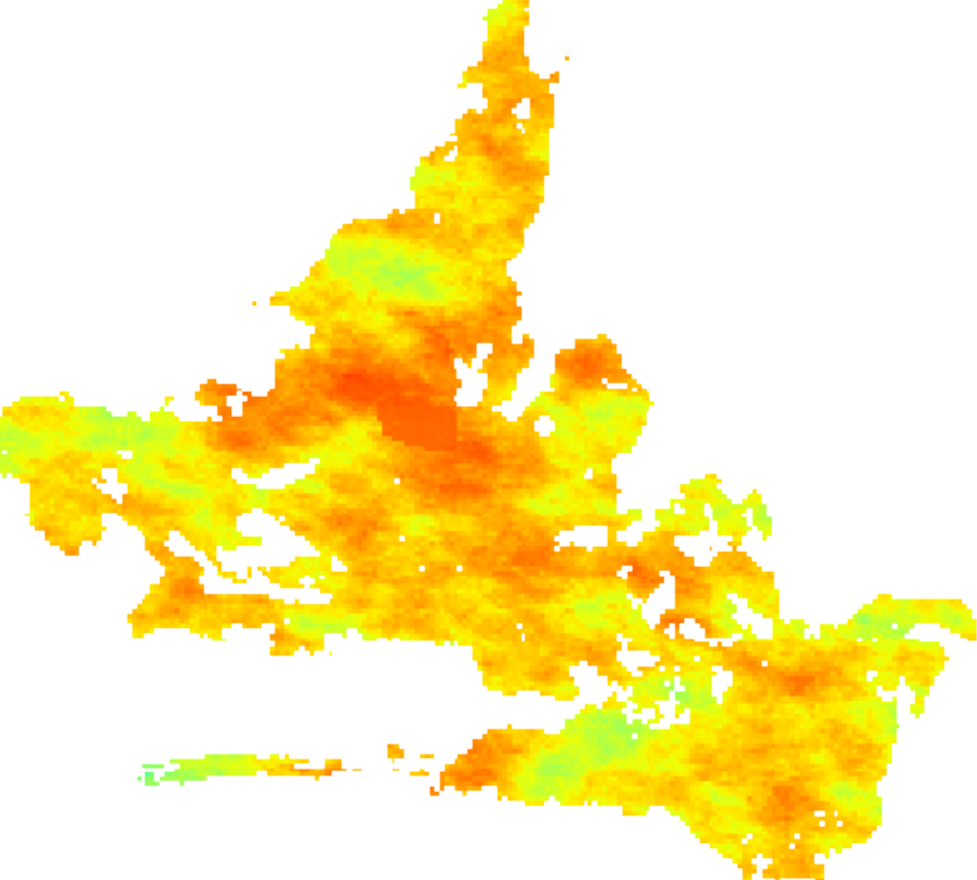}
    }
    \subfloat[$4\times$-CONST]{
      \includegraphics[width=0.5\textwidth]{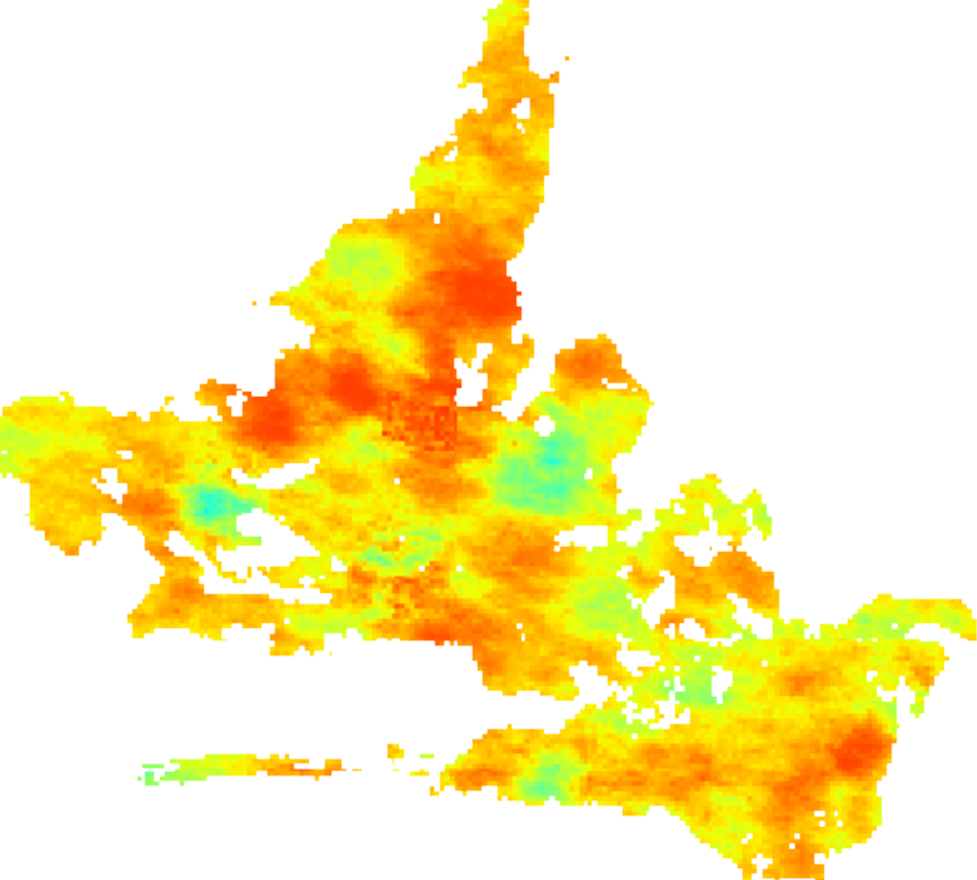}
    }
    \linebreak
    \subfloat[$4\times$-GEO1]{
      \includegraphics[width=0.5\textwidth]{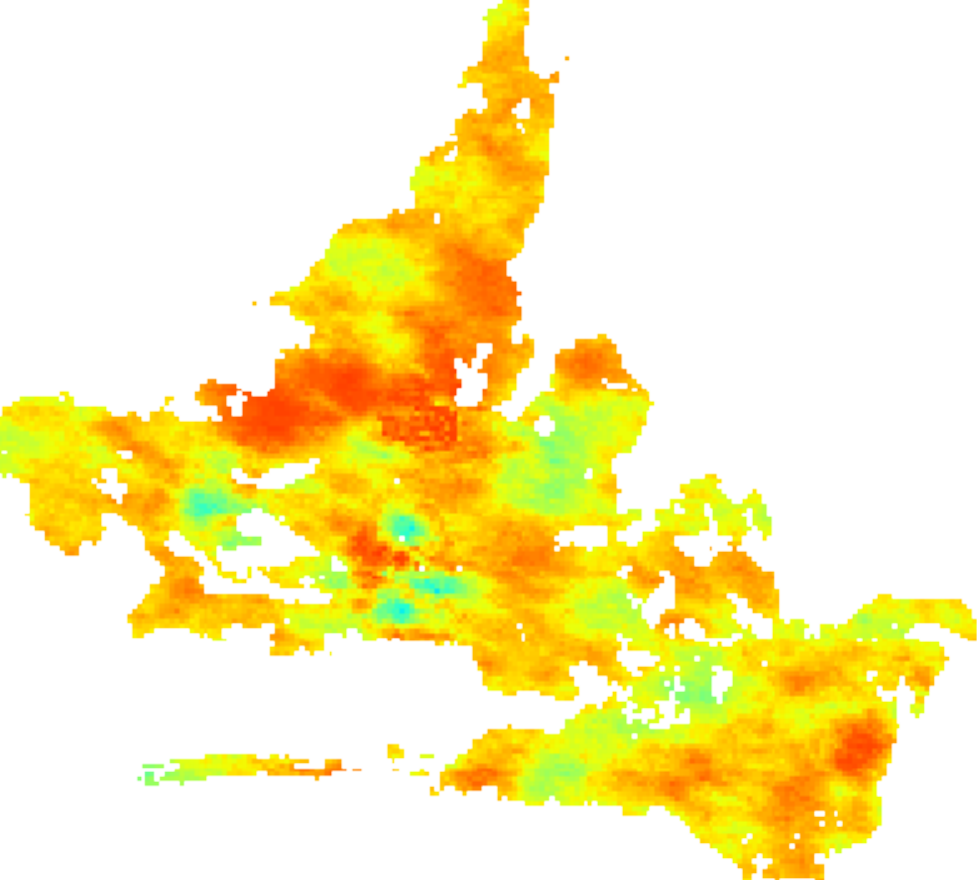}
    }
    \subfloat[$4\times$-GEO2]{
      \includegraphics[width=0.5\textwidth]{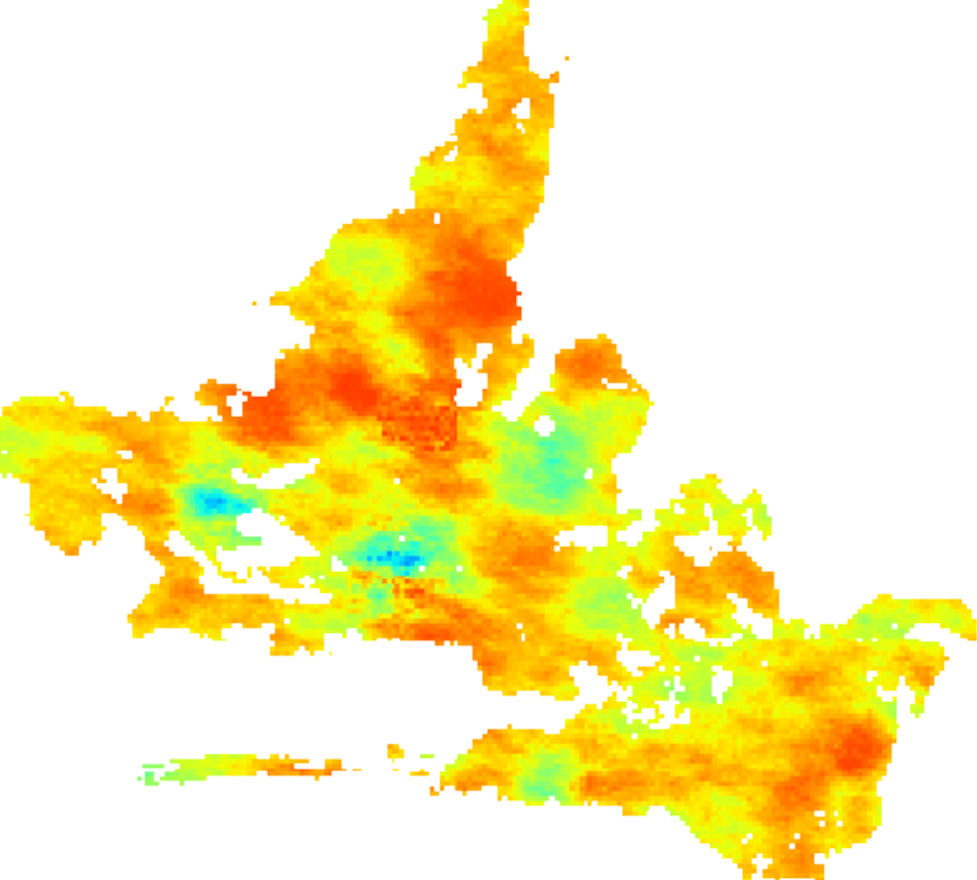}
    }
    \linebreak
      \includegraphics[width=0.5\textwidth]{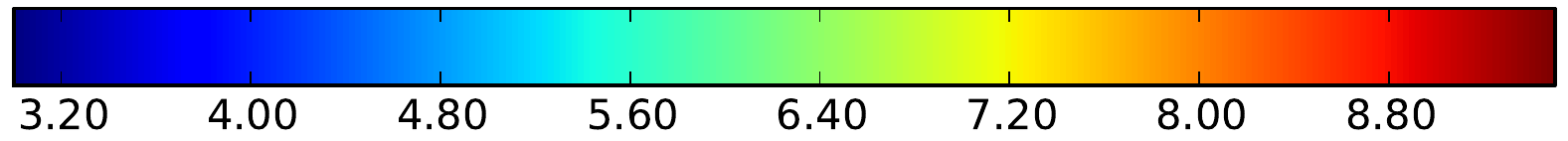}
    \captionsetup{justification=justified}
\caption{Mean log-permeability (Field Case, layer 17).}
\label{Fig:Field-PERMI1}
\end{figure}

\begin{figure}
\centering
    \captionsetup{justification=centering}
    \subfloat[$4\times$-CONST]{
      \includegraphics[width=0.5\textwidth]{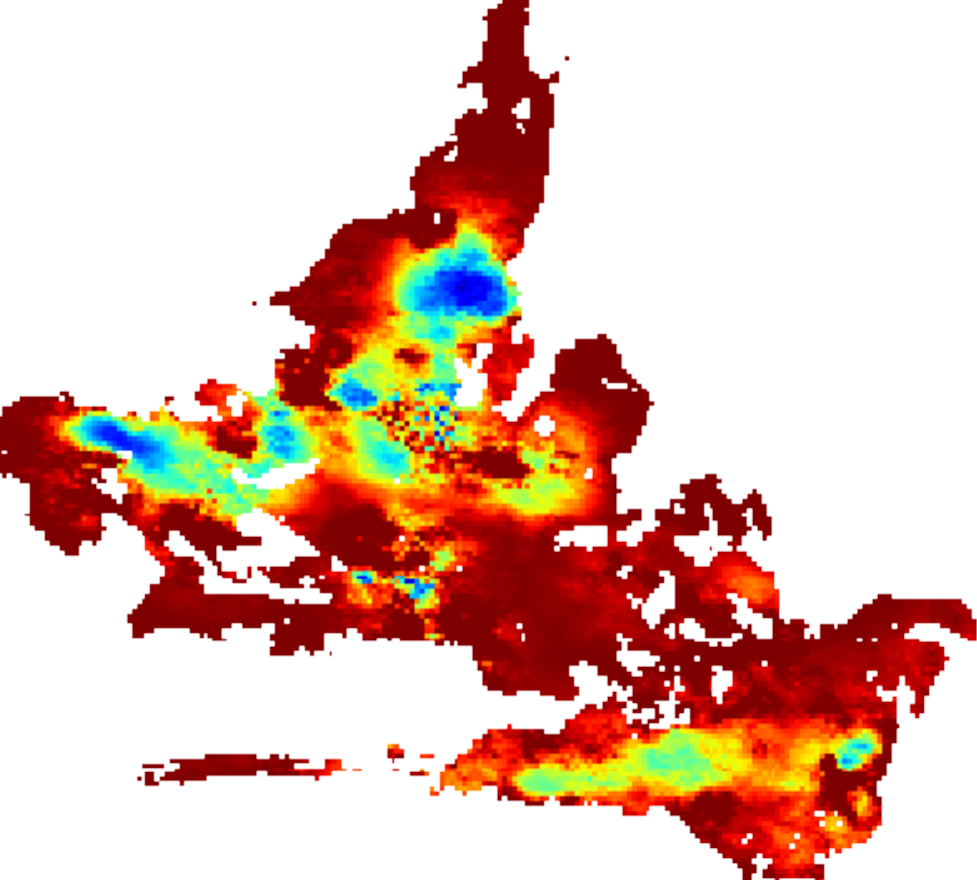}
    }
    \subfloat[$4\times$-GEO1]{
      \includegraphics[width=0.5\textwidth]{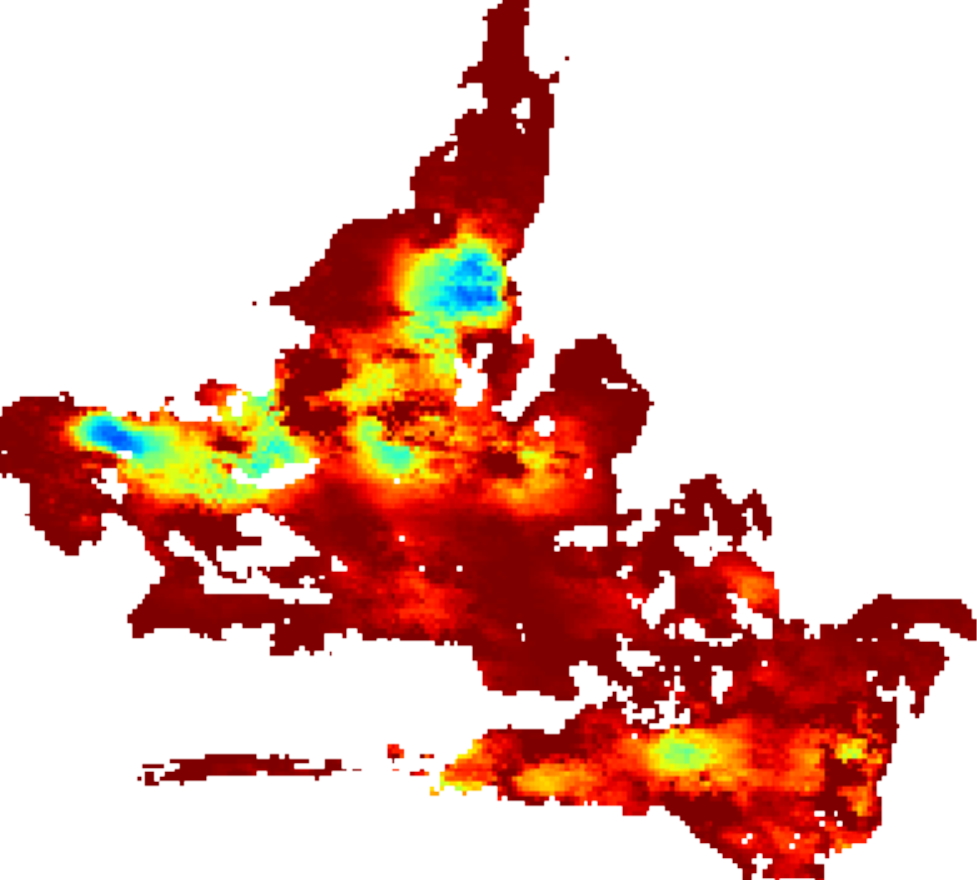}
    }
    \linebreak
    \subfloat[$4\times$-GEO2]{
      \includegraphics[width=0.5\textwidth]{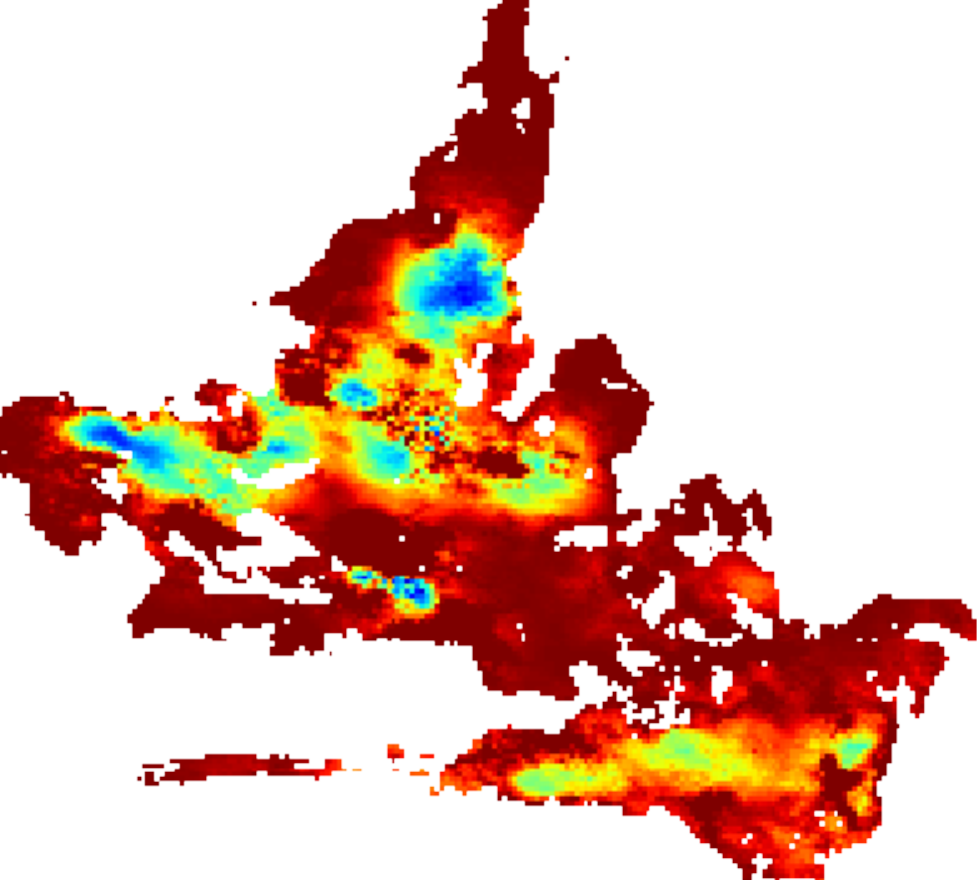}
    }
    \linebreak
      \includegraphics[width=0.5\textwidth]{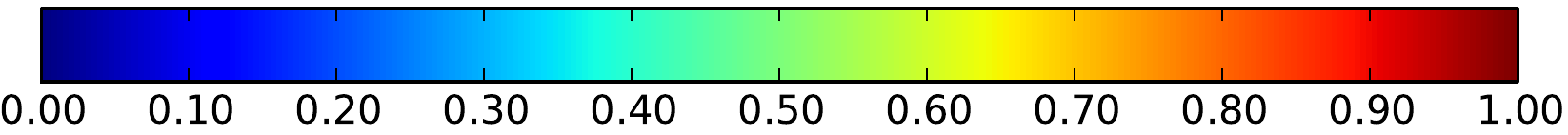}
    \captionsetup{justification=justified}
\caption{Normalized variance of log-permeability for (Field Case, layer 17).}
\label{Fig:Field-NV}
\end{figure}

\begin{figure}
\centering
    \captionsetup{justification=centering}
    \subfloat[$4\times$-CONST, Water cut]{
      \includegraphics[width=0.33\textwidth]{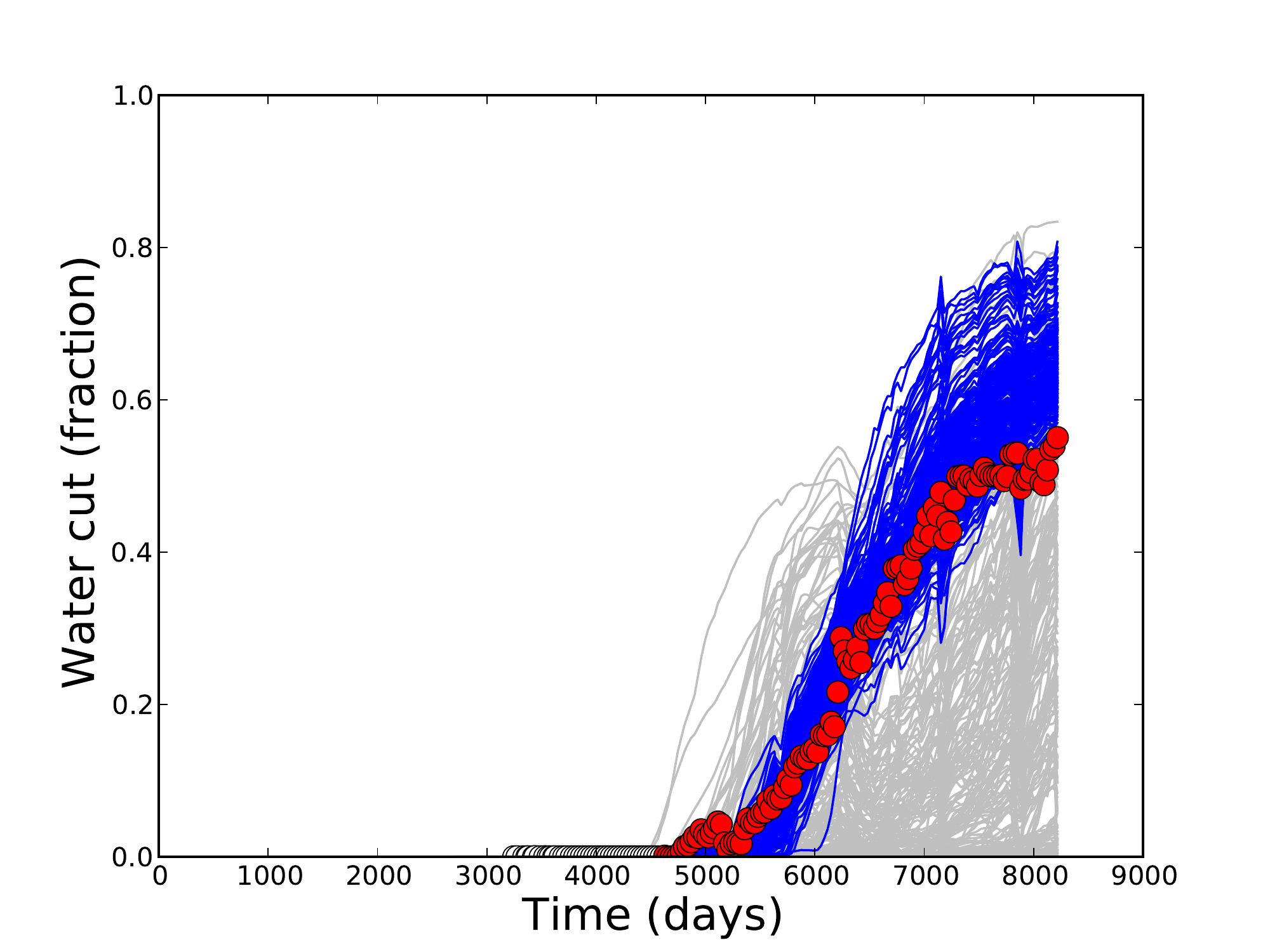}
    }
    \subfloat[$4\times$-CONST, GOR]{
      \includegraphics[width=0.33\textwidth]{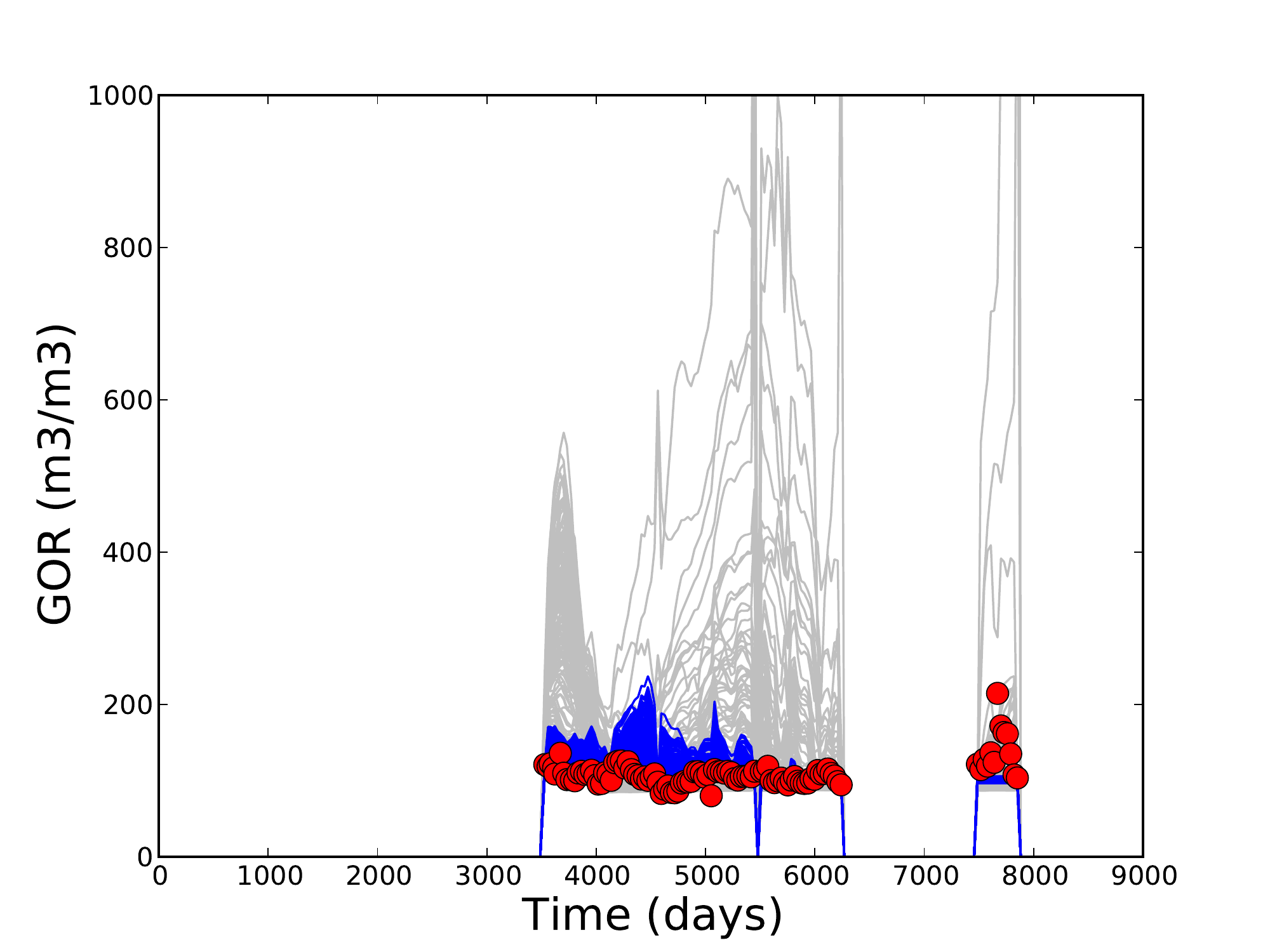}
    }
    \subfloat[$4\times$-CONST, BHP]{
      \includegraphics[width=0.33\textwidth]{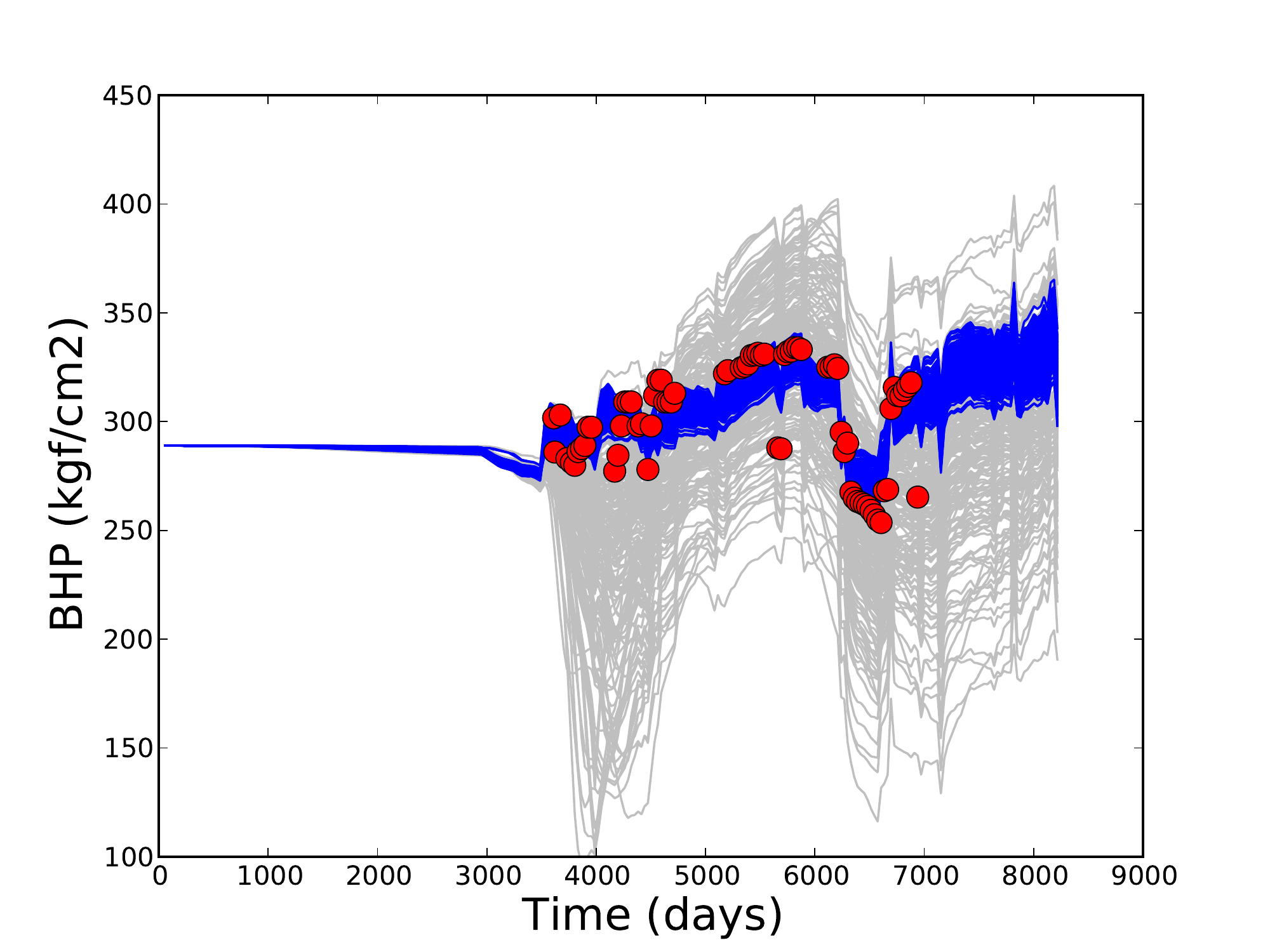}
    }
    \linebreak
    \subfloat[$4\times$-GEO1, Water cut]{
      \includegraphics[width=0.33\textwidth]{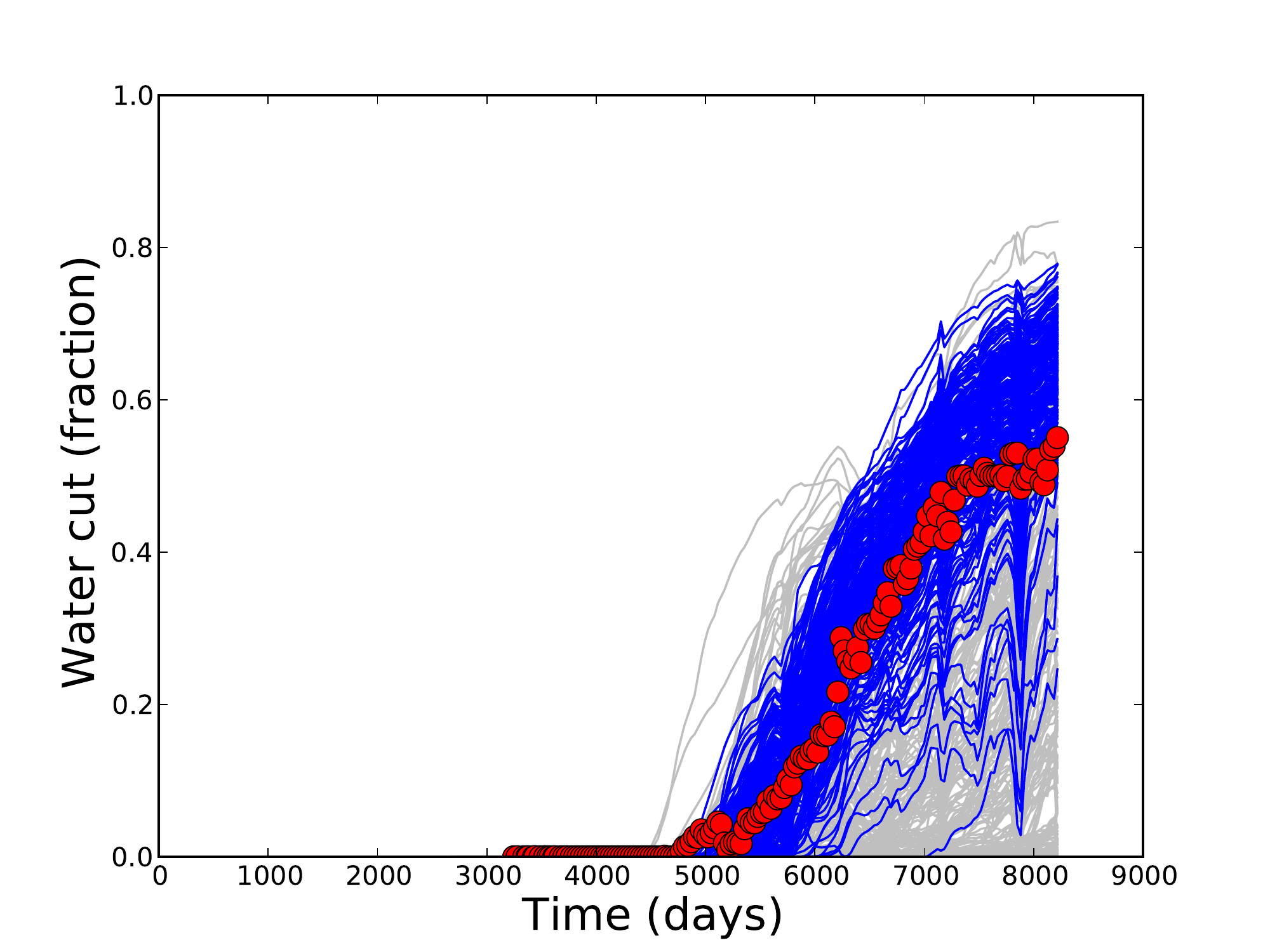}
    }
    \subfloat[$4\times$-GEO1, GOR]{
      \includegraphics[width=0.33\textwidth]{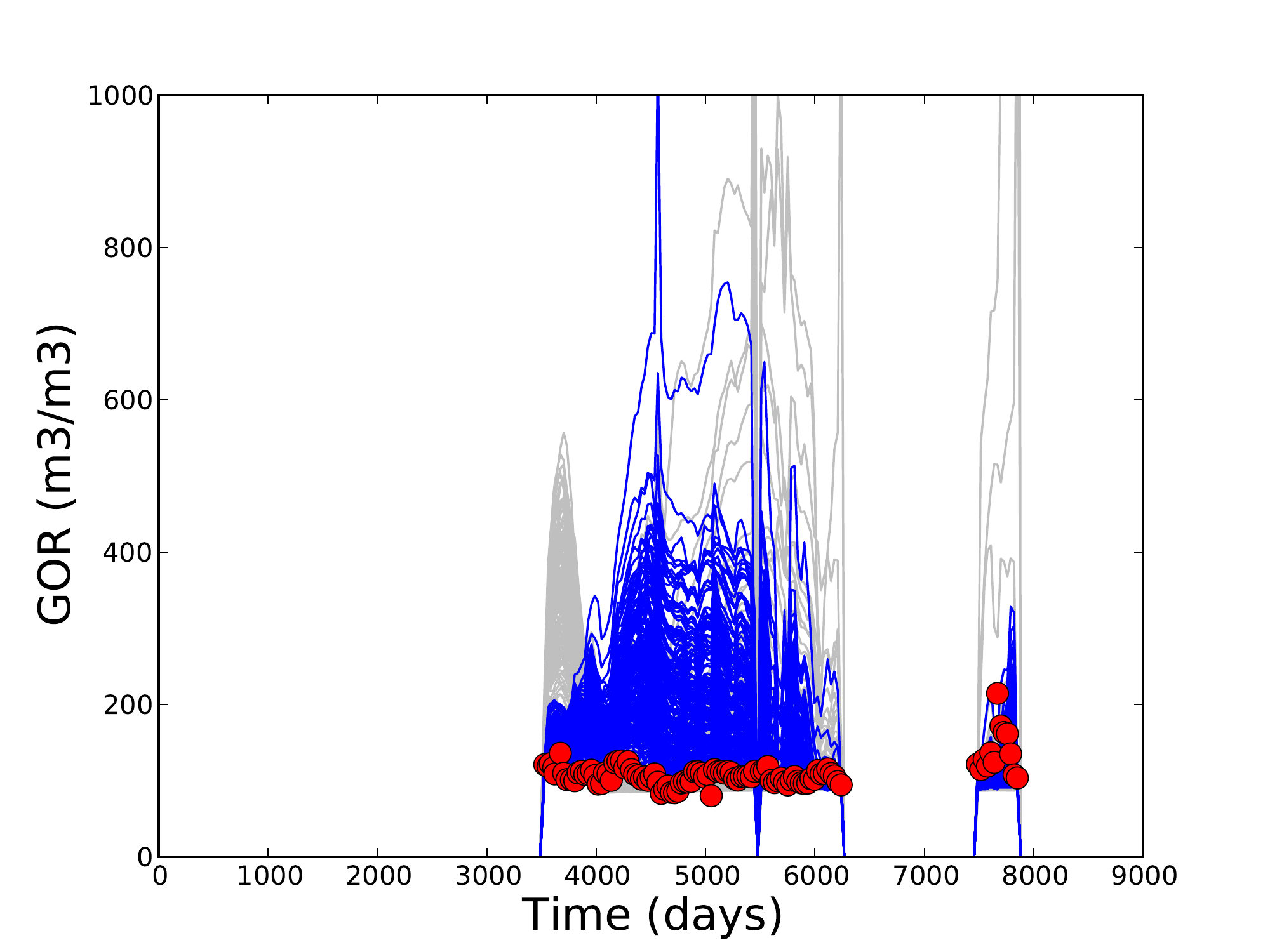}
    }
    \subfloat[$4\times$-GEO1, BHP]{
      \includegraphics[width=0.33\textwidth]{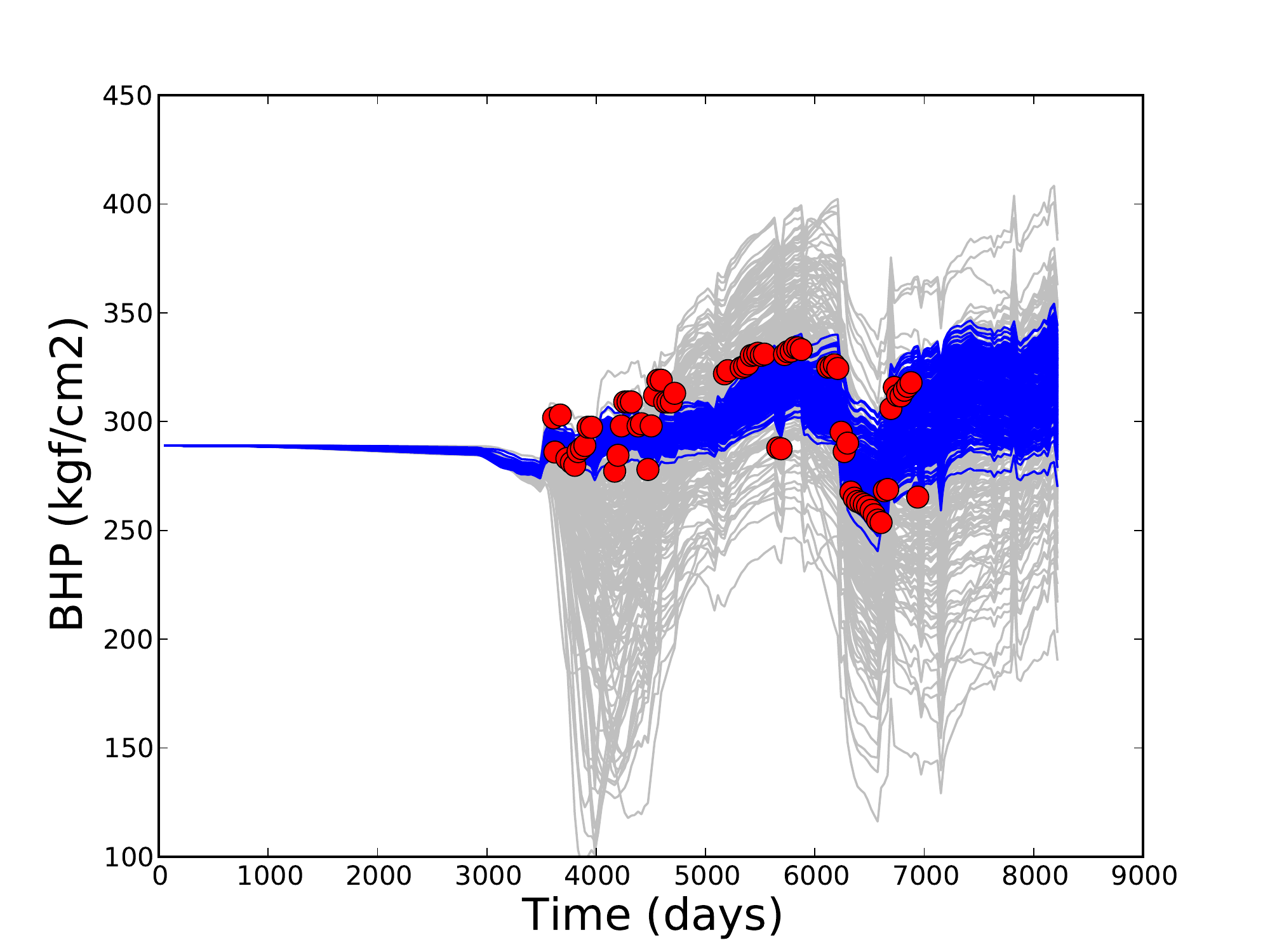}
    }
    \linebreak
    \subfloat[$4\times$-GEO2, Water cut]{
      \includegraphics[width=0.33\textwidth]{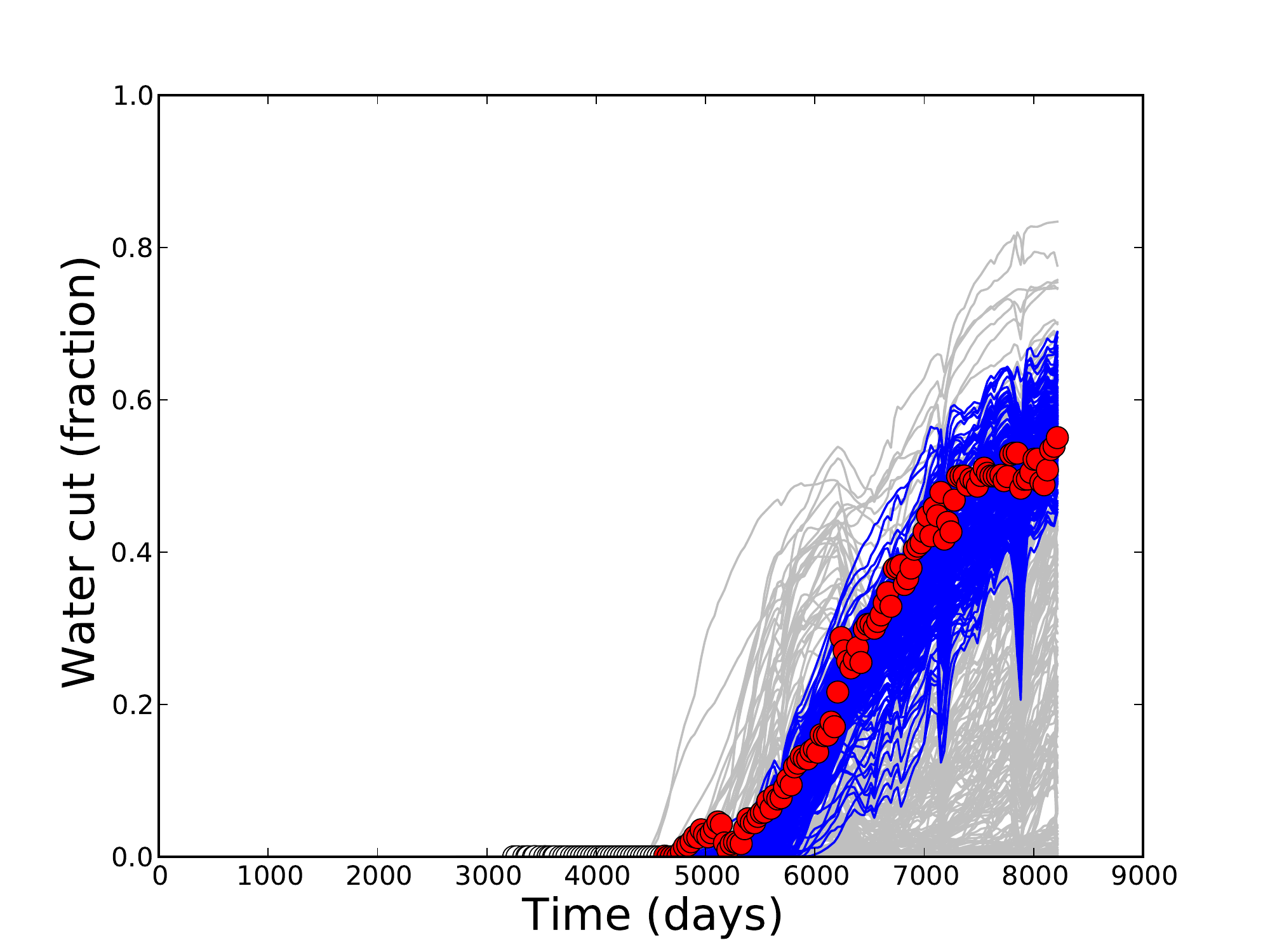}
    }
    \subfloat[$4\times$-GEO2, GOR]{
      \includegraphics[width=0.33\textwidth]{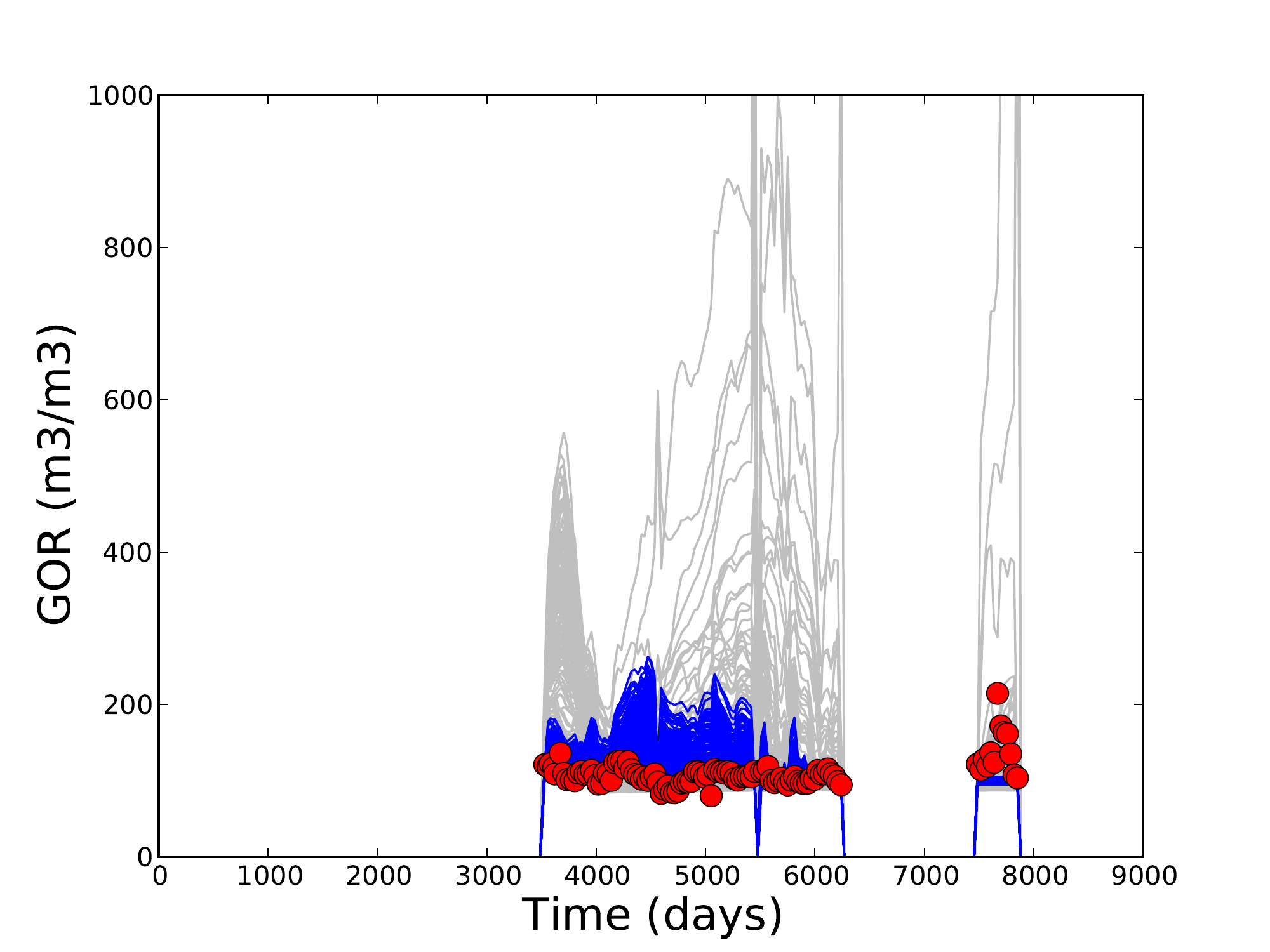}
    }
    \subfloat[$4\times$-GEO2, BHP]{
      \includegraphics[width=0.33\textwidth]{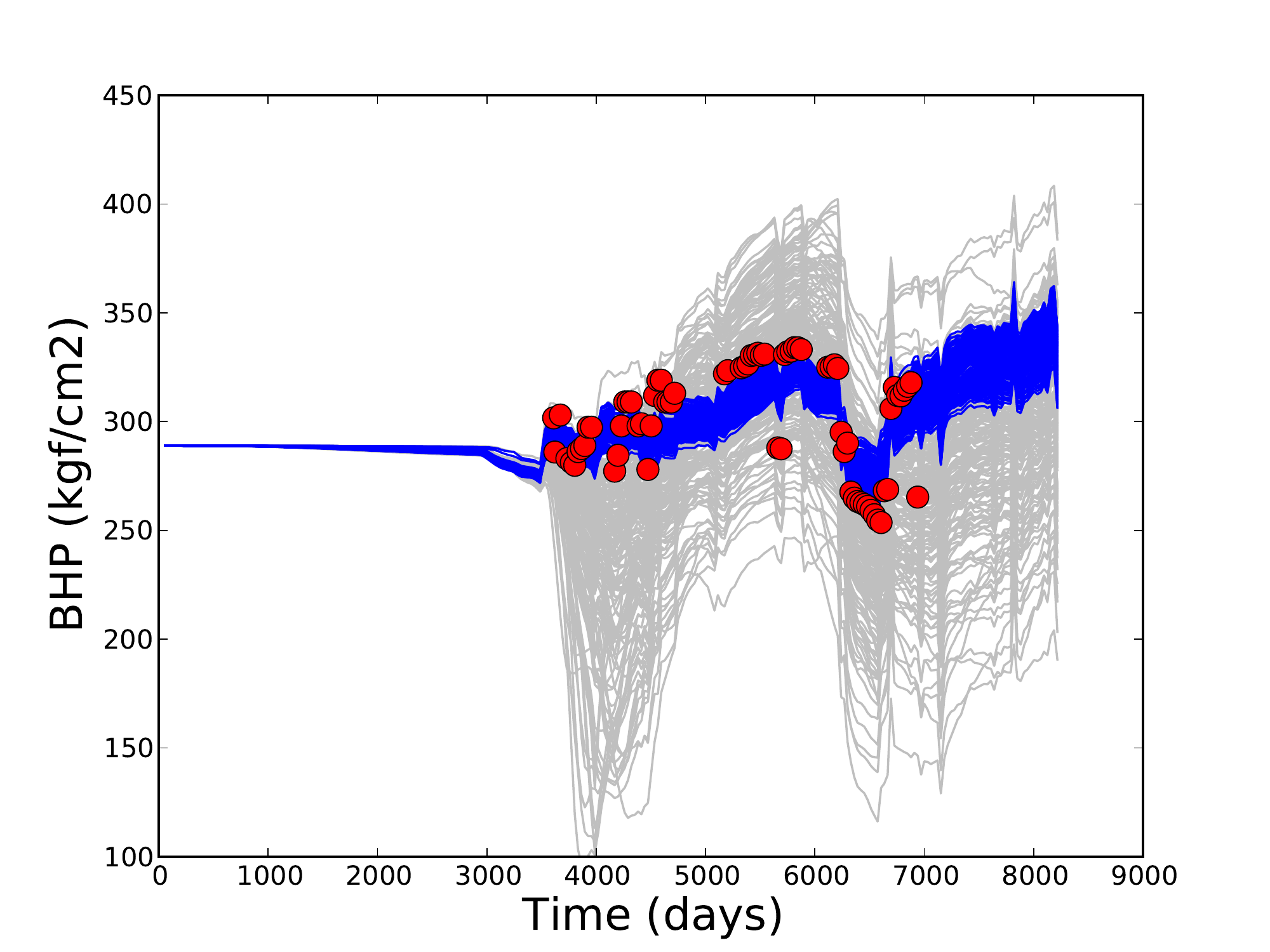}
    }
    \captionsetup{justification=justified}
\caption{Predicted production data for three wells. Red dots are the observations, gray curves are the predictions from the prior and blue curves are the predictions from the posterior ensemble (Field Case).}
\label{Fig:Field-Wells}
\end{figure}

Fig.~\ref{Fig:Field-Norms} shows the plot of the average data-mismatch and model-change norms for the three cases. The error bars in this figure correspond to one standard deviation for each side. The results in this figure shows that $4\times$-CONST resulted in the best data match and $4\times$-GEO1 resulted in the lowest model changes. The proposed approach obtained a compromise between both norms. Fig.~\ref{Fig:Field-PERMI1} shows the ensemble mean log-permeability indicating that all cases resulted in reasonable model estimates, which preserve the main characteristics of the prior realization. In terms of the posterior variance, the case $4\times$-GEO1 resulted in the highest values as illustrated in Fig.~\ref{Fig:Field-NV}. This results is consistent with the fact that this case resulted in the smallest changes in the prior realizations. However, the data matches after $4\times$-GEO1 are clearly of inferior quality as indicated in Fig.~\ref{Fig:Field-Wells}, which shows the predicted data for three different wells of the field.

\section{Conclusions}
\label{Sec:Conclusion}

This paper presented the results of an investigation on the use of ES-MDA with geometric inflation factors. The results warranted the following conclusions:
\begin{itemize}
  \item For the three test cases, the use of constant inflation factors resulted in the lowest values for the data-mismatch norm and the highest values for the model-change norm.
  \item The procedure proposed by \citet{rafiee:17a} obtained a good balance between data match and model change as longs as an appropriate number of data assimilation is select a priori. For the three test cases, this method selected a large value for the initial inflation. As a result, it is necessary to use a larger number of data assimilations to obtain acceptable results.
  \item The procedure introduced in this paper avoided this problem by selecting the geometric sequence based on the inflation at the last data assimilation. Overall, the initial inflation required to satisfy the MDP is significantly smaller than the initial inflation obtained by \citep{rafiee:17a}.
  \item The procedure introduced in this paper resulted in a trade-off between data match and model change for the three test problems.
  \item The use of geometric sequences also resulted in posterior ensembles with larger variances, which is typically desirable for uncertainty quantification, especially considering the fact that ensemble data assimilation methods tend to underestimate the posterior variance.
  \item Overall, the results are in agreement with the conjecture that starting with large inflations and gradually decreasing during the data assimilation steps improves the final model estimates.
\end{itemize}

\section*{Acknowledgement}
\label{Sec:Acknowledgement}

The author would like to thank Petrobras for supporting this research and for the permission to publish this paper.


\end{document}